\documentclass[11pt,a4paper]{article}
\usepackage[utf8]{inputenc}
\usepackage{amsmath,bm}
\usepackage{enumerate}
\usepackage{amsmath,nccmath}
\usepackage{amsthm}
\usepackage{amsfonts}
\usepackage[square, comma, sort&compress,numbers]{natbib}
\usepackage{amssymb}
\usepackage{setspace}
\usepackage{graphicx}
\usepackage{epstopdf}
\usepackage[shortlabels]{enumitem} 
\usepackage{verbatim}
\usepackage{color}
\usepackage[backref=page]{hyperref}
\hypersetup{
     colorlinks=true,    
     linktocpage=true,
    }
\renewcommand*{\backref}[1]{}
\renewcommand*{\backrefalt}[4]{~~{\tiny%
    \ifcase #1 Not cited.%
          \or [Cited on page~#2.]%
          \else [Cited on pages #2.]%
    \fi%
    }}
\usepackage[marginal]{footmisc}
\usepackage[section]{placeins}
\usepackage{indentfirst}
\usepackage[left=2cm,right=2cm,top=2cm,bottom=2cm]{geometry}
\author{}
\usepackage{float}
\setlist[itemize]{noitemsep,topsep=0pt}

\theoremstyle{plain}
\newtheorem{prop}{Proposition}[section]
\newtheorem*{prop*}{Proposition}
\newtheorem{cor}[prop]{Corollary}
\newtheorem{lem}[prop]{Lemma}
\newtheorem{thm}{Theorem}
\newtheorem{cormain}[thm]{Corollary}
\newtheorem*{lem*}{Lemma}
\newtheorem*{thm*}{Theorem}

\newtheorem{thms}[prop]{Theorem}

\theoremstyle{remark}

\theoremstyle{definition}
\newtheorem{rem}[prop]{Remark}
\newtheorem{defi}{Definition}[section]
\newtheorem*{defi*}{Definition}

\newtheorem*{ques*}{Question}
\newtheorem*{conj*}{Conjecture}
\newtheorem*{conven*}{Convention}

\renewcommand{\geq}{\geqslant}
\renewcommand{\leq}{\leqslant}

\definecolor{darkred}{rgb}{.56,0,0}

\allowdisplaybreaks
\usepackage{mathtools}

\usepackage{authblk}
\usepackage{titlesec} 

\numberwithin{equation}{section}
\titleformat{\subsubsection}[runin]
  {\normalfont\normalsize\bfseries}{\thesubsubsection}{1em}{}	
\makeatletter
\newcommand*{\rom}[1]{\expandafter\@slowromancap\romannumeral #1@}

\makeatother
\title{{\bf Robust heterodimensional dynamics in two-parameter unfoldings of homoclinic tangencies}}

\author{Dongchen Li\footnote{dongchenli@fudan.edu.cn, Shanghai Center for Mathematical Sciences, Fudan University, China}, 
Xiaolong Li\footnote{lixl@hust.edu.cn, School of Mathematics and Statistics, Huazhong University of Science and Technology, China}, 
Katsutoshi Shinohara\footnote{ka.shinohara@r.hit-u.ac.jp, Graduate School of Business Administration, Hitotsubashi University, Japan},  
Dmitry Turaev\footnote{d.turaev@imperial.ac.uk, Department of Mathematics, Imperial College London, United Kingdom}
}

\usepackage{tocloft}
\begin{document}

\def\D{\mathrm{D}}
\def\rank{\mathrm{rank}\,}
\def\d{\mathrm{d}}
\def\tr{\mathrm{tr}\,}
\def\diff{Di\!f\!f}
\def\eps{\varepsilon}
\def\Poincare{Poincaré\,\,}

\bibliographystyle{plainnat}

\maketitle

\noindent{\bf Abstract.}
We establish a necessary and  sufficient   condition for the birth of heterodimensional cycles from a generic homoclinic tangency to a hyperbolic periodic orbit. We prove for $C^r$ ($r=3,\dots,\infty,\omega$) dynamical systems on a manifold $\mathcal{M}$, with $\dim \mathcal{M}\geqslant 3$ for diffeomorphisms and with $\dim \mathcal{M}\geqslant 4$ for flows, that $C^1$-robust heterodimensional dynamics of coindex one appear in any generic two-parameter $C^r$ unfolding of a homoclinic tangency to a periodic orbit such that  at least one central multiplier is not real and the central dynamics are not sectionally dissipative. The heterodimensional dynamics also involve a blender exhibiting $C^1$-robust homoclinic tangencies. As a corollary, any system with a homoclinic tangency of the class described above belongs to the $C^r$ closure of the $C^1$-open Newhouse domain.




\tableofcontents

\section{Introduction}
Homoclinic tangencies and heterodimensional cycles  are the most fundamental objects in non-hyperbolic dynamics.
They are fragile, but also persistent. That is,  they can be destroyed by an arbitrarily small perturbation, but whenever a given one disappears due to a perturbation, another one can immediately emerge. Examples of such phenomenon were discovered by Newhouse \cite{Ne70} for homoclinic tangencies, and by Abraham and Smale \cite{AS70} and Simon \cite{Si72} for heterodimensional cycles. This persistence property implies the existence of open regions in the space of smooth dynamical systems where
every system is non-hyperbolic and structurally unstable (i.e., arbitrarily close to any system from such region there exist
systems which are not topologically equivalent to it).

Remarkably, it is proven that the persistence property in fact holds for all typical homoclinic tangencies and coindex-1 heterodimensional cycles. Namely, every system with a homoclinic tangency is in the $C^r$ ($r=2,\dots,\infty,\omega$) closure of the  $C^2$-open {\em Newhouse domain}, i.e., the  region  where
systems with homoclinic tangencies are $C^r$-dense (see \citep{Ne79,GTS93,PV94,Ro95}), and every system with a coindex-1 heterodimensional cycle is in the $C^r$ ($r=1,\dots,\infty,\omega$) closure of the $C^1$-open {\em Bonatti-D\'iaz domain}, i.e., the  region where systems with heterodimensional cycles are $C^r$-dense (see \citep{BD08,BDK12,LT21}). It was conjectured by Palis \citep{Pa00} that, along with the {\em Lorenz domain} where smooth flows with homoclinic loops to an equilibrium state are dense, these regions cover all robustly non-hyperbolic systems.


In this paper, we focus on homoclinic tangencies whose local dynamics are not ``essentially two-dimensional'' (see Definition \ref{defi:d_eff}). 
We prove (Theorem \ref{thm:eff}) that if a system has such homoclinic tangency,
 then it is  in the $C^r$ ($r=1,\dots,\infty,\omega$) closure of the  Bonatti-D\'iaz domain. Besides establishing a relation between the Newhouse and Bonatti-D\'iaz domains, the importance of our result is based on the finding in \citep{LT21} that, $C^r$-densely, the system in Bonatti-D\'iaz domain has an infinite sequence\footnote{For a detailed study of the blender sequences near heterodimensional cycles, we refer the readers to \citep{Li24a} for their generation mechanisms and to \citep{Li24b} for their structures and properties.} of {\em blenders} accumulating on a heterodimensional cycle, which are special hyperbolic sets responsible for various non-trivial robust phenomena that are fragile at  first glance (see the discussion in Section \ref{sec:intro3}). Thus, our result allows for possible extensions of blender-based results to higher dimensional systems having a homoclinic tangency.

In particular, as shown in \citep{Li24b},  the mentioned blender sequences contain  a generalization of the {\em blender-horseshoes} introduced  by Bonatti and D\'iaz in \citep{BD12}. Hence, based on the results in \citep{BD12,Li24b}, we obtain (Corollary \ref{thm:C1tangency}) that every system of Theorem \ref{thm:eff}  also lies in the $C^r$ closure of the {\em $C^1$-open} Newhouse domain (augmenting  the above-mentioned general result that the $C^2$-open Newhouse domain is $C^r$-close to every system with a homoclinic tangency). 


\subsection{Effective dimension and robust heterodimensional cycles}\label{sec:intro1}
A {\em homoclinic tangency} is a non-transverse intersection of the stable and unstable manifolds of a hyperbolic periodic orbit.  Systems from the Newhouse domain exhibit the behavior very different from that of the uniformly hyperbolic systems, such as the presence of
zero Lyapunov exponents \citep{DG09},
super-exponential growth of the number of periodic points \citep{Ka00},
divergence of Birkhoff averages \citep{KS17,KNS21,Ba21,BB21},
coexistence of infinitely many sinks \citep{Ne74,GS86,GST08,Be16} or other attractors or repellers \cite{GTS93,GTS97,GST08} and, in general, the dynamics of ultimate richness and complexity \cite{Tu10,Tu15}. Similar dynamical phenomena also take place in the Bonatti-D\'iaz domain (see e.g.  \citep{GIKN05,Be16,AST17,AST21}).

In \citep{GTS93,GST08}, the coexistence of saddle  periodic orbits with different dimensions of unstable manifolds was highlighted as the most basic property of systems from the Newhouse domain, provided that the dynamics have enough ``effective dimension'' for these orbits to exist, i.e., the dynamics are not reducible to those of
a two-dimensional map.
For a generic homoclinic tangency to a periodic orbit, the effective dimension $d_{\rm eff}$ is determined by the {\em central multipliers} of the periodic orbit. Multipliers are the eigenvalues of the linearization of the return map near the periodic orbit. The central multipliers are those which are the closest to the unit circle from outside and those closest to the unit circle from inside. In general position, the central multipliers are simple and there are either two real central multipliers -- one less than $1$ in the absolute value and the other larger than $1$ in the absolute value, or one real and one pair of complex-conjugate, or two pairs of complex-conjugate eigenvalues. According to \citep{GST08}, under some genericity conditions on the orbit of the homoclinic tangency, we have  
$d_{\rm eff}>1$ if there are two different central multipliers such that the absolute value of their product is larger than 1 and also two different central multipliers for which the absolute value of their product is smaller than $1$, i.e., the return map near the periodic orbit contracts some two-dimensional areas and also expands some. 
One of the main results of \cite{GST08} is that when
$d_{\rm eff}>1$,  {\em a generic unfolding of the homoclinic tangency creates coexisting hyperbolic periodic orbits with different dimensions of their unstable manifolds.}

In the present paper we significantly strengthen this result by showing that unfolding such a tangency  creates ``robustly connected'' hyperbolic sets with different indices. Recall that the index of a transitive hyperbolic set is the rank of its unstable bundle (i.e., the number of the positive Lyapunov exponents). We use the term {\em system} to refer to either a diffeomorphism of a manifold of dimension at least three, or a flow on  a manifold of dimension at least four.

\begin{defi}[Heterodimensional cycles]\label{def:hdc}
We say that a  system $f$ has a heterodimensional cycle involving two hyperbolic periodic orbits $L_1$ and $L_2$ of different indices if $W^u(L_1)\cap W^s(L_2)\neq \emptyset$ and  $W^u(L_2)\cap W^s(L_1)\neq \emptyset$.
\end{defi}

\begin{defi}[Robust heterodimensional dynamics]\label{def:hdd}
We say that a  system $f$ has a {\em $C^1$-robust heterodimensional dynamics}\footnote{What we call (robust) heterodimensional dynamics are often called (robust) heterodimensional cycles. To facilitate the discussion, we reserve the term ``heterodimensional cycle'' for the most basic case where the two hyperbolic sets are trivial, i.e., they are two hyperbolic periodic orbits, see Definition \ref{def:hdc}.}
 if there are two uniformly-hyperbolic basic (i.e., compact, transitive and locally maximal) sets $\Lambda_1(f)$ and $\Lambda_2(f)$ with different indices, and a $C^1$-neighborhood $\mathcal{U}$ of $f$ such that $W^u(\Lambda_{1}(\tilde f))\cap W^s(\Lambda_{2}(\tilde f))\neq\emptyset$ and $W^u(\Lambda_{2}(\tilde f))\cap W^s(\Lambda_{1}(\tilde f))\neq\emptyset$ for all $\tilde f\in\mathcal{U}$, where $\Lambda_1(\tilde f)$ and $\Lambda_2(\tilde f)$ are the hyperbolic continuations of $\Lambda_1(f)$ and $\Lambda_2(f)$, respectively.
\end{defi}

Note that, due to the Kupka-Smale Theorem, intersections between invariant manifolds of periodic orbits of a generic system are transverse. Hence, at least one of the sets $\Lambda_1$ and $\Lambda_2$ must be non-trivial.  The heterodimensionality implies that, within either $W^u(\Lambda_{1})\cap W^s(\Lambda_{2})$ or $W^u(\Lambda_{2})\cap W^s(\Lambda_{1})$, any particular intersection of the corresponding stable and unstable leaves can be destroyed by a $C^r$-small perturbation of the system.
However, as at least one of the hyperbolic sets is non-trivial, say $\Lambda_{1}$, the laminations $W^s(\Lambda_{1})$ and $W^u(\Lambda_{1})$ contain uncountably many leaves, so the possibility to destroy the non-transverse intersection for any given pair of leaves does not preclude the existence, for any $C^1$-close system, of a pair of leaves which do have such intersection. 
\begin{defi}[Homoclinic relations]
Two transitive hyperbolic sets of a diffeomorphism are {\em homoclinically related} if they have the same index, and the unstable manifold of each set intersects the stable manifold of the other set transversely.
\end{defi}

\begin{thm}\label{thm:eff}
Let $f$ be a system of class $C^r$ ($r=1,\dots,\infty,\omega$) which has a generic homoclinic tangency to a hyperbolic periodic orbit $L$ with $d_{\rm eff}>1$. Then, arbitrarily $C^r$-close to $f$, there exists a system $\tilde f$ having  $C^1$-robust heterodimensional dynamics, which involve two non-trivial uniformly-hyperbolic basic  sets $\Lambda_{cs}$ and $\Lambda_{cu}$. Moreover, one of the hyperbolic sets is homoclinically related to the continuation of $L$.
\end{thm}

A precise definition of the effective dimension is given in Section \ref{sec:intro2}, Definition \ref{defi:d_eff}. The discussion leading to this definition explains that
the condition $d_{\rm eff}>1$ is also {\em necessary} for the birth of heterodimensional cycles/dynamics from the bifurcation of a generic homoclinic tangency to a hyperbolic periodic orbit.

We also prove that the sets $\Lambda_{cs}$ and $\Lambda_{cu}$ in Theorem \ref{thm:eff} are {\em center-stable (cs)} and, respectively, {\em center-unstable (cu) blenders}. This is the key fact responsible for the robustness of the heterodimensional cycle. The concept of blender and the theory of the $C^1$-robustness of the associated heterodimensional intersections were introduced by Bonatti and D\'iaz in \citep{BD96,BD08}; see the discussion in Section \ref{sec:intro3}.

Theorem \ref{thm:eff} is an immediate consequence of a more general statement, Theorem \ref{thm:sf} (see Section \ref{sec:intro2}), where we show that {\em given a $C^r$ ($r=3,\dots,\infty,\omega$) family  of perturbations of the original system with a generic homoclinic tangency with $d_{\rm eff}>1$, if the family is proper, i.e., it depends on at least two parameters and satisfies certain  genericity conditions, then there are open regions of parameter values for which the system has  $C^1$-robust heterodimensional dynamics involving two blenders with different indices,
$\Lambda_{cs}$ and $\Lambda_{cu}$}.

Since  periodic orbits are dense in a hyperbolic basic set, and  the stable/unstable manifold of each periodic orbit is dense in the stable/unstable lamination of the set (see e.g. \citep[Theorem 6.4.15]{KH95}),  Theorem \ref{thm:eff} implies that $f$ is in the $C^r$-closure of  the $C^1$-open Bonatti-D\'iaz domain. 

We stress that conditions of Theorem \ref{thm:sf} are explicitly verifiable (see conditions C1--C5 in Section~\ref{sec:gc}), which allows for finding robust heterodimensional dynamics in systems with high regularity and with various additional structures or symmetries, or in any specific model of applied or theoretical interest. From the applied point of view, this result implies the {\em inhomogeneity of chaos}: for the parameter values from the above described open regions, any attractor containing any of these hyperbolic sets must also contain the other hyperbolic set, with a different number of positive Lyapunov exponents. Thus, the result allows for a significant extension of the class of known attractors with such inhomogeneity property; for earlier examples see e.g. in \cite{TS98,LT17,LT20}.    

The link between the existence of the homoclinic tangencies and heterodimensional cycles has been discussed before in several papers. For instance, it is shown in \citep{DR92} that homoclinic tangencies can be created from the so-called non-connected heterodimensional cycles. In \citep{Sh11}, the creation of heterodimensional cycles by $C^1$-small perturbations was studied for systems with no domination, which is closely related to the
existence of homoclinic tangencies. This direction was further pursued
in \citep{BB12,BCDG13,Lix17}. In \citep{BDP21}, the birth of robust heterodimensional dynamics from homoclinic tangencies with $d_{\rm eff}>1$ was shown under an additional assumption that the homoclinic tangency is connected to a pre-existing blender.
The birth of heterodimensional cycles from a pair of homoclinic tangencies with $d_{\rm eff}=1$ (e.g., due to a periodic perturbation of the Lorenz attractor) and from a pair of Shilnikov homoclinic loops was established in \cite{LT20} and \cite{Li16,LT17}, respectively. Heterodimensional cycles were  obtained in \citep{To19} from a transverse homoclinic to a saddle-focus that undergoes an Andronov-Hopf bifurcation. 
More recently, this result was applied to a class of homoclinic tangencies in \citep{To25}.

\subsection{Blenders and $C^1$-robust homoclinic tangencies}\label{sec:intro3}
We say that a transitive hyperbolic set has a {\em $C^1$-robust homoclinic tangency}
if its stable and unstable laminations have a non-transverse intersection for
all $C^1$-close systems. 
Note that the existence of a $C^1$-robust homoclinic tangency  implies that the system belongs to an $C^1$-open set where those having a homoclinic tangency to a saddle are $C^r$-dense ($r\geq 1$), i.e., it belongs to the $C^1$-open Newhouse domain.

Specific examples of $C^1$-robust homoclinic tangencies were constructed in \citep{Si72,As08}. In \citep{BD12}, Bonatti and D\'iaz showed  that (see \citep[Theorem 4.9]{BD12})  a $C^1$-robust homoclinic tangency appears at a generic unfolding of a quadratic homoclinic tangency to a periodic point which is homoclinically related to the so-called {\em blender-horseshoe}, a special blender with certain dynamical structure. The notion of blender is introduced by Bonatti and D\'iaz in \citep{BD96}.  Roughly speaking, it is a hyperbolic set
whose projection to a certain subspace contains an open subset of this subspace: with additional assumptions, this makes
an intersection of  the unstable or stable lamination of the hyperbolic set with a manifold of a lower than complementary dimension robust. In particular, this blender property persists for all $C^1$-small perturbations. There are several variations of the precise notion of a blender, see e.g. \citep{BD96,BDV05, BKR14,BBD16,Bi22}. 

In this paper, we use a constructive definition of a blender, as initiated in \citep{LT21,Li24b}, which we call {\em standard}.
The blender-horseshoe is, by definition, topologically conjugate to a Bernoulli scheme with
2 symbols, whereas the standard blender is allowed to be described by an arbitrary finite Markov chain. So, a blender-horseshoe is a center-unstable   standard blender  such that the finite Markov chain describing it has exactly 2 elements. One can see from the proof of  \citep[Theorem 4.9]{BD12} that this particular choice of the Markov chain
is immaterial, and the result holds for any  standard blender (both center-stable and center-unstable), see  \citep[Theorem 1]{Li24b}.

Let $f$ have a generic homoclinic tangency of effective dimension larger than 1. By \citep{GTS93},  $f$ is in the $C^r$ closure of an open subset $\mathcal{U}$ of the $C^2$-open Newhouse domain such that it contains a dense subset $\mathcal{U}'$ where each system has a homoclinic tangency to the continuation of $L$. By Theorem \ref{thm:eff}, given any $\tilde f\in\mathcal{U}'$, one  can perturb  $\tilde f$ such that a standard  blender appears, which is homoclinically related to the continuation of $L$. Since $\tilde f$ is in the $C^2$-open Newhouse domain, the continuation of $L$ acquires a homoclinic tangency by an additional arbitrarily $C^r$-small perturbation, which does not destroy the blender
nor the homoclinic relation between the blender and $L$.  As we mentioned, the construction in the proof of \citep[Theorem 4.9]{BD12} is applicable to such configuration. As a result, from
our Theorem, we can deduce the results in \citep[Section 3]{Li24b}, namely,

\begin{cormain}\label{thm:C1tangency}
The system $\tilde f$ in Theorem \ref{thm:eff} can be taken such that at least one of the two blenders $\Lambda_{cs}$ and $\Lambda_{cs}$ exhibits a $C^1$-robust homoclinic tangency.
\end{cormain}

In fact, a stronger result is established in \citep[Corollary 4]{Li24b} that, under the same assumption of Corollary \ref{thm:C1tangency}, one obtains the $C^1$-robust presence of  {\em uncountably many} orbits of homoclinic tangency to the blender.

 In terms of the geometry of the space of dynamical systems, Corollary \ref{thm:C1tangency} implies that an arbitrarily small, in $C^r$ for any $r=1,2,\dots, \infty,\omega$, perturbation of any system with a homoclinic tangency with $d_{\rm eff}>1$ puts the system inside the intersection of the $C^1$-open Bonatti-D\'iaz domain and  the $C^1$-open Newhouse domain.
Note that, by the works of Ures \citep{Ur95} and Moreira \citep{Mo11}, no hyperbolic basic sets in dimension two can produce  $C^1$-robust  tangencies, so the ``higher dimensionality'' requirement in Corollary \ref{thm:C1tangency} is indeed essential for the result.

It is noteworthy  that, by \citep[Theorem 5.5]{Li24b}, all blenders obtained in this paper, in particular, the blenders $\Lambda_{cs}$ and $\Lambda_{cu}$ in Theorem~\ref{thm:eff} are {\em arrayed blenders} in the terminology of \citep{Li24b}, which  can generate  parablenders introduced by Berger in \citep{Be16}. Parablenders were used there to build an  example  showing the existence of an open region, known as the ``Berger domain'', within the space of finite-parameter families of diffeomorphisms. In this domain, generic families exhibit infinitely many sinks for an open set of parameter values. For recent advances in this direction, see \citep{BR21,BCP22}. In a forthcoming paper \citep{LT24}, we prove that any diffeomorphism  with a  homoclinic  tangency with $d_{\rm eff}>1$ lies in the $C^r$ closure of the Berger domain, provided that the associated periodic point is of ``sink-producing'' type \citep{Tu96}.
\\
 
The  paper is organized as follows. In Section \ref{sec:intro2}, we make the precise setting of the problem and  state our main result, Theorem \ref{thm:sf}, the parametric version of Theorem \ref{thm:eff}. For that, we  classify generic homoclinic tangencies of effective dimension larger than one into two cases: {\em saddle-focus} and {\em bi-focus}.  Sections \ref{sec:firstreturn}--\ref{sec:robusthdc} are devoted to the saddle-focus case, where the main goal is to prove that a heterodimensional cycle involving the periodic point $L$  can be obtained within a two-parameter unfolding of its homoclinic tangency. A plan of the proof is given in the end of Section \ref{sec:prop}.
The bi-focus case is dealt with in Section \ref{sec:df} by reducing it to the saddle-focus case.
 
\section{Heterodimensional dynamics in two-parameter unfolding families}\label{sec:intro2}
\subsection{General setup}
We consider two closely related settings. In the first one, we take a $C^r$-diffeomorphism ($r=3,\dots,\infty,\omega$) $f$ of a $d$-dimensional ($d\geqslant 3$) manifold
and assume that $f$ has a hyperbolic periodic point $O$. In the second setting, we consider
a $C^r$ flow defined on a $(d+1)$-dimensional manifold. We denote the flow also by $f$, and assume that it has a periodic orbit $L$. When we do not distinguish between these two cases, we call $f$ {\em a system}.

Let the system $f$ have an orbit of homoclinic tangency, denoted by $\Gamma$, to the hyperbolic periodic orbit, i.e.,
the stable and unstable invariant manifolds $W^s$ and $W^u$ of the periodic orbit intersect each other non-transversely at the points of $\Gamma$. Our goal is to create heterodimensional cycles in a neighborhood of $L\cup \Gamma$. As shown in \citep{GST08}, the bifurcations of $\Gamma$ strongly depend on the effective dimension of the problem, which is determined, in a generic case, by the central multipliers of the periodic orbit.

In the case where $f$ is a diffeomorphism, let $\tau$ be the period of $O$, i.e., $f^\tau(O)=O$; we denote as
$L=\{O,f(O),\dots, f^{\tau-1}(O)\}$ the orbit of $O$. We define the local map $T_0$ as a restriction of $f^\tau$ to a small neighborhood $U_0$ of $O$. If $f$ is a flow, we pick a point $O\in L$ and take a small $d$-dimensional cross-section $U_0$
to $L$ though the point $O$, then the orbits of the flow $f$ starting in $U_0$ define the Poincar\'e return map to $U_0$, which we also call the local map and denote as $T_0$.

In both settings, $O$ is a hyperbolic fixed point of $T_0$. Let its local stable and unstable manifolds with respect to $T_0$ have dimensions $d^s$ and $d^u$, respectively; note that $d^s+d^u=d$.
By definition, $d^u$ is the index of the hyperbolic orbit $L$.

The multipliers $\lambda_1,\dots, \lambda_{d^s}, \gamma_1,\dots, \gamma_{d^u}$ of $O$ are defined as the eigenvalues of the derivative $\D T_0$ at $O$; we order them as follows:
\begin{equation}\label{eq:eigenvalue}
|\lambda_{d^{s}}|\leqslant \dots\leqslant|\lambda_1|<1<|\gamma_1|\leqslant\dots\leqslant |\gamma_{d^u}|.
\end{equation}
The center-stable and center-unstable multipliers are those nearest to the unit circle from the inside and, respectively, from the outside. Up to an arbitrarily $C^r$-small perturbation, we can always assume that the central multipliers are just $\lambda_1$ and $\gamma_1$ along with their complex conjugates, if any.

One, then, distinguishes the following three generic cases:
\begin{itemize}
\item saddle:\\
$\lambda_1=\lambda$ and $\gamma_1=\gamma$ are real;
\item saddle-focus:\\
(2,1): $\lambda_1=\lambda^*_2=\lambda e^{i\omega}$ ($0<\omega<\pi$) and $\gamma_1=\gamma$ is real, so
$|\lambda_3|<\lambda<1<|\gamma|<|\gamma_2|$, or \\[5pt]
(1,2): $\lambda_1=\lambda$ is real and $\gamma_1=\gamma^*_2=\gamma e^{i\omega}$ ($0<\omega<\pi$), so
 $|\lambda_2|<|\lambda|<1<\gamma<|\gamma_3|$;
\item bi-focus (or saddle-focus of type (2,2) in the terminology of \citep{GST08}):\\
 $\lambda_1=\lambda^*_2=\lambda e^{i\omega_1}$ ($0<\omega_1<\pi$) and $\gamma_1=\gamma^*_2=\gamma e^{i\omega_2}$ ($0<\omega_2<\pi$), so $|\lambda_3|<\lambda<1<\gamma<|\gamma_3|$.
\end{itemize}

Under some genericity assumptions (conditions C3, C4 in Section \ref{sec:gc}), for the system $f$ and for all $C^1$-close systems, there is an invariant manifold  $\mathcal{M}^c$ transverse to the strong-stable and strong-unstable directions, such that every orbit that never leaves a small neighborhood $U$ of $L\cup \Gamma$ lies in $\mathcal{M}^c$ \cite{GTS93b}.
In the case of saddle, this manifold is two-dimensional when $f$ is a diffeomorphism and three-dimensional when $f$ is a flow. Therefore, since two-dimensional maps or three-dimensional flows cannot have saddle periodic orbits with different dimensions of the unstable manifolds, no heterodimensional cycles can be born in a small neighborhood of $\Gamma\cup L$ in such case.\footnote{The situation changes when the two-dimensional reduction is prevented due to a violation of the genericity conditions, or by extending the neighborhood, see e.g. \citep{LT20}.}  
In this paper we, therefore, focus on the case where at least one of the central multipliers is non-real.

However, in the saddle-focus case, there is another restriction which prevents from
the birth of heterodimensional cycles out of a homoclinic tangency. Namely, when
$|\lambda\gamma|<1$, one cannot expect the emergence of heterodimensional cycles in a small neighborhood $U$ of $L\cup\Gamma$ in the case of the saddle-focus $(2,1)$ without additional assumptions. Indeed, when
$|\lambda\gamma|<1$, the dynamics on the invariant manifold $\mathcal{M}^c$ is sectionally dissipative, i.e., the return map to the neighborhood $U_0$ of $O$ contracts two-dimensional areas. The contraction of areas implies that, restricted to $\mathcal{M}^c$, every hyperbolic periodic orbit that can be born near
$\Gamma\cup L$ is either a sink or has the same index as $L$, see \citep{GST08}, which means no heterodimensional cycles in $\mathcal{M}^c\cap U$.
On the other hand, since $\mathcal{M}^c$ contains all the orbits which never leave $U$, if the system had a heterodimensional cycle in $U$, then there would be  heteroclinic orbits that connect the saddles of the cycle lying entirely in $U$ and hence in $\mathcal{M}^c$, giving a heterodimensional cycle in $\mathcal{M}^c$. 
Therefore, in this case the creation of heterodimensional cycles near $\Gamma\cup L$ is impossible. The case of saddle-focus (1,2) reduces to the case of saddle-focus (2,1) by the time reversal (in the discrete-time case this means the change of $f$ to $f^{-1}$). Therefore, no heterodimensional cycle can be born in a small neighborhood of a generic homoclinic tangency to a saddle-focus (1,2) when $|\lambda\gamma|>1$ (when the time is reversed, the multipliers are replaced by their inverse).

As defined in \cite{GST08}, in the remaining cases, the effective dimension $d_{\rm eff}$  associated to the bifurcation of a generic homoclinic tangency (the maximal possible number of zero Lyapunov exponents for periodic orbits that can emerge from this bifurcation) is larger than $1$. 
\begin{defi}[Effective dimension]\label{defi:d_eff}
The bifurcation of a generic homoclinic tangency corresponds to $d_{\rm eff}>1$ if $L$ is either 
\begin{itemize}[nosep]
\item a saddle-focus of type (2,1) with $|\lambda\gamma|\geq 1$, or
\item a saddle-focus of type (1,2) with $|\lambda\gamma|\leq 1$, or
\item a bi-focus.
\end{itemize}
\end{defi}

We consider only homoclinic tangencies satisfying $d_{\rm{eff}}>1$. Moreover, we assume that 
$$|\lambda\gamma|>1.$$
Note that if $|\lambda\gamma|=1$, it can be made larger than $1$ by a $C^r$-small perturbation of the system; the case $|\lambda\gamma|<1$ is reduced to the case $|\lambda\gamma|>1$ by the time reversal.

\subsection{Proper unfoldings}\label{sec:prop}
Our main results on the birth of heterodimensional cycles are obtained by considering
{\em proper} two-parameter unfoldings of generic homoclinic tangencies with
$d_{\rm eff}>1$. The genericity conditions C1--C5 are described
in Section~\ref{sec:gc}; they can be fulfilled for any homoclinic tangency after an arbitrarily $C^r$-small perturbation (see Section \ref{sec:gc}). The notion of a proper unfolding is defined as follows.

Assume the tangency of $W^u(L)$ and $W^s(L)$ at the points of $\Gamma$ is quadratic, i.e., the genericity conditions C1 and C2 are satisfied. Then, for every system $C^r$-close to $f$, one can define the {\em splitting parameter} $\mu$ -- a smooth functional equal to the (signed) distance between $W^u(O)$ and $W^s_{loc}(O)$ near a certain point of $W^s_{loc}(O)$; this point is chosen such that, for any system $g$ which is close to $f$,
an orbit of homoclinic tangency close to $\Gamma$ exists if and only if $\mu(g)=0$.

Other important parameters are the arguments of the complex central multipliers of $L$, i.e.,
$\omega$ in the saddle-focus case and $\omega_1,\omega_2$ in the bi-focus
case. For the bi-focus with $\lambda\gamma>1$, we use $\omega_1$ as a bifurcation control parameter and denote $\omega=\omega_1$ (when $\lambda\gamma<1$, one should choose $\omega=\omega_2$). Since the hyperbolic point $O$ admits unique continuations for all nearby systems and the dependence of $O$ on the system is smooth, the parameters $\omega$ all depend smoothly on the system.

A family $\{f_\varepsilon\}$ of $C^r$ systems with $f_{\varepsilon^*}=f$ is a {\em proper unfolding} if it depends on at least two parameters and
\begin{equation}\label{eq:para_sf}
\rank \left.
\dfrac{\partial (\mu(f_\varepsilon),\omega(f_\varepsilon))}{\partial \varepsilon}
\right|_{\varepsilon=\varepsilon^*}=2.
\end{equation}
\begin{thm}\label{thm:sf}
Suppose $|\lambda\gamma|>1$ and genericity conditions C1--C5 of Section~\ref{sec:gc} are satisfied.  Assume that either
\begin{itemize}
\item $L$ is a saddle-focus of type (2,1),  or
\item $L$ is a bi-focus, and $\omega_1/\pi,\omega_2/\pi,1$ are rationally independent.
\end{itemize}
Then, for any proper unfolding $\{f_\varepsilon\}$ and any neighborhood $U$ of $O\cup L$, there exist open regions of $\varepsilon$ values accumulating on $0$, such that the corresponding system  $f_{\varepsilon}$ has  $C^1$-robust heterodimensional dynamics in $U$, which involve a standard cs-blender $\Lambda_{cs}$ of index $d^u$ containing the continuation of $L$ and a standard cu-blender $\Lambda_{cu}$ of index $d^u+1$.
\end{thm}

\begin{rem}
Theorem~\ref{thm:sf} shows that the blender in Corollary~\ref{thm:C1tangency} that exhibits a robust tangency  is   homoclinically related to the continuation of the original periodic orbit $L$. Moreover, the blender is center-stable  if $\lambda\gamma>1$ and is center-unstable if $\lambda\gamma<1$.
\end{rem}

Condition (\ref{eq:para_sf}) means that we can always choose the values of $\mu$ and $\omega$ as
the first two components of the vector of the parameters $\varepsilon$, i.e., we can assume that 
$\varepsilon=(\mu,\omega,\hat \varepsilon)$, where $\hat\eps$ represents the remaining components of $\eps$. It is obvious that it is enough to prove Theorem \ref{thm:sf}
only for two-parameter proper families $f_{\mu,\omega}$: if a proper family depends on more than two parameters, one finds the corresponding region of robust heterodimensional dynamics  in the two-parameter subfamily
given by constant $\hat\varepsilon$ (the robustness implies that the heterodimensional dynamics persist also when $\hat\varepsilon$ is perturbed).

\begin{rem}
$C^1$-robust tangencies can emerge within the family in Theorem \ref{thm:sf} if one shows that conditions C1--C5 can be achieved by changing $\mu$ and $\omega$. Indeed, as discussed above Corollary \ref{thm:C1tangency}, there are four steps for achieving robust homoclinic tangencies: (1) perturb the original system into the Newhouse domain; (2) obtain conditions C1--C5; (3) create robust heterodimensional cycles and hence blenders; (4) apply the Bonatti-D\'iaz construction \citep{BD12} (or the modified one in \citep{Li24b}). By \citep{GTS93} and \citep{BD12},  steps (1) and, respectively, (4) can be done within a generic one-parameter family $f_\eps$ with $d\mu/d\eps \neq 0$ while step (3) follows from Theorem \ref{thm:sf}.
\end{rem}

Note that the robustness implies that  the heterodimensional dynamics   persist for an open set of parameters $\varepsilon$ around $\varepsilon=\varepsilon_j$.

Reversing the time, we obtain the same result when $|\lambda\gamma|<1$: one should change a saddle-focus of type (2,1) to a saddle-focus (1,2) and, in the case of bi-focus, interchange $\omega_1$ and $\omega_2$, replace condition C5 by those applied to the time-reversed system,
and replace the cs-blender $\Lambda_{cs}$ by a cu-blender $\Lambda_{cu}$ and the cu-blender
$\Lambda_{cu}$ of index $d^u+1$ by a cs-blender $\Lambda_{cs}$ of index $d^u-1$.
Since the genericity conditions on the homoclinic tangency can be achieved by an arbitrarily $C^r$-small perturbation (see Section \ref{sec:gc}), Theorem \ref{thm:sf} implies
Theorem \ref{thm:eff}.

The proof of Theorem \ref{thm:sf} for the saddle-focus case is done in Sections \ref{sec:firstreturn}--\ref{sec:robusthdc}. As it is enough to do the proof only for 2-parameter proper unfoldings, so this will be the case we consider.
Recall that $U$ is a small neighborhood of the periodic orbit $L$ and the orbit $\Gamma$ of homoclinic tangency. In Section \ref{sec:firstreturn}, we find a formula for the first-return map along the orbits in $U$. Then, in Section \ref{sec:QQ}, following \cite{GST08}, we find values of $\mu$ for which the first-return map has a hyperbolic fixed point $Q\in U$ with index different from that of $O$. After that, we create a non-transverse intersection $W^u(O)\cap W^s(Q)$ by adjusting the value of $\omega$ (Section \ref{sec:nontr}). We also use the area expansion
(guaranteed by the condition $|\lambda\gamma|>1$) to show the existence of the transverse intersection $W^s(O)\cap W^u(Q)$ (Section \ref{sec:trans}). This gives us a heterodimensional cycle involving the orbit $L$ of $O$ and the orbit of $Q$, summarized as Theorem \ref{thm:sfsimple}. Finally, we finish the proof in Section \ref{sec:robusthdc} by establishing that this heterodimensional cycle
 satisfies conditions of   \citep[Theorem 7]{LT21}, which implies that an additional change of $\omega$ creates the blenders
$\Lambda_{cs}$ and $\Lambda_{cu}$ with $C^1$-robust heterodimensional intersections.
The theorem for the bi-focus case is proved in Section \ref{sec:df} where we show that the first-return map in the bi-focus case can be brought to the same form as in the saddle-focus case.

It should be noted that, although achieving robust heterodimensional dynamics from the heterodimensional cycle found in Theorem \ref{thm:sfsimple} is based on the existing results of \citep{LT21}, doing so in the two-parameter family poses a difficulty -- one needs to keep track of a secondary parameter ($t$ in Lemma \ref{lem:index2}) in the series of perturbations used to obtain the cycle, so that the conditions required by \citep{LT21} can be verified.  This leads to intricate estimates on the derivatives with respect to the parameters, as will be seen through the computations. 

\subsection{Genericity conditions}\label{sec:gc}
We now describe precisely the genericity conditions of Theorem \ref{thm:sf}. In this section, all the invariant manifolds are associated with $O$, and we hence omit `$O$' in their notations for simplicity. Recall that orbits of the system $f$ define the return map (the local map) $T_0$ of a small $d$-dimensional disc $U_0$ around the periodic point $O$, the hyperbolic fixed point of $T_0$. The homoclinic orbit $\Gamma$ intersects $U_0$ in two sequences of points, converging to $O$ -- one sequence lies in $W^u_{loc}$ and the other in $W^s_{loc}$.  We take two points $M^+\in W^s_{loc} \cap\Gamma$ and $M^-\in W^u_{loc}\cap\Gamma$. Since $M^-$ and $M^+$ belong to the same orbit $\Gamma$, the orbits that pass near $\Gamma$ define the global map $T_1$ from a small neighborhood of $M^-$ in $U_0$ to a small neighborhood of $M^+$ in $U_0$. We assume that $\Gamma$ is an orbit of a homoclinic tangency, which means that the image $T_1(W^u_{loc})$ is tangent at $M^+$ to $W^s_{loc}$. We require that

\noindent {\bf C1.} the tangent spaces to the manifolds $T_1(W^u_{loc})$ and $W^s_{loc}$ at $M^+$ have only one common vector (up to a multiplication to a scalar), that is, $\dim( \mathcal{T}_{M^+}T_1(W^u_{loc})\cap \mathcal{T}_{M^+}(W^s_{loc}))=1$.

\noindent {\bf C2.} the tangency between $T_1(W^u_{loc})$ and $W^s_{loc}$ at $M^+$ is {\em quadratic}.

\noindent These two conditions can be formulated as follows \citep{NPT83}: there exist smooth coordinates $(s_1,s_2,t_1,t_2)\in\mathbb R\times\mathbb R^{d^s-1}\times\mathbb R\times\mathbb R^{d^u-1}$ near $M^+$ such that $W^s_{loc}=\{t_1=0,t_2=0\}$ and $T_1(W^u_{loc})=\{t_1=h(s_1),s_2=0\}$ for some function $h$ satisfying $h(0)=0,h'(0)=0,h''(0)\neq 0$.\\

Recall that the multipliers $\lambda_1,\dots, \lambda_{d^s}, \gamma_1,\dots, \gamma_{d^u}$ of $O$ are the eigenvalues of the derivative $\D T_0$ at $O$; the center-stable and center-unstable multipliers are those nearest to the unit circle from the inside and, respectively, from the outside. Denote by $d^{cs}$ and $d^{cu}$ the numbers of center-stable and center-unstable multipliers. Then, $d^{cs}=2,d^{cu}=1$ in the saddle-focus case (as mentioned before we only consider the $(2,1)$ case, since the result for the $(1,2)$ case is obtained by the time reversal), and $d^{cs}=2,d^{cu}=2$ in the bi-focus case. These are the only cases we consider in this paper.

If $d^s>d^{cs}=2$ or $d^u>d^{cu}$, we need more genericity assumptions. Recall that (see e.g. \citep{SSTC1}) there exists a strong-stable foliation $\mathcal{F}^{ss}$ of $W^s$ consisting of $(d^s-d^{cs})$-dimensional leaves. In particular, the leaf through $O$ is the strong-stable manifold $W^{ss}$ which is tangent at $O$ to
the invariant subspace of $DT_0$ that corresponds to the multipliers $\lambda_{d^{cs}+1},\dots,\lambda_{d^s}$. There also exists a $(d^{cs}+d^u)$-dimensional extended unstable invariant manifold $W^{uE}$, which contains $W^u$ and is tangent at $O$ to the eigenspace corresponding to $\lambda_1,\dots,\lambda_{d^{cs}},\gamma_1,\dots,\gamma_{d^u}$. The manifold $W^{uE}$ is not unique, but any two of them are tangent to each other at the points of $W^u$. Similarly, when $d^u>d^{cu}$, there exists a strong-unstable foliation $\mathcal{F}^{uu}$ of $W^u$ consisting of $(d^u-d^{cu})$-dimensional leaves, which includes the strong-unstable manifold $W^{uu}$ corresponding to the multipliers $\gamma_{d^{cu}+1},\dots,\gamma_{d^u}$. We also have a $(d^{cu}+d^s)$-dimensional extended stable invariant manifold $W^{sE}$, which contains $W^s$ and is tangent at $O$ to the eigenspace corresponding to $\gamma_1,\dots,\gamma_{d^{cu}},\lambda_1,\dots,\lambda_{d^s}$. It is not unique but any two of such manifolds are tangent to each other at points of $W^s$. We require that

\noindent {\bf C3.} if $d^{cs}<d^s$, then $\D T_1(\mathcal{T}_{M^-}W^{uE})$ is transverse to the leaf of $\mathcal{F}^{ss}$ through $M^+$; if $d^{cu}<d^u$, then $\D T^{-1}_1(\mathcal{T}_{M^+}W^{sE})$ is transverse to the leaf of $\mathcal{F}^{uu}$ through $M^-$.

\noindent {\bf C4.} the homoclinic tangency does not belong to the strong-stable or strong-unstable manifolds, namely, $\Gamma\not\subset W^{ss}\cup W^{uu}$.

Conditions C1--C4 are quite standard, see \cite{NPT83,GST08}. These conditions are open in the space of systems having the homoclinic tangency (i.e., on the codimension-1 surface $\{\mu=0\}$); they do not depend on the choice of the points $M^+$ and $M^-$, and condition C3 does not depend on the choice of the manifolds $W^{uE}$ and
$W^{sE}$.  

The next genericity condition is new for the theory of homoclinic tangency and is specific for the case of complex center multipliers (similar conditions were employed for the study of heterodimencional cycles in \cite{LT21}). To describe this condition, we introduce, following \citep{GST08}, $C^r$ coordinates $(x,y,u,v)\in\mathbb{R}^{d^{cs}}\times\mathbb{R}^{d^{cu}}\times\mathbb{R}^{d^s-d^{cs}}\times\mathbb{R}^{d^u-d^{cu}}$
straightening the local invariant manifolds and the leaves of the local invariant foliations $\mathcal{F}^{ss}$ and $\mathcal{F}^{uu}$ in $U$,
see (\ref{eq:maps:T0}) and (\ref{eq:maps:T0_nonlinear}). Thus, in the new coordinates, 
\begin{align*}
&W^s_{loc}=\{y=0,v=0\},\qquad W^{ss}_{loc}=\{x=0,y=0,v=0\},\\
&W^u_{loc}=\{x=0,u=0\},\qquad W^{uu}_{loc}=\{y=0,x=0,u=0\},
\end{align*}
and the foliations $\mathcal{F}^{ss}$ and $\mathcal{F}^{uu}$ are , respectively, $\{x=const,y=0,v=0\}$ and $\{y=const,x=0,u=0\}$. Also, in these coordinates, the extended manifold $W^{uE}$ is tangent to $\{u=0\}$ at $\{x=0,u=0\}$ and $W^{sE}$ is tangent to $\{v=0\}$ at $\{y=0,v=0\}$. Thus, we have $M^+=(x^+,0,u^+,0)$ and $M^-=(0,y^-,0,v^-)$, and condition C4 reads
\begin{equation}\label{eq:condition4}
x^+=(x_1^+,x_2^+)\neq 0 \qquad \mbox{and}\qquad y^-\neq 0.
\end{equation}

We choose coordinates in such a way that $T_0$ acts linearly on the foliations $\mathcal{F}^{ss}$
and $\mathcal{F}^{uu}$, i.e., if we denote $T_0|_{W^s_{loc}}:(x,u)\mapsto (\bar x,\bar u)$  and $T_0|_{W^u_{loc}}:(y,v)\mapsto (\bar y,\bar v)$, then
\begin{equation}\label{eq:intro:2}\!\!\!\!
\begin{array}{ll}
\begin{array}{c}\mbox{Saddle-focus}\\\mbox{case:}\end{array}&\!\bar{x} = \lambda\begin{pmatrix}
\cos \omega &-\sin\omega \\
\sin \omega & \cos\omega
\end{pmatrix}\!x \quad\mathrm{and}\quad \bar{y}=\gamma y;\\[15pt]
\begin{array}{c}\mbox{Bi-focus}\\\mbox{case:}\end{array}&\!\!\bar{x} = \lambda\!\begin{pmatrix}
\!\cos \omega_1 &\!\!-\sin\omega_1 \\
\!\sin \omega_1 &\!\!\cos\omega_1
\end{pmatrix}\!x \;\;\mathrm{and}\;\; \bar{y}=\gamma\!\begin{pmatrix}
\!\cos \omega_2 &\!\!-\sin\omega_2 \\
\!\sin \omega_2 &\!\!\cos\omega_2
\end{pmatrix}\! y.\!\!\!\!\!\!\!\!\!\!\!
\end{array}
\end{equation}
This can be always achieved, see e.g. Lemma 6 of \citep{GST08}. 

Now assume that condition C1 is satisfied, and let $\beta\in \mathbb{R}^2$ be the $(x_1,x_2)$-component of the tangent vector between $T_1(W_{loc}^{u})$
and $W^s_{loc}$ at the point $M^+$. We impose that

\noindent {\bf C5.} the vector $(x^+_1,x^+_2)$ (the $x$-component of $M^+$) is not parallel to the vector $\beta$.

This condition is independent of the choice of coordinates that straighten the local manifolds and foliations and linearise the central dynamics as in \eqref{eq:intro:2}. Indeed, for any $C^r$ transformation which keeps the central dynamics linear, its restriction to $W^s_{loc}$ and $W^u_{loc}$ acts
linearly in the central coordinates, which means that the vectors $a$ and $x^+$ are multiplied to the same matrix and remain non-parallel.

It is not hard to show that all these genericity conditions are $C^r$-open and if they are not fulfilled for a given orbit of homoclinic tangency,
then they can be achieved by an arbitrarily $C^r$-small ($r=2,\dots, \infty,\omega$) perturbation of the system (in the smooth case
one adds a local perturbation to $f$; in the analytic case one uses the scheme described in \citep{BT86,GTS07}).

\section{First-return map in the saddle-focus case}\label{sec:firstreturn}
We start to prove Theorem \ref{thm:sf} when $O$ is a saddle-focus of type (2,1), where the center-stable multipliers are $\lambda_1=\lambda^*_2=\lambda e^{i\omega}$ $(0<\omega<\pi)$, and the center-unstable multiplier $\gamma_1=:\gamma$ is real, with the assumptions $|\lambda\gamma|>1$ and C1--C5. Note that conditions C3 and C4 imply that there is a three-dimensional $C^1$-smooth invariant manifold containing $O$ and the orbit $\Gamma$ of homoclinic tangency \cite{GTS93b,Tu96}. However, we do not consider reducing the problem to a three-dimensional one, since the method developed here is also applied to the bi-focus case, where no such 3-dimensional reduction exists.

In this section, we derive necessary formulas for the first-return maps near $\Gamma$.

\subsection{Local and global maps}\label{sec:maps}
A first-return map is a composition of some iteration of the local map $T_0$
defined in a small neighborhood $U_0$ of its fixed point $O$ and the global map $T_1$ which takes a small, of size $\delta^-$, cube $\Pi^-$  centered at the homoclinic point $M^-$ inside a size $\delta^+$ cube $\Pi^+$ around the homoclinic point $M^+$. Recall that the system  smoothly depends on parameters $\eps$, and the same holds true for the maps $T_0$ and $T_1$. It is enough to
 consider the case where $\eps=(\mu,\omega)$, where the splitting parameter $\mu$
 varies near $0$ and $\omega$ varies near some $\omega^*\in (0,\pi)$. 
Throughout the paper, we use $O(\cdot)$ to denote a term which is smaller than the quantity inside the bracket multiplied by a constant independent of 
the choice of $\eps$ near $\eps^*=(0,\omega^*)$ and the choice of $\delta^{\pm}$. Similarly, we use the notation $o(1)$ to denote terms which tend to $0$ (in some limit) uniformly for small $\eps-\eps^*$ and $\delta^{\pm}$.

By Lemma 6 in \citep{GST08}, there exists a $C^r$-coordinate transformation in $U_0$, which is, along with its derivatives up to order 2, of class $C^{r-2}$ with respect to parameters, and the local map $T_0$ after this transformation can be written in the following form:
\begin{equation}\label{eq:maps:T0}
\begin{aligned}
\bar x_1& = \lambda(\varepsilon) x_1\cos\omega(\varepsilon) - \lambda(\varepsilon) x_2\sin\omega(\varepsilon)  + p_1(x,y,u,v,\varepsilon),\\ 
\bar x_2& = \lambda(\varepsilon) x_1  \sin\omega(\varepsilon)  + \lambda(\varepsilon) x_2 \cos\omega(\varepsilon) +  p_2(x,y,u,v,\varepsilon),\\ 
\bar y  & = \gamma(\varepsilon)  y +  p_3(x,y,u,v,\varepsilon),\\
\bar u  & = A(\varepsilon)  u +  p_4(x,y,u,v,\varepsilon),\quad
\bar v  = B(\varepsilon)  v + p_5(x,y,u,v,\varepsilon),
\end{aligned}
\end{equation}
where all coefficients depend on $\eps$, $x=(x_1,x_2)\in \mathbb{R}^2$, the eigenvalues of $A$ and $B$ are $\lambda_3,\dots,\lambda_{d^s}$ and $\gamma_2,\dots,\gamma_{d^u}$, respectively, and functions $p$ satisfy
\begin{equation}\label{eq:maps:T0_nonlinear}
\begin{array}{l}
p_{1,2,3,4}(0,y,0,v,\varepsilon)=0,\quad p_{1,2,3,5}(x,0,u,0,\varepsilon)=0, \\[10pt]
\dfrac{\partial p_{1,2,4}}{\partial x}(0,y,0,v,\varepsilon)=0, \quad  \dfrac{\partial p_{3,5}}{\partial y}(x,0,u,0,\varepsilon)=0,
\end{array}
\end{equation}
for all sufficiently small $x,y,u,v$. One can see from (\ref{eq:maps:T0_nonlinear}) that as we mentioned in Section \ref{sec:gc}, in these coordinates the local invariant manifolds $W^s_{loc}$ and $W^u_{loc}$ are $\{y=0,v=0\}$ and $\{x=0,u=0\}$ the leaves of the local invariant foliations $\mathcal{F}^{ss}$ and $\mathcal{F}^{uu}$ are , respectively, $\{x=const,y=0,v=0\}$ and $\{y=const,x=0,u=0\}$, and the extended manifold $W^{uE}$ is tangent to $\{u=0\}$ at $\{x=0,u=0\}$ and $W^{sE}$ is tangent to $\{v=0\}$ at $\{y=0,v=0\}$. See Figure~\ref{fig:manifolds} for an illustration.

\begin{figure}[!h]
\begin{center}
\includegraphics[scale=1.2]{figs/fig1.pdf}
\caption{A schematic picture for the local manifolds and the strong-stable foliation, where in particular the transversality condition is drawn (i.e.,  condition C3 of Section~\ref{sec:gc}).}
\label{fig:manifolds}
\end{center}
\end{figure}

By Lemma 7 in \citep{GST08}, when $r\geq 3$, there exist functions $\hat p_i(x_1,x_2,\tilde y,u,\tilde v,\eps)$ $(i=1,\dots,5)$, which along with their first derivatives\footnote{More specifically, two derivatives with respect to parameters can be lost in the above coordinates, see \citep[Lemma 7]{GST08} for details.} with respect to coordinates and parameters  are uniformly bounded, such that
two points are related by
$(\tilde x_{1},\tilde x_{2},\tilde y,\tilde u,\tilde v)= T_0^k(x_{1},x_{2},y,u,v)$ if and only if 
\begin{equation}\label{eq:maps:T0k}
\begin{aligned}
 \tilde x_{1}&= \lambda^k(\eps) x_1\cos k\omega(\eps)  - \lambda^k(\eps) x_2\sin k\omega(\eps) + \hat \lambda^k \hat p_1(x_1,x_2,\tilde y,u,\tilde v,\eps),\\ 
 \tilde x_{2} &=\lambda^k(\eps) x_1  \sin k\omega(\eps)  + \lambda^k(\eps) x_2 \cos k\omega(\eps)  + \hat \lambda^k\hat p_2(x_1,x_2,\tilde y,u,\tilde v,\eps),\\ 
 y   &= \gamma^{-k}(\eps) \tilde y +  \hat \gamma^{-k}\hat p_3(x_1,x_2,\tilde y,u,\tilde v,\eps),\\
  \tilde u   &= \hat \lambda^k \hat p_4(x_1,x_2,\tilde y,u,\tilde v,\eps) ,\quad
 v  = \hat \gamma^{-k} \hat p_5(x_1,x_2,\tilde y,u,\tilde v,\eps),
\end{aligned}
\end{equation}
where $\hat \lambda < \lambda$, $\hat\gamma>|\gamma|$. We take $\hat \lambda$ and $\hat\gamma$  close to $\lambda$ 
and $|\gamma|$, respectively, so that
\begin{equation}\label{eq:lambdagamma}
\hat{\lambda}>\sqrt{\frac{\lambda}{|\gamma|}} > |\gamma|^{-1}>\hat\gamma^{-1},\qquad
\hat{\lambda}> |\lambda\gamma|\hat\gamma^{-1},\qquad 
\hat \lambda >\lambda^2,\qquad
\lambda<|\gamma|\hat\gamma^{-1}.
\end{equation}
Below we omit the dependence of the coefficients on the parameters when there is no ambiguity, but we will need 
to evaluate derivatives of coefficients with respect to parameters.

Recall that $M^+=(x^+,0,u^+,0)$ and $M^-=(0,y^-,0,v^-)$. The global map $T_1:(\tilde x_{1},\tilde x_{2},\tilde y,\tilde u,\tilde v)\in \Pi^-\mapsto (x_{1},x_{2},y,u,v)\in \Pi^+$, is given by
\begin{align*}
{x}_1 - x_1^+(\eps) &=a'_{11} \tilde x_{1} +a'_{12} \tilde x_{2} +b'_1(\tilde y - y^-) + a'_{14} \tilde u + a'_{15} (\tilde v-v^-) + \dots, \\
{x}_2 - x_2^+(\eps) &=a'_{21} \tilde x_{1} +a'_{22} \tilde  x_2 +b'_2(\tilde  y - y^-) + a'_{24} \tilde  u + a'_{25} (\tilde  v-v^-) + \dots, \\
y -y^+(\eps)&=  c'_1 \tilde  x_1 +c'_2 \tilde  x_2 +b'_3(\tilde  y - y^-)+ c'_4 \tilde  u + c'_5 (\tilde  v-v^-) + \dots, \\
u - u^+(\eps) &=a'_{41} \tilde  x_1 +a'_{42} \tilde  x_2 +b'_4(\tilde  y - y^-) + a'_{44} \tilde  u + a'_{45} (\tilde  v-v^-) + \dots ,\\
 v -v^+(\eps)&=a'_{51} \tilde  x_1 +a'_{52} \tilde  x_2 +b'_5(\tilde  y - y^-) + a'_{54} \tilde  u + a'_{55} (\tilde  v-v^-) + \dots ,
\end{align*}
where the dots here and in similar formulas below denote the second and higher order terms in the Taylor expansions. Here all coefficients depend on
$\eps$; note that $y^+(\eps^*)=0$, $v^+(\eps^*)=0$. 

Since $W^{sE}$ is tangent at $M^+$ to $\{v = 0\}$ and the leaf of $\mathcal{F}^{uu}$ through $M^-$ is given by
$(\tilde x,\tilde u, \tilde y-y^-)=0$, the transversality of $T_1^{-1} (W^{sE})$ to this leaf at $M^-$ (as given by condition C3) means that $\det a'_{55}\neq 0$. Therefore,
$\tilde v$ can be expressed as a function of $(v, \tilde x, \tilde u,\tilde y ,\eps)$ and the above formula for the map $T_1$ can be rewritten as
follows:
\begin{equation}\label{eq:maps:T1_cross00}
\begin{aligned}
{x}_1 - x_1^+ &=a_{11} \tilde  x_1 +a_{12} \tilde  x_2 +b_1(\tilde  y - y^-) + a_{14} \tilde  u + a_{15}  v + \dots, \\
{x}_2 - x_2^+ &=a_{21} \tilde  x_1 +a_{22} \tilde  x_2  + b_2(\tilde  y - y^-)+ a_{24} \tilde  u + a_{25}  v + \dots, \\
 y &= y^+ + c_1 \tilde  x_1 +c_2 \tilde  x_2 +b_3(\tilde  y - y^-)+ c_4 \tilde  u + c_5  v + \dots, \\
u - u^+ &=a_{41} \tilde  x_1 +a_{42} \tilde  x_2 +b_4(\tilde  y - y^-) + a_{44} \tilde  u + a_{45}  v + \dots ,\\
\tilde  v - v^- &=a_{51} \tilde  x_1 +a_{52} \tilde  x_2 +b_5(\tilde  y - y^-) + a_{54} \tilde  u + a_{55}  v + \dots ;
\end{aligned}
\end{equation}
as before, all coefficients are at least $C^{r-2}$ with respect to 
$\eps$.  
Note that we can perform any rotation of the $x$-coordinates without breaking the form (\ref{eq:maps:T0}) and (\ref{eq:maps:T0_nonlinear}) 
of the map $T_0$. We do such rotation in order to make 
$$b_2=0.$$ 

Since $W^s_{loc}$ near $M^+$ is given by $(y,v)=0$ and $W^u_{loc}$ near $M^-$ is given by $(\tilde x,\tilde u)=0$, the tangency means
of $T_1 (W^u_{loc})$ and $W^s_{loc}$ for $\varepsilon=\varepsilon^*$ (condition C1) 
$$b_3(\eps^*)=0.$$
Indeed, if we denote the tangents  to $W^s_{loc}$ at $M^+$ and to $W^u_{loc}$ at $M^-$  by $(x,0,u,0)$ and, respectively, $(0, \tilde y,0,\tilde v)$, then $b_3(\eps^*)\neq 0$ would imply that the zero vector 
is the only solution to the equation $( 0, \tilde y,0,\tilde v)=D T_1(x,0,u,0)$, i.e., the transversality of $T_1 (W^u_{loc})\cap W^s_{loc}$.

Under the condition $b_2=b_3=0$, we find that   $(b_1,0, 0, b_4, 0)$  is the common tangent vector between $T_1(W^u_{loc})$ and $W^s_{loc}$ at $M^+$. Its $x$-component is 
$\beta=(b_1,0)$, and condition C5 reads as
\begin{equation}\label{eq:x+2}
x^+_2 \neq 0, \qquad b_1\neq 0.
\end{equation}

Condition C2 (the quadraticity of the tangency) means that the quadratic terms in the $y$-equation of (\ref{eq:maps:T1_cross00})
contain the term $d (\tilde y-y^-)^2$ with non-zero $d$. Since $d\neq 0$, one can always choose $y^-(\eps)$ (smoothly dependent on $\eps$) such that the coefficient $b_3(\eps)$ remains zero for all $\eps$ close to $\eps^*$. 

Thus, we can find from \eqref{eq:maps:T1_cross00} the equations for $T_1(W^u_{loc})$ as
\begin{equation*}
\begin{aligned}
 y &= y^+  +c_5  v + O((x-x^+_1)^2+v^2), \\
{x}_2 - x_2^+ &= b_2b^{-1}_1({x}_1 - x_1^+) + (a_{25}-a_{15}b^{-1}_1)  v +O((x-x^+_1)^2+v^2), \\
u - u^+ &=b_4b^{-1}_1({x}_1 - x_1^+)   + (a_{45}-a_{15}b^{-1}_1)  v + O((x-x^+_1)^2+v^2) .
\end{aligned}
\end{equation*}
Obviously, $y^+(\eps)$ can serve as  the signed distance between $T_1(0,0,y^-(\eps),0,v^-(\eps))$ and $W^s_{loc}$. Therefore, we take it equal to the splitting parameter  $\mu$.

Since $W^{uE}$ is tangent to $\tilde u=0$ at $M^-$, the transversality of 
$T_1 (W^{uE})$ at $M^+$ to the leaf $(x-x^+, y, v)=0$ of $\mathcal{F}^{ss}$ for $\eps=\eps^*$ (condition C3) means that, if we denote the tangents at $M^-$ and $M^+$ by $(x,y,0,v)$ and $(\tilde x, \tilde y,0,\tilde v)$, then $0$ is the only solution to the equation $(\tilde x, \tilde y,0,\tilde v)=D T_1(x,y,0,v)$. This means that
\begin{equation}\label{eq:det}
\det \left(\begin{array}{ccc} a_{11} & a_{12} & b_1\\ a_{21} & a_{22} & 0 \\ c_1 & c_2 & 0\end{array} \right) \neq 0
\qquad\Rightarrow\qquad
b_1 \cdot \det \left(\begin{array}{cc} a_{21} & a_{22}\\
c_1 & c_2\end{array}\right) \neq 0.
\end{equation}
In particular,
\begin{equation}\label{eq:a_31}
c_1^2+c_2^2\neq 0.
\end{equation}

Summarizing, we obtain the following

\begin{lem}[cf. {\citep[Corollary 1]{GST08}}]
In the above coordinates, the map $T_1$ is given by
\begin{equation}\label{eq:maps:T1_cross}
\begin{aligned}
{x}_1 - x_1^+ &=a_{11} \tilde  x_1 +a_{12} \tilde  x_2 +b_1(\tilde  y - y^-) + a_{14} \tilde  u + a_{15}  v + \dots, \\
{x}_2 - x_2^+ &=a_{21} \tilde  x_1 +a_{22} \tilde  x_2  +  a_{24} \tilde  u + a_{25}  v + \dots, \\
 y &= \mu + c_1 \tilde  x_1 +c_2 \tilde  x_2 +c_4 \tilde  u + c_5  v + d (\tilde  y - y^-)^2  
  + \dots, \\
u - u^+ &=a_{41} \tilde  x_1 +a_{42} \tilde  x_2 +b_4(\tilde  y - y^-) + a_{44} \tilde  u + a_{45}  v + \dots ,\\
\tilde  v - v^- &=a_{51} \tilde  x_1 +a_{52} \tilde  x_2 +b_5(\tilde  y - y^-) + a_{54} \tilde  u + a_{55}  v + \dots ,
\end{aligned}
\end{equation}
where the dots in the $y$-equation stand for $O(|\tilde  y - y^-|^3 + (v, \tilde x, \tilde u)^2 + 
|\tilde  y - y^-| \cdot \|v, \tilde x,\tilde u\|)$. 
\end{lem}


\subsection{Normal form for the first-return map}\label{norforfir}

We define for each large $k$ the first-return map $T_k:= T_1\circ T^k_0$ that acts on $\Pi^+$:
$$
(x_1,x_2,y,u,v)
\xmapsto{T^k_0}(\tilde x_{1},\tilde x_{2},\tilde y,\tilde u,\tilde v)
\xmapsto{T_1}(\bar x_1,\bar x_2,\bar y,\bar u,\bar v).
$$
That is, for any 
$(x_1,x_2,y,u,v)\in\Pi^+$,  we have $(\bar x_1,\bar x_2,\bar y,\bar u,\bar v)=T_k(x_1,x_2,y,u,v)\in \Pi^+$ if and only if there is 
$(\tilde x_{1},\tilde x_{2},\tilde y,\tilde u,\tilde v)\in \Pi^-$ such that
$(\tilde x_{1},\tilde x_{2},\tilde y,\tilde u,\tilde v)=T_0^k(x_1,x_2,y,u,v)$ and $(\bar x_1,\bar x_2,\bar y,\bar u,\bar v)=T_1(\tilde x_{1},\tilde x_{2},\tilde y,\tilde u,\tilde v)$.
In this section, in several steps, we obtain a normal form for the maps $T_k$ -- the formula   \eqref{eq:maps:Tk}.


\subsubsection{Composition.}
Combining \eqref{eq:maps:T0k} and \eqref{eq:maps:T1_cross} with replacing $(x_1,x_2,y,u,v)$ by $(\bar x_1,\bar x_2,\bar y,\bar u,\bar v)$ in \eqref{eq:maps:T1_cross}, yields 
\begin{equation}\label{eq:maps:Tk_0}
\begin{aligned}
\bar x_1 -x_1^+ &= \lambda^k \hat\alpha_1x_1 +\lambda^k \hat \beta_1 x_2 + (b_1+O(\lambda^k)) (\tilde y -y^-)+ O((\tilde y -y^-)^2)+O(\bar v) +O(\hat\lambda^k),\\
\bar x_2 -x_2^+&= \lambda^k \hat\alpha_2 x_1 +\lambda^k \hat \beta_2 x_2 +  O(\lambda^k)(\tilde y -y^-)+ O((\tilde y -y^-)^2)+O(\bar v) +O(\hat\lambda^k),\\
\bar y&= \mu + \lambda^k \alpha^*x_1 +\lambda^k \beta^*x_2 +  O(\lambda^k)(\tilde y -y^-)+ 
d  (\tilde y -y^-)^2
+ O((\tilde y -y^-)^3) +O(\bar v)+ O(\hat\lambda^k),\\
\bar u - u^+ &= \lambda^k \hat\alpha_4x_1 +\lambda^k \hat \beta_4 x_2 + (b_4+O(\lambda^k))  (\tilde y -y^-) + O((\tilde y -y^-)^2)+O(\bar v) +O(\hat\lambda^k),\\
\tilde v-v^-&=\lambda^k \hat\alpha_5 x_1 +\lambda^k \hat \beta_5 x_2 + (b_5+O(\lambda^k))  (\tilde y -y^-) + O((\tilde y -y^-)^2)+O(\bar v) +O(\hat\lambda^k),
\end{aligned}
\end{equation}
where
\begin{equation}\label{eq:maps:alphabeta}
\begin{aligned}
\alpha^* &= c_1\cos k\omega+c_2\sin k\omega,
\qquad
& \beta^* &= - c_1\sin k\omega+c_2\cos k\omega,
\\
\hat \alpha_i &= a_{i1}\cos k\omega+a_{i2}\sin k\omega,
\qquad
&\hat \beta_i &= -a_{i1}\sin k\omega+a_{i2}\cos k\omega,
\end{aligned}
\end{equation}
with $i=1,2,4,5$. 
Here the terms $O(\bar v)$ are from the terms involving $v$ in \eqref{eq:maps:T1_cross}. The terms $O(\lambda^k)$ and $O(\hat\lambda^k)$ are  from \eqref{eq:maps:T0k} (so their derivatives are  also estimated as $O(\hat{\lambda}^k)$). More specifically, the former  are from products of $(\tilde y-y^-)$ with either $\tilde x$ or $\tilde u$ in the dots of \eqref{eq:maps:T1_cross}, and the latter are from the remaining terms in the dots (which do not contain $v$) -- they are $O(\lambda^{2k})$ and are absorbed by the $O(\hat\lambda^k)$ terms because $\lambda^2 <\hat\lambda$. All the terms $O(v),O(\lambda^k),O(\hat\lambda^k)$ are  functions of $(x,\tilde {y},u,\tilde v)$.

\subsubsection{Shilnikov coordinates.}\label{sec: Shilnikov}
Formula \eqref{eq:maps:T0k} gives $y,v$ as functions of $(x_1,x_2,\tilde y,u,\tilde v)$. Therefore, following \citep{GS72,GTS07,GST08}, we can use the following coordinates (called
``Shilnikov coordinates'' in \citep{CDF90}) on $\Pi^+\cap T_0^{-k} (\Pi^-)$:
\begin{equation}\label{eq:maps:coor1}
X_1 = x_1 - x_1^+,\quad
X_2 = x_2 - x_2^+,\quad
U = u - u^+,\qquad
Y = \tilde y - y^-,\quad
V= \tilde v-v^-.
\end{equation}
In the new coordinates, the domain of the map
$T_k:(X_1,X_2,Y,U,V)\mapsto (\bar X_1,\bar X_2,\bar Y,\bar U,\bar V)$ is a small neighborhood of 0 whose size $\delta$ is independent of $k$; in particular, it does not shrink to 0 as $k\to\infty$. Below we find a formula for $T_k$ in these coordinates in the ``cross-form'' $(X_1,X_2,Y,U,\bar V)\mapsto (\bar X_1,\bar X_2,\bar Y,\bar U, V)$. 

\noindent\textbf{Step 1}
By \eqref{eq:maps:T0k} one has 
\begin{align*}
\bar y &= \gamma^{-k}(\bar {Y}+y^-) + \hat{\gamma}^{-k} \hat g_3(\bar X +x^+,\bar {Y}+y^-, \bar U+u^+, \bar {V}+v^-),\\
\bar v&=\hat\gamma^{-k}\hat g_5(\bar X+x^+,\bar {Y}+y^-, \bar U+u^+, \bar {V}+v^-),
\end{align*}
in \eqref{eq:maps:Tk_0}. Applying \eqref{eq:maps:coor1} to \eqref{eq:maps:Tk_0}, we obtain
\begin{equation}\label{eq*}
\begin{aligned}
\bar X_1 &= \lambda^k \hat\alpha_1(X_1+ x_1^+) +\lambda^k \hat \beta_1(X_2+x_2^+) + (b_1+ O(\lambda^k)) Y + O(Y^2) +
O(\hat \lambda^k) + O(\hat\gamma^{-k}),\\
\bar X_2 &= \lambda^k \hat\alpha_2(X_1+ x_1^+) +\lambda^k \hat \beta_2(X_2+x_2^+) + O(\lambda^k) Y + O(Y^2) +
O(\hat \lambda^k) + O(\hat\gamma^{-k}),\\
\bar Y &= \gamma^k \mu-y^- + \lambda^k\gamma^k \alpha^*(X_1+ x_1^+) +\lambda^k\gamma^k  \beta^*(X_2+x_2^+)  \\
 & \qquad\qquad\qquad\qquad + \gamma^k Y^2(d + O(Y))+ \gamma^k O(\lambda^k) Y+ \gamma^k O(\hat \lambda^k) + \gamma^k O(\hat\gamma^{-k}),\\
\bar U &= \lambda^k \hat\alpha_4(X_1+ x_1^+) +\lambda^k \hat \beta_4(X_2+x_2^+) + (b_4+ O(\lambda^k)) Y + O(Y^2) +
O(\hat \lambda^k) + O(\hat\gamma^{-k}),\\
V&=\lambda^k \hat\alpha_5(X_1+ x_1^+) +\lambda^k \hat \beta_5(X_2+x_2^+) + (b_5+ O(\lambda^k)) Y +  O(Y^2) +
O(\hat \lambda^k) + O(\hat\gamma^{-k}),
\end{aligned}
\end{equation}
where the terms $O(\lambda^k)$ are affine functions of $X$, the terms $O(\hat\lambda^k)$ are functions of $(X,Y,U,V)$, and the terms $O(\hat\gamma^{-k})$ are functions of $(\bar X, \bar Y, \bar U, \bar V)$, coming from the terms $O(\bar v)$ in \eqref{eq:maps:Tk_0}   and the relation $\bar v= \hat\gamma^{-k}\hat p_5(\bar X_1+x^+_1,\bar X_2+x^+_2, \bar Y+y^-,\bar V+v^-, \bar U+u^+)$ by \eqref{eq:maps:T0k}.

\noindent\textbf{Step 2}
The above formula describes the implicit relations between $(X,Y,U,V)$ and $(\bar X, \bar Y, \bar U, \bar V)$ (because the right-hand sides depend on $(\bar X, \bar Y, \bar U)$ and $V$). But, since the derivatives of the right-hand sides of the equations for $\bar X$, $\bar U$ and $V$ with respect to $(\bar X, \bar U, V)$ are small ($O(\hat \gamma^{-k})$ or $O(\hat \lambda^k)$), one can get rid of the dependence on $(\bar X, \bar U, V)$ in the right-hand sides of (\ref{eq*}) (i.e., express $(\bar X, \bar Y, \bar U, V)$ as functions of
$(X, Y, \bar Y, U, \bar V)$ by the Implicit Function Theorem). It is easy to see
that the equations (\ref{eq*}) keep their form, just the $O(\hat \lambda^k)$ and $O(\hat \gamma^{-k})$ terms are now functions of 
$(X,Y,\bar Y, U, \bar V)$ and the parameters $\varepsilon=(\mu,\omega)$, with the derivatives with respect to $(X,U,Y,\varepsilon)$ of order $O(\hat \lambda^k)$ and the derivatives with respect to $(\bar Y,\bar V)$ of
the smaller order $O(\hat \gamma^{-k})$ (by \eqref{eq:lambdagamma}, $\hat\lambda^k \gg \hat \gamma^{-k}$ for large $k$).

\noindent
\textbf{Step 3}
Next, we consider the equation for $\bar Y$ obtained in Step 2, which has the same form of the thrid equation in \eqref{eq*} but the right-hand side is now a function of $(X,Y,\bar Y,U,\bar V)$ with estimates on derivatives described in Step 2. Let us denote its right-hand side as $F$, namely,
\begin{align*}
F(X,Y,\bar Y,U,\bar V)=&
\gamma^k \mu-y^- + \lambda^k\gamma^k \alpha^*(X_1+ x_1^+) +\lambda^k\gamma^k  \beta^*(X_2+x_2^+)\\
 &\qquad\qquad+ \gamma^k Y^2(d + O(Y))+  \gamma^k O(\lambda^k) Y+ \gamma^k O(\hat \lambda^k) + \gamma^k O(\hat\gamma^{-k}).
\end{align*}
Note that by Step 2 one has $\partial F/\partial \bar Y = O(\gamma^k\hat\gamma^{-k})\ll 1$. Thus,  by applying the Implicit Function Theorem to the equation $\bar Y - F=0$ and using \eqref{eq:lambdagamma}, we can express $\bar Y$ as a function of $X,Y,U,\bar V$ as
\begin{equation}\label{eq:revise1}
\begin{aligned}
\bar Y =& \gamma^k \mu-y^- + \lambda^k\gamma^k \alpha^*(X_1+ x_1^+) +\lambda^k\gamma^k  \beta^*(X_2+x_2^+)
 + d\gamma^k Y^2 \\
 &\qquad\qquad+ \gamma^k(O(\hat\lambda^k)+O(\lambda^k)Y+O(Y^3)),
\end{aligned}
\end{equation}
where
\begin{align*}
\dfrac{\partial \bar Y}{\partial Y}
&=\dfrac{\partial F}{\partial Y}  
\left/\left(1-\dfrac{\partial F}{\partial \bar Y}\right)\right.
=\dfrac{\gamma^k (2dY +O(\lambda^k)+O(Y^2))}{1-\gamma^k\hat \gamma^{-k}+\dots}\\
&=
2d\gamma^k Y + \gamma^k(O(\lambda^k)+O(\gamma^k\hat\gamma^{-k})Y+O(Y^2)).
\end{align*}
Estimating the other partial derivatives  in the same way,  we obtain
\begin{equation}\label{ybarexp}
\bar Y = \gamma^k \mu-y^- + \lambda^k\gamma^k \alpha^*(X_1+ x_1^+) +\lambda^k\gamma^k  \beta^*(X_2+x_2^+)
 + d\gamma^k Y^2 + \gamma^k\hat h_3(X,Y,U,\bar V),
\end{equation}
where
\begin{equation}\label{yderest}
\begin{array}{l}
\hat h_3=O(\hat\lambda^k)+O(\lambda^k)Y+O(Y^3),\\
\dfrac{\partial \hat h_3}{\partial (X,U)} 
=O(\hat\lambda^k)+O(\lambda^k) Y,\quad
\dfrac{\partial \hat h_3}{\partial Y} 
=O(\lambda^k)+O(\gamma^k\hat\gamma^{-k})Y+O(Y^2),\quad
\dfrac{\partial \hat h_3}{\partial \bar V}
=O(\hat\gamma^{-k}),
\\
\dfrac{\partial \hat h_3}{\partial \mu}
=O(k\gamma^k\hat\gamma^{-k})+o(Y^2), \quad 
\dfrac{\partial \hat h_3}{\partial \omega}=k\gamma^k\hat\gamma^{-k} 
O(|\mu| + \lambda^k+ Y^2 ) + o(Y^2).
\end{array}
\end{equation}
Note that  the factor $\gamma^k\hat\gamma^{-k}$ in the derivatives with respect to parameters originates from the $O(\hat\gamma^{-k})$ terms in $F$ (like in the above computation for $\partial \bar Y/\partial Y$), and
 the factor $k$  appears from the differentiation of $\lambda^k$ and $\gamma^k$ -- recall that $\lambda$ and $\gamma$ are functions of the parameters $\mu$ and $\omega$ (so, for example, $d\lambda^k/d \eps=k\lambda^{k-1} \cdot d\lambda/d\eps$).

\noindent\textbf{Step 4}
Substituting the expression (\ref{ybarexp}) into the $O(\hat\gamma^{-k})$ terms in the equations for other variables, which are found in Step 2 with the same forms as the ones  in \eqref{eq*}, we find by the chain rule that their first derivatives with respect to $V$ are small. Hence, we can finally express $\bar X_1,\bar X_2, \bar U,V$ as functions of $X,Y,U,\bar V$ and
 obtain
\begin{equation}\label{eq:maps:Tk_1}
\begin{aligned}
\bar X_1 &= \lambda^k \hat\alpha_1(X_1+ x_1^+) +\lambda^k \hat \beta_1(X_2+x_2^+) + b_1 Y + \hat h_1(X,Y,U,\bar V),\\
\bar X_2 &= \lambda^k \hat\alpha_2(X_1+ x_1^+) +\lambda^k \hat \beta_2(X_2+x_2^+) + \hat h_2(X,Y,U,\bar V),\\
\bar U &= \lambda^k \hat\alpha_4(X_1+ x_1^+) +\lambda^k \hat \beta_4(X_2+x_2^+) + b_4 Y + \hat h_4(X,Y,U,\bar V),\\
V&=\lambda^k \hat\alpha_5(X_1+ x_1^+) +\lambda^k \hat \beta_5(X_2+x_2^+) + b_5 Y + \hat h_5(X,Y,U,\bar V),
\end{aligned}
\end{equation}
where
\begin{equation}\label{eq:hath3old}
\begin{array}{l}
\hat h_i =O(\hat \lambda^k)+O(\lambda^k)Y+O(Y^2), \\
\dfrac{\partial \hat h_i}{\partial(X,U)} 
=O(\hat\lambda^k)+O(\lambda^k)Y,\qquad
\dfrac{\partial \hat h_i}{\partial Y}
=O(\lambda^k)+O(Y), \qquad \dfrac{\partial \hat h_i}{\partial \bar V}
=O(\hat\gamma^{-k}),\\
\dfrac{\partial \hat h_i}{\partial \mu}
=O(k\gamma^k\hat\gamma^{-k})+O(Y^2), \qquad
\dfrac{\partial \hat h_i}{\partial \omega}
=O(k\lambda^k\gamma^k\hat\gamma^{-k})+O(Y^2), 
\qquad i=1,2,4,5.
\end{array}
\end{equation}

\subsubsection{Removing linear in $Y$ terms.}\label{sec:coor1}
Since $b_1\neq 0$, we can let
\begin{equation}\label{eq:maps:coor2}
 U^{new} = U - \dfrac{b_4}{b_1} X_1\quad\mbox{and}\quad V^{new} = V - b_5 Y.
\end{equation}
After this coordinate change, formulas \eqref{ybarexp}, \eqref{eq:maps:Tk_1}  become
\begin{equation}\label{eq:maps:Tk_2}
\begin{aligned}
\bar X_1 &= \lambda^k \hat\alpha_1 (X_1+ x_1^+) +\lambda^k \hat \beta_1 (X_2+x_2^+) + b_1 Y + \hat h_1(X,Y,U,\bar V),\\
\bar X_2 &= \lambda^k \hat\alpha_2 (X_1+ x_1^+) +\lambda^k \hat \beta_2 (X_2+x_2^+) + \hat h_2(X,Y,U,\bar V),\\
\bar Y& =  \gamma^k \mu-y^- + \lambda^k\gamma^k \alpha^*(X_1+ x_1^+) +\lambda^k\gamma^k  \beta^*(X_2+x_2^+)  + d \gamma^k Y^2  + \gamma^k\hat h_3(X,Y,U,\bar V),\\
\bar U &= \lambda^k (\hat\alpha_4-\dfrac{b_4}{b_1}\hat\alpha_1) (X_1+ x_1^+) +\lambda^k (\hat \beta_4-\dfrac{b_4}{b_1}\hat\beta_1) (X_2+x_2^+)  + \hat h_4(X,Y,U,\bar V),\\
V&=\lambda^k \hat\alpha_5 (X_1+ x_1^+) +\lambda^k \hat \beta_5 (X_2+x_2^+)+ \hat h_5(X,Y,U,\bar V).
\end{aligned}
\end{equation}
where the  functions $\hat h$ are different from those in \eqref{eq:maps:Tk_1}, but we keep the same notations for simplicity. Moreover, these new $\hat h$ satisfy the same estimates  \eqref{yderest} and \eqref{eq:hath3old}. Removing the linear in $Y$ terms helps to control the cone structure described in Lemma~\ref{lem:cones}.

\subsubsection{Shift of the origin.}
Since $d\neq 0$, it is not hard to see that one can shift the coordinate origin to $O(\lambda^k)$ so that the right-hand sides of the equations for $(\bar X,\bar U,V)$ would not have constant terms, i.e., they would vanish
at $(X,U,Y,\bar V)=0$, and, simultaneously, the linear in $Y$ term in $\hat h_3$ would vanish. Thus, we perform such coordinate transformation:
\begin{equation}\label{eq:maps:coor2.1}
(X,Y,U,V)^{new}=(X,Y,U,V)+O(\lambda^k),
\end{equation}
for some constant $O(\lambda^k)$, and formula \eqref{eq:maps:Tk_2} assumes the form
\begin{equation}\label{eq:maps:Tk_2.1}
\begin{aligned}
\bar X_1 &= \lambda^k \hat\alpha_1X_1 +\lambda^k \hat \beta_1X_2+ b_1 Y + \hat h_1(X,Y,U,\bar V),\\
\bar X_2 &= \lambda^k \hat\alpha_2X_1+\lambda^k \hat \beta_2 X_2 + \hat h_2(X,Y,U,\bar V),\\
\bar Y &= \hat\mu
+ \lambda^k\gamma^k \alpha^*X_1 +\lambda^k\gamma^k  \beta^*X_2 + d\gamma^k Y^2 + \gamma^k\hat h_3(X,Y,U,\bar V),\\
\bar U &=\lambda^k (\hat\alpha_4-\dfrac{b_4}{b_1}\hat\alpha_1) X_1 +\lambda^k (\hat \beta_4-\dfrac{b_4}{b_1}\hat\beta_1) X_2  + \hat h_4(X,Y,U,\bar V),\\
V&=\lambda^k \hat\alpha_5 X_1 +\lambda^k \hat \beta_5 X_2  + \hat h_5(X,Y,U,\bar V),
\end{aligned}
\end{equation}
where the constant term $\hat \mu$ of the $\bar Y$ equation is given by 
\begin{equation}\label{eq:maps:Tk_const}
\hat\mu=\gamma^k \mu-y^- +\lambda^k\gamma^k(\alpha^* x^+_1 + \beta^* x^+_2)+O(\hat\lambda^k\gamma^k).
\end{equation}
Here the functions $\hat h$ are different from those before and we keep the same notations for simplicity.  These new  $\hat h$ satisfy
\begin{equation}\label{hathnew} 
\begin{array}{l}
\hat h_i =O(\hat \lambda^k)(X,U,\bar V)+O(\lambda^k)Y+O(Y^2), \\
\dfrac{\partial \hat h_i}{\partial(X,U)} 
=O(\hat\lambda^k)+O(\lambda^k)Y,\qquad
\dfrac{\partial \hat h_i}{\partial Y}
=O(\lambda^k)+O(Y), \qquad \dfrac{\partial \hat h_i}{\partial \bar V}
=O(\hat\gamma^{-k}),\\
\dfrac{\partial \hat h_i}{\partial \mu}
=O(\gamma^k\hat\gamma^{-k})+O(Y^2), 
\qquad
\dfrac{\partial \hat h_i}{\partial \omega}
=O(\lambda^k\gamma^k\hat\gamma^{-k})+O(Y^2) + O(k\lambda^k)Y, 
\quad i=1,2,4,5,\\ \\
\hat h_3=O(\hat\lambda^k)(X,U,\bar V) + O(\lambda^k) XY +O(\lambda^k) Y^2+O(Y^3),\\
\dfrac{\partial \hat h_3}{\partial (X,U)} 
=O(\hat\lambda^k)+O(\lambda^k) Y,\quad
\dfrac{\partial \hat h_3}{\partial Y} 
=O(\hat\lambda^k)+O(\lambda^k)X+ O(\gamma^k\hat\gamma^{-k})Y+O(Y^2),\\
\dfrac{\partial \hat h_3}{\partial \bar V}
=O(\hat\gamma^{-k}),\quad
\dfrac{\partial \hat h_3}{\partial \mu}
=O(\gamma^k\hat\gamma^{-k})+o(Y^2), \\
\dfrac{\partial \hat h_3}{\partial \omega}=\gamma^k\hat\gamma^{-k} 
O(|\mu| + Y^2 + \lambda^k) + o(Y^2)+o(\lambda^k)Y,
\end{array}
\end{equation}
where we replaced $k\hat\gamma^{-k}$ by $\hat\gamma^{-k}$ in the estimates for the derivatives with respect to parameters by taking $\hat\gamma$ slightly smaller.

\subsubsection{Selecting the center coordinates.}\label{sec:coor3}
Note that $\alpha^*$ and $\beta^*$ cannot vanish at the same time since their phases have a difference $\pi/2$ by \eqref{eq:maps:alphabeta}. In fact, we will further consider $k$ such that 
\begin{equation}\label{eq:alpha*}
\alpha^*\neq 0.
\end{equation}
More precisely, we assume that $\alpha^*$ stays bounded away from zero -- for every $\omega\in (0,\pi)$ there are infinitely many such $k$.

Introduce the coordinates
\begin{equation}\label{eq:maps:coor3}
W=(X_2, U)^T,\quad
Z=\frac{1}{b_1} X_1 + \frac{\beta^*}{\alpha^* b_1} X_2,
\end{equation}
where the superscript $T$ denotes transpose, 
and rewrite formula \eqref{eq:maps:Tk_2.1} as
\begin{equation}\label{eq:maps:Tk}
\begin{aligned}
\bar Z &= \lambda^k \alpha_1 Z+ Y+\lambda^k\beta_1W+ h_1(Z,Y,W,\bar V),\\
\bar Y &= \hat\mu + \alpha_2 \lambda^k\gamma^k Z  +d \gamma^kY^2+ \gamma^k h_2(Z,Y,W,\bar V),\\
\bar W &=  \lambda^k \alpha_3 Z+\lambda^k\beta_3W+ h_3(Z,Y,W,\bar V),\\
V&= \lambda^k \alpha_4 Z+\lambda^k\beta_4W+ h_4(Z,Y,W,\bar V),
\end{aligned}
\end{equation}  
where
the coefficients $\alpha$ and $\beta$ are given by
\begin{equation}\label{eq:maps:Tk_coef}
\begin{array}{l}
\alpha_1 = \hat\alpha_1 + \hat\alpha_2\dfrac{\beta^*}{\alpha^*}
,
\quad 
\beta_1 =\dfrac{1}{b_1\alpha^*}
\left(-\hat\alpha_1\beta^*+\hat\beta_1\alpha^*-\hat\alpha_2\dfrac{(\beta^*)^2}{\alpha^*}+\hat\beta_2\beta^* ,\; 0\right),\\[10pt]
\alpha_2= b_1\alpha^*\neq 0,\quad
\alpha_3 =\begin{pmatrix}
 b_1\hat\alpha_2  \\[10pt]
 b_1\hat\alpha_4- b_4\hat\alpha_1
\end{pmatrix},
\quad
\beta_3 = \begin{pmatrix}
\hat\beta_2  -\hat\alpha_2\dfrac{\beta^*}{\alpha^*} & 0\\[10pt]
 \hat\beta_4 -\dfrac{b_4}{b_1}\hat\beta_1 -(\hat\alpha_4-\dfrac{b_4}{b_1}\hat\alpha_1)\dfrac{\beta^*}{\alpha^*} & 0
\end{pmatrix},\\[20pt]
\alpha_4 = b_1\hat\alpha_5,
\quad
\beta_4 =\begin{pmatrix}
 \hat\beta_5 -\hat\alpha_5\dfrac{\beta^*}{\alpha^*} 
& 0\\ 0 & 0
\end{pmatrix}.
\end{array}
\end{equation}
The functions $h$ are given by
\begin{equation}\label{eq:maps:Tk_h}
\begin{array}{l}
h_1 = \dfrac{1}{b_1} \hat h_1 + \dfrac{\beta^*}{\alpha^* b_1} \hat h_2,
\qquad
h_2=\hat h_3,\qquad
h_3 =(\hat h_2,\hat h_4)^T,
\qquad
h_4 = \hat h_5,
\end{array}
\end{equation}
and, by (\ref{hathnew}), they satisfy 
\begin{equation}\label{eq:maps:Tk_nonlinear}
\begin{aligned}
& h_{1,3,4}
=O(\hat \lambda^k)(Z,W,\bar V)+O(\lambda^k)Y+O(Y^2),\\
&h_2=O(\hat\lambda^k)(Z,W,\bar V)+O(\lambda^k)(Z,W)Y+O(\lambda^k)Y^2+O(Y^3),
\\ \\
&\dfrac{\partial  h_{1,3,4}}{\partial \mu}=O(\gamma^k\hat\gamma^{-k})+O(Y^2),\quad
\dfrac{\partial  h_2}{\partial \mu}
=O(\gamma^k\hat\gamma^{-k})+o(Y^2),
\\ 
&\dfrac{\partial h_{1,3,4}}{\partial \omega}
=O(\lambda^k\gamma^k\hat\gamma^{-k})+O(Y^2) + O(k\lambda^k)Y, \\
&\dfrac{\partial h_2}{\partial \omega}=\gamma^k\hat\gamma^{-k} 
O(|\mu| + Y^2 + \lambda^k)+ o(Y^2)+o(\lambda^k)Y,
\\ \\
& \dfrac{\partial  h_{i}}{\partial(Z,W)} 
=O(\hat\lambda^k)+O(\lambda^k)Y,\quad  i=1,2,3,4\\
& 
\dfrac{\partial  h_{1,3,4}}{\partial Y}=O(\lambda^k)+O(Y),
\\
&\dfrac{\partial  h_2}{\partial Y} 
=O(\hat\lambda^k) + O(\lambda^k)(Z,W)+ O(\gamma^k\hat\gamma^{-k})Y+O(Y^2).
\end{aligned}
\end{equation}
We also have
$$\dfrac{\partial  h_i}{\partial \bar V}=O(\hat\gamma^{-k}), \quad i=1,2,3,4.$$
However, we will only use the weaker estimates
\begin{equation}\label{eq:h_df}
\dfrac{\partial  h_{1,3,4}}{\partial \bar V}=O(\gamma^{-k})
\quad\mbox{and}\quad
\dfrac{\partial  h_2}{\partial \bar V}=O(\gamma^{-k})Y+O(\hat\gamma^{-k}), 
\end{equation} 
for conformity with the bi-focus case (see Section \ref{sec:df}).\\

Recall that the domain of $T_k$ does not shrink to 0 as $k\to\infty$. We will (by possibly diminishing the value of 
$\delta$) consider the restriction of $T_k$ to the cube
\begin{equation}\label{eq:domain}
\Pi:=[-\delta,\delta]\times [-\delta,\delta]\times [-\delta,\delta]^{d^s-1}\times [-\delta,\delta]^{d^u-1} \subset \hat\Pi,
\end{equation}
in the coordinates $(Z,Y,W,V)$.

\begin{conven*}
Throughout the rest of the paper, we assume $\delta$ small. By saying that  a result holds  for all sufficiently small $\delta$ and sufficiently large $k$, we mean that there exists $\delta_0>0$ such that for each fixed $\delta<\delta_0$ the result holds for all $k>k_\delta$ for some $k_\delta$ depending on $\delta$. Sometimes for brevity we may only say `for all  sufficiently large $k$'.
\end{conven*}

\subsection{Invariant cone fields}
\begin{lem}\label{lem:cones} There exists a constant $K>0$
such that for all sufficiently small $\delta>0$ and for all sufficiently large $k$ such that the coefficient $\alpha^*$ defined in (\ref{eq:maps:alphabeta}) is bounded away from zero, the cone fields in $\Pi$ 
\begin{align}
\mathcal{C}^{cs}&=\{( z, y,  w,  v): \| v\| < K\delta(| z|+| y|+\|w\|) \},\label{eq:cone_cs}\\
\mathcal{C}^{ss}&=\{( z, y,  w, v) : |z|+|y|+\| v\| < K\delta\|w\|
\},  \label{eq:cone_ss}
\end{align}
are backward-invariant in the sense that if a point $\bar M\in \Pi$ has its pre-image $M=T^{-1}_k(\bar M)$ in $\Pi$, then the cone at $\bar M$ is mapped into the cone at $M$ by $\D T^{-1}_k$; and the cone field
\begin{align}
\mathcal{C}^{cu}=\{(z,y,w,v): \|w\|< K(\delta|z|+\gamma^{-k}|y|+\gamma^{-k}\|v\|)   
\},\label{eq:cone_cu}
\end{align} 
is forward-invariant in the sense that if a point $M\in \Pi$ has its image $\bar M=T_k(M)$ in $\Pi$, then the cone at $M$ is mapped into the cone at $\bar M$ by $\D T_k$. There also exists a  forward-invariant cone field
\begin{equation}\label{eq:coneuu_df}
\mathcal{C}^{uu}= \left\{(z,y,w,v): 
|z| + \|w\|<K\gamma^{-k}\|v\|,\;
|y|< K\delta^2\|v\|
\right\}
\end{equation}
defined on a smaller region
\begin{equation}\label{eq:PiY}
\Pi_Y:=[-\delta,\delta]\times [-2\delta^2,2\delta^2]\times [-\delta,\delta]^{d^s-1}\times [-\delta,\delta]^{d^u-1} \subset \Pi.
\end{equation}
\end{lem}

\begin{proof}
Since $(Z,Y,W,\bar V)=O(\delta)$, it follows from  \eqref{eq:maps:Tk}, \eqref{eq:maps:Tk_nonlinear}, \eqref{eq:h_df}  that the differential $\D T_k: (z,y,w,v)\mapsto (\bar z,\bar y,\bar w,\bar v)$
at a point $(Z,Y,W,V)$ satisfies, for every given small $\delta$ and sufficiently large $k$, 
\begin{align}
\bar z&=O(\lambda^k)z+(1+O(\delta))y+O(\lambda^k)w+O( \gamma^{-k})\bar v,\label{eq:cone:1}\\[5pt]
\bar y&=\alpha_2 \lambda^k\gamma^k(1+O(\delta) )z+\gamma^kO(\delta)y +\lambda^k\gamma^k O(\delta^2) w+O(\delta)\bar v,\label{eq:cone:2}\\[5pt]
\bar w&=O(\lambda^k)z+O(\delta)y +O(\lambda^k)w+O( \gamma^{-k})\bar v,\label{eq:cone:3}\\[5pt]
v&=O(\lambda^k)z+O(\delta)y +O(\lambda^k)w+O( \gamma^{-k})\bar v.\label{eq:cone:4}
\end{align}
Note that the coefficient of $\bar v$ in \eqref{eq:cone:2} is $O(Y)+O(\gamma^k\hat\gamma^{-k})$ by \eqref{eq:h_df}. Since
$|Y|\leq \delta$ in $\Pi$ and $|\gamma|<\hat\gamma$, this coefficient indeed can be replaced by $O(\delta)$ for large enough $k$.

\noindent{\bf (1) The cone field $\mathcal{C}^{ss}$}. Choose some $K>0$. 
Our goal is to show that for a proper choice of $K$,
for any vector $(\bar z,\bar y,\bar w,\bar v)\in\mathcal{C}^{ss}$ its preimage 
$(z,y,w,v)=DT_k^{-1} (\bar z,\bar y,\bar w,\bar v)$ also lies in $\mathcal{C}^{ss}$.

Since, $\|\bar v\| < K \delta \|\bar w\|$ by \eqref{eq:cone_ss},
it follows from \eqref{eq:cone:3} that 
\begin{equation}\label{ss:barw}
\|\bar w\|<C_1\lambda^k(|z|+\|w\|)+C_1\delta |y|
\end{equation}
for some constant $C_1>0$ independent of $K$ (for  small $\delta$).
Since $\max\{|\bar z|,|\bar y|,\|\bar v\|\}< K \delta \|\bar w\|$, and $b_1\alpha^*$ is bounded away from zero, 
it follows from (\ref{ss:barw}) that, upon expressing $y$ as a function of $\bar z,z,w,\bar v$ from \eqref{eq:cone:1},   
\begin{equation}\label{yeqss}
|y|\le C_2\lambda^k(|z|+\|w\|)
\end{equation}
for some constant $C_2$  independent of  $K$. Substituting this into \eqref{eq:cone:1}, we find that 
\begin{equation}\label{zeqss}
|z|\le C_3\delta \|w\|
\end{equation}
for some $C_3$  independent of $K$. The substitution of these inequalities into  \eqref{eq:cone:4} gives
$$\|v\|\le C_4\lambda^k\|w\|$$
for some $C_4$  independent of $K$.  Thus, we have $(z,y,w,v)\in \mathcal{C}^{ss}$ if we take $K>\max\{C_2,C_3,C_4\}$ (with $\delta$ sufficiently small).

\noindent{\bf (2) The cone field $\mathcal{C}^{cs}$}. Take any vector $(\bar z,\bar y,\bar w,\bar v)\in\mathcal{C}^{cs}$ defined
by some constant $K>0$. Let $(z,y,w,v)=DT_k^{-1} (\bar z,\bar y,\bar w,\bar v)$. If $\delta$ is small enough, then
$$\|\bar v\|< |\bar z|+|\bar y|+\|\bar w\|,$$
and it follows from \eqref{eq:cone:1} and \eqref{eq:cone:3} that
$$
|\bar z|+\|\bar w\|=O(\lambda^k)(|z|+\|w\|) + O(y) + O(\gamma^{-k})(|\bar z|+\|\bar w\|)
+O(\gamma^{-k}) |\bar y|,$$
which implies
$$|\bar z|+\|\bar w\| \leq C_5 (\lambda^k (|z|+\|w\|) + |y| +|\gamma|^{-k} |\bar y|),$$
for a constant $C_5$ independent of $K$. The substitution into \eqref{eq:cone:2} gives
$$
|\bar y| \leq C_6 (|\lambda\gamma|^k(|z|+\|w\|) + |\gamma|^k \delta |y|),$$
so
$$|\gamma|^{-k} \|\bar v\|\leq  C_7 (\lambda^k(|z|+\|w\|) + \delta |y|)$$
if $k$ is large enough; the constants $C_{6,7}$ are independent of the choice of $K$.
By \eqref{eq:cone:4}, we finally get
$$\|v\| \leq  C_8 (\lambda^k(|z|+\|w\|) + \delta |y|),$$
implying $(z,y,w,v)\in\mathcal{C}^{cs}$ for $K > C_8$, as desired.

\noindent{\bf (3) The cone field $\mathcal{C}^{cu}$}. Take any vector $(z, y, w, v)\in\mathcal{C}^{cu}$ for some constant $K>0$.
If $\delta$ is small enough, then
\begin{equation*}
\|w\|< |z|+|y|+\|v\|,
\end{equation*}
This together with \eqref{eq:cone:4} yields
\begin{equation}\label{wcu}
w= O(|z|+|y|+\gamma^{-k}\|\bar v\|),
\end{equation}
which is $K$-independent.  Substituting \eqref{wcu} into
 equations \eqref{eq:cone:1} -- \eqref{eq:cone:3} gives
\begin{align}
|y|&=O(\lambda^k)z+O(1)|\bar z|+O(\gamma^{-k})\|\bar v\|, \nonumber\\[5pt]
\lambda^k |z|&=O(\gamma^{-k})|\bar y|+O(\delta)|y|+O(\gamma^{-k}\delta)\| \bar v\|, \label{cu:6}\\[5pt]
\|\bar w\|&=O(\lambda^k)|z| +O(\delta)| y| + O(\gamma^{-k})\| \bar v\|. \nonumber
\end{align}
Using \eqref{cu:6}, one finds that
$$y=O(\gamma^{-k})|\bar y|+O(1)|\bar z|+O(\gamma^{-k})\| \bar v\|,$$
and
$$\|\bar w\|\leq C_9 (|\gamma|^{-k} |\bar y| + \delta |\bar z| + |\gamma|^{-k} \|\bar v\|)$$
for some constant $C_9$, so $(\bar z, \bar y, \bar w, \bar v)\in\mathcal{C}^{cu}$ for $K>C_9$, as desired.

\noindent{\bf (4) The cone field $\mathcal{C}^{uu}$}. Take any vector $(z, y, w, v)\in\mathcal{C}^{uu}$ for some constant $K$. By \eqref{eq:coneuu_df}, if $\delta$ is small enough and $k$ are large enough,
$$|z|+\|w\|\leq \|v\|, \qquad |y|\leq \delta \|v\|.$$
By \eqref{eq:cone:4}, this implies
\begin{equation}\label{uu:5}
\|v\|\leq C |\gamma|^{-k}\|\bar v\|,
\end{equation}
for a constant $C>0$ which is independent of $K$ for small $\delta$ and large enough $k$. It also follows that
\begin{equation}\label{uu:5bis}
|z|+\|w\|=O(\gamma^{-k}) \bar v, \qquad y = \gamma^{-k} O(\delta) \bar v.
\end{equation}
Note that equation \eqref{eq:cone:3} can be written as
$$\bar y= \gamma^k [O(\lambda^k)z+O(\delta)y + O(\lambda^k) w] + (O(Y)+ O(\gamma^k\hat\gamma^{-k}))\bar v$$
(see \eqref{eq:h_df}). In $\Pi_Y$, we have $|Y|\le \delta^2$ by \eqref{eq:PiY}, so
\begin{equation}\label{eq:cone:2bis}
\bar y= \gamma^k [O(\lambda^k)z+O(\delta)y + O(\lambda^k) w] + O(\delta^2)\bar v.
\end{equation}
Now, substituting \eqref{uu:5bis} into \eqref{eq:cone:1}, \eqref{eq:cone:3} and \eqref{eq:cone:2bis} yields
$$\bar z=O(\gamma^{-k}) \bar v, \qquad \bar w=O(\gamma^{-k}) \bar v,
\qquad\bar y=O(\delta^2) \bar v,$$
i.e.,  as desired, $(\bar z, \bar y, \bar w, \bar v) \in \mathcal{C}^{uu}$, if $K$ was chosen large enough.
\end{proof}

\begin{rem}\label{rem:cone_contraction} By \eqref{uu:5}, vectors in $\mathcal{C}^{uu}$ are uniformly expanded by $\D T_k$
with the expansion factor of order $\gamma^{k}$ at least. Also,
vectors in $\mathcal{C}^{ss}$ are uniformly contracted by $\D T_k$: substituting \eqref{yeqss},\eqref{zeqss} into \eqref{ss:barw} gives
\begin{equation}\label{barw}
\bar w=O(\lambda^k) w.
\end{equation}
This will be used in the proof of Lemma \ref{lem:index2}.
\end{rem}

\section{Periodic point $Q$ with a higher index and its local invariant manifolds}\label{sec:QQ}
In this section, we find, by unfolding the tangency to $O$, a fixed point $Q$ of the first-return map $T_k$, along with formulas for its local invariant manifolds. In the next section, we describe how a further change of parameters creates
 a heterodimensional cyle involving $O$ and $Q$.


Throughout this section we assume that $k$ is taken such that $\alpha^*$ in \eqref{eq:maps:Tk} is bounded away from zero.
We also use that $|\lambda\gamma|>1$, so $|\gamma|^{-k} \ll \hat\lambda^k \ll \lambda^k$ for large $k$.
\subsection{Existence of periodic points with a higher index} \label{sec:Q}

Recall that our system in at least $C^3$, so even though we loose two derivatives with respect to parameters when introducing the coordinates of   \citep[Lemma 6]{GST08} in Section~\ref{sec:maps}, all coefficients in our formulas are still smooth functions of $\eps=(\mu,\omega)$. We will not indicate the dependence of $\eps$ unless it is necessary. In the statement of the lemma below all the coefficients in the right-hand sides are taken at $\mu=0$, so they are smooth functions of $\omega$ (except for the constant $\hat\lambda$). Recall that $\alpha^*$ is the coefficient defined in \eqref{eq:maps:alphabeta}.

\begin{lem}\label{lem:index2}
 For all sufficiently large $k$ such that $\alpha^*$ is bounded away from zero, there exist smooth functions 
 $\rho_k(\omega)$ and $\hat\rho_k(\omega,t)$ $(t\in \mathbb R)$ with bounded first derivatives uniformly for all $k$ and bounded $t$, such that, at
\begin{equation}\label{eq:mu_Q}
\begin{aligned}
\mu=\mu_k(\omega,t):=-\lambda^k(\alpha^*x_1^++\beta^*x_2^+) +\hat\lambda^k\rho_k(\omega)+\dfrac{(b_1\alpha^*)^2}{2d}\lambda^{2k} (t-\dfrac{t^2}{2})+
\hat\lambda^k\lambda^{k} \hat \rho_k(\omega,t),
\end{aligned}
\end{equation}
the map $T_k$ has a fixed point $Q=(Z_Q,Y_Q,W_Q,V_Q)$ with
\begin{equation}\label{eq:Q_coor}
Z_Q=-\dfrac{b_1\alpha^*}{2d}\lambda^k t+O(\hat\lambda^k),\quad Y_Q=-\dfrac{b_1\alpha^*}{2d}\lambda^k t,\quad (W_Q,V_Q)=O(\hat\lambda^k),
\end{equation}
where the $O(\hat \lambda^k)$ terms have the derivatives with respect to $t$ and $\omega$ of order
$O(\hat \lambda^k)$. The point $Q$ has index $d^u+1$ for $t\in (t_-,t_+)$ where
\begin{equation}\label{tpm}
t_+ = 1 +O(\hat\lambda^k\lambda^{-k}), \qquad t_- = -1 +O(\hat\lambda^k\lambda^{-k});
\end{equation}
at $t=t_+$ the point $Q$ has a multiplier $\nu=1$ and at $t=t_-$ it has a multiplier $\nu=-1$.
\end{lem}

\begin{rem}
Here we introduce an artificial parameter $t$ (the rescaled value of the $Y$-coordinate of the fixed point $Q$) as the most convenient way for us to control the position and the index of $Q$.
\end{rem}

\begin{proof} By \eqref{eq:maps:Tk} and \eqref{eq:maps:Tk_const},
the map $T_k$ has a fixed point $(Z_Q,Y_Q,W_Q,V_Q)$ with $Y_Q$ given by \eqref{eq:Q_coor} when
\begin{equation}\label{eq:maps:Tkfixed}
\begin{aligned}
Z_Q &= \lambda^k \alpha_1 Z_Q-\dfrac{b_1\alpha^*}{2d}\lambda^k t+\lambda^k\beta_1W_Q+ h_1(Z_Q,-\dfrac{b_1\alpha^*}{2d}\lambda^k t,W_Q,V_Q),\\
W_Q &=  \lambda^k \alpha_3 Z_Q+\lambda^k\beta_3W_Q+ h_3(Z_Q,-\dfrac{b_1\alpha^*}{2d}\lambda^k t,W_Q,V_Q),\\
V_Q&= \lambda^k \alpha_4 Z_Q+\lambda^k\beta_4W_Q+ h_4(Z_Q,-\dfrac{b_1\alpha^*}{2d}\lambda^k t,W_Q,V_Q),
\end{aligned}
\end{equation} 
and
\begin{equation}\label{muyqfixed} 
\begin{aligned}
\mu = -\lambda^k(\alpha^* x^+_1 + \beta^* x^+_2) - b_1\alpha^* \lambda^k Z_Q  - \dfrac{(b_1\alpha^*)^2}{4d} \lambda^{2k} t^2 - h_2(Z_Q,-\dfrac{b_1\alpha^*}{2d}\lambda^kt,W_Q,V_Q)\\
\qquad \qquad \qquad +O(\hat\lambda^k)+O(\lambda^k\gamma^{-k}),
\end{aligned}
\end{equation}
where the term $O(\hat\lambda^k)$ depends only on $\omega$ and comes from the terms in \eqref{eq:maps:Tk_const} other than $(\alpha^* x^+_1 + \beta^* x^+_2)$, and the term $O(\lambda^k\gamma^{-k})$ comes from $Y_Q$ depending on both $\omega$ and $t$: recall that $\gamma^{-k}\ll \hat\lambda^k\ll \lambda^k$. By \eqref{eq:maps:Tk_nonlinear} and \eqref{eq:h_df}, the functions $h$ satisfy the following estimates for $(Z,Y,W,V,\mu)=O(\lambda^k)$:
\begin{equation}\label{hagain}
\begin{aligned}
&\dfrac{\partial  h_{i}}{\partial \mu}=O(\gamma^k\hat\gamma^{-k}),
\quad
\dfrac{\partial h_{i}}{\partial \omega}
=O(\lambda^k\gamma^k\hat\gamma^{-k}), 
\quad
 \dfrac{\partial  h_{i}}{\partial(Z,W,\bar V)} 
=O(\hat\lambda^k),\quad i=1,2,3,4,\\
&\dfrac{\partial  h_{1,3,4}}{\partial Y}=O(\lambda^k),\quad
\dfrac{\partial  h_2}{\partial Y} 
=O(\hat\lambda^k), \qquad
h_{1,3,4}=O(\lambda^{2k}), \quad h_2=O(\hat\lambda^k\lambda^k).
\end{aligned}
\end{equation}
Since the right-hand sides of \eqref{eq:maps:Tkfixed} and  \eqref{muyqfixed} are small with first derivatives, we can resolve these equations with respect to $(Z_Q,W_Q,V_Q,\mu)$ as functions of $t$ and $\omega$ by the Implicit Function Theorem.
Using estimates \eqref{hagain}, we obtain formulas \eqref{eq:mu_Q} and \eqref{eq:Q_coor} (like we did it before, the coefficient $k$ emerging when differentiating $\lambda^k$ and $\gamma^k$ with respect to the parameters, is eliminated by slightly increasing $\hat\lambda$).

Thus, to prove the lemma, it remains to check that $Q$ has index $d^u+1$ when $t=\left|\dfrac{2d Y_Q}{b_1\alpha^*}\right|\lambda^{-k}$ varies between $t_-$ and $t_+$ satisfying \eqref{tpm}.
By Lemma \ref{lem:cones}, in the tangent space at any periodic point in $\Pi$, there exist a forward-invariant subspace $E^{cu}\subset \mathcal{C}^{cu}$ and a backward-invariant subspace $E^{cs}\subset \mathcal{C}^{cs}$. Let us write $E=E^{cs}\cap E^{cu}$. Vectors in $E$ have the form $(z,y,S_1(z,y),S_2(z,y))$, where $S_{1,2}$ are linear maps satisfying that for any given $K^*>0$ we have $\|S_{1,2}(z,y)\|<K^*(|z|+|y|)$ for all sufficiently small $\delta$ and all sufficiently large $k$. To see this, one notes that for a vector $(z,y,w,v)\in E$, the $w$-component satisfies
\begin{align*}
\|w\|&<K\delta(|z|+|y|+\|v\|)\qquad\qquad &\mbox{by \eqref{eq:cone_cu}}\\
&<K\delta(|z|+|y|+K\delta(|z|+|y|+\|w\|)), &\mbox{by \eqref{eq:cone_cs}}
\end{align*}
or, equivalently,
$$\|w\|<\dfrac{K\delta(1+K\delta)}{1-(K\delta)^2}(|z|+|y|)<K^*(|z|+|y|),$$
since $\delta$ can be taken sufficiently small. The $v$-component can be checked in the same way.

Consider $DT_k|_E$ as the linear transformation of $\mathbb{R}^2$ defined in the following way:
\[DT_k|_E(z,y)=(\bar z,\bar y) \] if and only if
\[DT_k(z,y,S_1(z,y),S_2(z,y))=(\bar z,\bar y, S_1( \bar z,\bar  y),S_2(\bar z,\bar y)).\]

Differentiating the first two equations in \eqref{eq:maps:Tk} with using $\alpha_2 = b_1\alpha^*$ in \eqref{eq:maps:Tk_coef} and the estimates in \eqref{eq:h_df} and \eqref{eq:maps:Tk_nonlinear}, yields
\begin{align*}
\bar z=&O(\lambda^k)z+y+O(\hat\lambda^k+\lambda^k |Y_Q|)(z+S_1(z,y)) +
O(\lambda^k+|Y_Q|)y
+ O(\gamma^{-k})S_2(\bar z,\bar y),\\
\bar y =&b_1\alpha^*\lambda^k\gamma^k z+ 2Y_Q\gamma^k d y
+
\gamma^k(O(\hat\lambda^k+ \lambda^k O(|Z_Q|+\|W_Q\|)+\gamma^k\hat\gamma^{-k}|Y_Q|+|Y_Q|^2) y \\
&\qquad\qquad\qquad + O(\hat\lambda^k+\lambda^k |Y_Q|)(z+S_1(z,y)) + O(\gamma^{-k}|Y_Q|+\hat\gamma^{-k}) S_2(\bar z,\bar y).
\end{align*}
Since $S_i$ are bounded and linear, using the Implicit Function Theorem, we obtain  the following formula for $DT_k|_E$ at the point $Q$:
\begin{equation}\label{dtke}
\begin{aligned}
\bar z=&O(\lambda^k)z+(1+O(\lambda^k) +O(Y_Q))y,\\
\bar y=&b_1\alpha^*\lambda^k\gamma^k(1+O(\hat\lambda^k\lambda^{-k}) +O(Y_Q))z\\
&+ \gamma^k (2d Y_Q+ \lambda^k O(|Z_Q|+\|W_Q\|)+ 
O(\gamma^k\hat\gamma^{-k})Y_Q+ O(Y_Q^2)+O(\hat\lambda^k))y.
\end{aligned}
\end{equation}
It follows that the eigenvalues of $DT_k|_E$, denoted by $\nu_1$ and $\nu_2$, satisfy 
\begin{align}
\nu_1+\nu_2&=- b_1\alpha^*\lambda^k\gamma^k t +O(\hat\lambda^k\gamma^k),\label{eigen:1}\\
\nu_1\nu_2&=-b_1\alpha^*\lambda^k\gamma^k +O(\hat\lambda^k\gamma^k),\label{eigen:2}
\end{align}
where we used \eqref{eq:Q_coor} and \eqref{eq:lambdagamma} to sort the $O(\cdot)$ terms. 

The matrix $DT_k(Q)|_E$ has two eigenvalues outside the unit circle if and only if $|\nu_1\nu_2|>1$ and 
$|\nu_1+\nu_2|< |\nu_1\nu_2+1|$ 
(see e.g. \citep[Section 2.3.1]{Liphd}); the boundary $\nu_1+\nu_2 = \nu_1\nu_2 +1$ corresponds to
one of the multipliers equal to $1$, and  $\nu_1+\nu_2 = - (\nu_1\nu_2 +1)$ corresponds to a multiplier equal to $(-1)$. In our case, the inequality $|\nu_1\nu_2|>1$ follows immediately from \eqref{eigen:2} for sufficiently large $k$ due to the assumption $|\lambda\gamma|>1$, and the equality 
$\nu_1+\nu_2 = \pm (\nu_1\nu_2+1)$ immediately gives \eqref{tpm} by \eqref{eigen:1}, \eqref{eigen:2}.
 
This implies that $Q$ has index $d^u+1$ for $t\in (t_-,t_+)$, since the remaining eigenvalues of $DT_k(Q)$ correspond to eigenvectors belonging to
either $\mathcal{C}^{ss}$ or $\mathcal{C}^{uu}$ and, by Remark \ref{rem:cone_contraction}, these are $d^s-1$ eigenvalues inside the unit circle and $d^u-1$ eigenvalues outside the unit circle.
\end{proof}

Note also that the above lemma implies that
\begin{equation}\label{eq:t_deri}
\begin{aligned}
&\dfrac{\partial \mu_k}{\partial t}=\dfrac{(b_1\alpha^*)^2}{2d}\lambda^{2k}(1-t)+O(\hat\lambda^k\lambda^{k}),
\quad
\dfrac{\partial Z_Q}{\partial t} = -\dfrac{b_1\alpha^*}{2d}\lambda^{k} + O(\hat\lambda^k),
\\
&\dfrac{\partial Y_Q}{\partial t} = -\dfrac{b_1\alpha^*}{2d}\lambda^{k},\quad
\dfrac{\partial (W_Q,V_Q)}{\partial t} =O(\hat\lambda^{k}),
\\
&\dfrac{\partial (\mu_k,Z_Q,Y_Q)}{\partial \omega}=O(k\lambda^{k}),\quad
\dfrac{\partial (W_Q,V_Q)}{\partial \omega} =O(\hat\lambda^{k})
\end{aligned}
\end{equation}

\begin{rem}\label{saddle1} 
It follows from \eqref{eigen:1} and \eqref{eigen:2} that among the multipliers of $Q$ (the eigenvalues of $DT_k(Q)$) those corresponding to $\mathcal{C}^{uu}$ have larger absolute values than $\nu_{1,2}$ (we have from \eqref{eigen:1} and \eqref{eigen:2} that $|\nu_{1,2}|$ is estimated as $O(\lambda^k\gamma^k)$, whereas the multipliers corresponding to $\mathcal{C}^{uu}$ are of order $O(\gamma^{k})$ at least). It also follows, that if $t$ is bounded away from zero, then $|\nu_1|\neq |\nu_2|$, so the nearest to the unit circle unstable multiplier of $Q$
is real. We denote this multiplier as $\nu$ and the other one as 
$\hat\nu$. We find from \eqref{eigen:1} and \eqref{eigen:2} that
\begin{equation}\label{eigen:3}
\nu = \frac{1}{t +O(\hat\lambda^k\lambda^{-k})}, \qquad \hat\nu =- b_1\alpha^*\lambda^k\gamma^k t + O(\hat\lambda^k\gamma^k).
\end{equation}
As it is seen from \eqref{dtke}, the corresponding eigenvectors of $DT_k|_E$ at the point $Q$ are given by
\begin{equation}\label{eigen:4}
(z,y) = ( 1 +O(\hat\lambda^k\lambda^{-k}),\nu) \qquad \mbox{and} \qquad  
(z,y) = ( 1 +O(\hat\lambda^k\lambda^{-k}),\hat\nu).
\end{equation}
This formula will be used in Section~\ref{sec:Wse}.
\end{rem}

\subsection{Local stable invariant manifold $W^s_{loc}(Q)$}
In this and the next sections, we derive equations for the local stable and unstable manifolds of the point $Q$ found at $\mu=\mu_k$ by Lemma \ref{lem:index2}. These formulas will be used for finding intersections of $W^s(Q)$ with
$W^u(O)$ and of $W^u(Q)$ with $W^s(O)$, which thus lead to a heterodimensional cycle involving $O$ and $Q$. This will be the first step for obtaining the robust heterodimensional dynamics in Theorem \ref{thm:sf}. 

Denote by $W^s_{loc}(Q)$ the connected component of $W^s(Q)\cap\Pi$ which contains $Q$.

\begin{lem}\label{lem:WsQ}
The local stable manifold $W^s_{loc}(Q)$ is given by 
\begin{equation}\label{eq:WsQ:0.1}
(Z,Y,V)=(Z_Q,Y_Q,V_Q) + \phi^s(W;\omega,t),
\end{equation}
where $\phi^s=(\phi^s_z,\phi^s_y,\phi^s_v)$ is defined for all $W\in [-\delta,\delta]^{d_s-1}$ and all 
$t\in(t_-,t_+)$, is at least $C^{r-2}$, and satisfies
$\phi^s(W_Q)=0$ and, for some $C>0$,
\begin{equation}\label{eq:WsQ:0.2}
\begin{aligned}
&|\phi^s_z| + \left\|\dfrac{\partial \phi^s_z}{\partial W}\right\|\leq C \lambda^{-k}\hat{\lambda}^{k},
\qquad
\|\phi^s_y,\phi^s_v\|+\left\|\dfrac{\partial (\phi^s_y,\phi^s_v)}{\partial W}\right\|\leq C \lambda^k,
\\
&\left\|\dfrac{\partial (\phi^s_z,\phi^s_y)}{\partial \omega}\right\|\leq C \lambda^{-k}\hat\lambda^k,
\qquad
\left\|\dfrac{\partial \phi^s_v}{\partial \omega}\right\|\leq Ck\lambda^{k},\\
&\left\|\dfrac{\partial (\phi^s_z,\phi^s_y)}{\partial t}\right\|\leq C (\lambda^k\|W\|+\hat\lambda^k),
\qquad
\left\|\dfrac{\partial \phi^s_v}{\partial t}\right\|\leq C\lambda^{2k}.
\end{aligned}
\end{equation}
\end{lem}

\begin{proof}
Let us move the origin to $Q$ by applying
\begin{equation}\label{newcoor}
Z^{new}=Z-Z_Q,\quad Y^{new}=Y-Y_Q,\quad W^{new}=W-W_Q,\quad V^{new}=V-V_Q,
\end{equation}
which allows us to rewrite $T_k$ (see \eqref{eq:maps:Tk}) at 
$\mu=\mu_k$ as
\begin{equation}\label{locmapq}
\begin{aligned}
\bar Z^{new} &= \lambda^k \alpha_1 Z^{new}+Y^{new}+\lambda^k\beta_1W^{new}+ h_1(M+Q)-h_1(Q),\\
\bar Y^{new} &=  b_1\alpha^*\lambda^k\gamma^k Z^{new}  +2d\gamma^k Y_QY^{new} +d\gamma^k(Y^{new})^2+ \gamma^k (h_2(M+Q)-h_2(Q)),\\
\bar W ^{new}&=  \lambda^k \alpha_3 Z^{new}+\lambda^k\beta_3W^{new}+ h_3(M+Q)-h_3(Q),\\
V^{new}&= \lambda^k \alpha_4 Z^{new}+\lambda^k\beta_4W^{new}+ h_4(M+Q)-h_4(Q),
\end{aligned}
\end{equation}
where we denote $M=(Z^{new},Y^{new},W^{new},\bar V^{new})$; recall that the functions $h$ satisfy the estimates \eqref{eq:maps:Tk_nonlinear} and \eqref{eq:h_df}. We rewrite these formulas as
\begin{equation}\label{eq:preWs}
\begin{aligned}
 Y^{new}&=  \bar Z^{new}-  \alpha_1 \lambda^k Z^{new}- \beta_1 \lambda^k W^{new}-\tilde h_1,\\
Z^{new} &= \frac{1}{b_1\alpha^*} (\lambda\gamma)^{-k} \bar Y^{new}+  tY^{new} - \frac{d}{b_1\alpha^*}\lambda^{-k} (Y^{new})^2- \frac{1}{b_1\alpha^*} \lambda^{-k} \tilde h_2,\\
\bar W^{new}&=  \lambda^{k} \alpha_3 Z^{new}+\lambda^k\beta_3W^{new}+ \tilde h_3,\\
V^{new}&= \lambda^k \alpha_4 Z^{new}+\lambda^k\beta_4W^{new}+ \tilde h_4,
\end{aligned}
\end{equation}
where we denoted $\tilde h_i(M) = h_i(M+Q)-h_i(Q)$, $i=1,2,3,4$; we will further drop the ``new'' label. By
\eqref{eq:maps:Tk_nonlinear}, \eqref{eq:h_df}, \eqref{eq:t_deri}, \eqref{eq:lambdagamma}, we have
that given any constant $C>0$, if $|Y|\leq C\lambda^k$, then
\begin{equation}\label{tldha}
\begin{aligned}
& \tilde h_{i}
=O(\hat \lambda^k),\quad \dfrac{\partial \tilde  h_{i}}{\partial(Z,W,\bar V,\omega)} 
=O(\hat\lambda^k),\qquad  i=1,2,3,4,
\\
&\dfrac{\partial  \tilde h_{1,3,4}}{\partial Y}=O(\lambda^k),\quad
\dfrac{\partial  \tilde h_{1,3,4}}{\partial t}=O(\lambda^{2k}),\\
&
\dfrac{\partial  \tilde h_2}{\partial Y} 
=O(\hat\lambda^k)+\lambda^k O(W),\quad
\dfrac{\partial \tilde h_{2}}{\partial t}=
\lambda^k (O(\hat\lambda^k)+O(W)).
\end{aligned}
\end{equation}
Note that the bounds in the $O(\cdot)$ terms here are independent of the choice of $C$, provided   $k$ is large enough.

As shown in  the last paragraph of the proof of Lemma~\ref{lem:index2}, 
 the tangents to the stable manifold $W^s_{loc}(Q)$ belong to
 the backward-invariant cone field $\mathcal{C}^{ss}$ given by \eqref{eq:cone_ss}. At the same time, 
 since vectors in $\mathcal{C}^{ss}$  are uniformly expanded by $DT^{-1}_k$, the tangents to $W^u_{loc}(Q)$ do not belong to $\mathcal{C}^{ss}$, and hence every manifold near $Q$ whose tangents belong to $\mathcal{C}^{ss}$ is transverse to $W^u_{loc}(Q)$. In particular, $\{(Z,Y,V)=0\}$ intersect $W^u_{loc}(Q)$ transversely. Therefore, according to the Palis lambda-lemma, we can find $W^s_{loc}(Q)$
as the limit of the backward iterations of $\{(Z,Y,V)=0\}$  by (\ref{locmapq}).

Thus, to prove the lemma, it suffices to show that the preimage of the graph of a function $(Z,Y,V)=\phi(W;\omega,t)$ satisfying \eqref{eq:WsQ:0.2} is again a graph of a function satisfying the same estimates. 
To this aim, one substitutes  
$(\bar Z,\bar Y,\bar V)=\phi(\bar W;\omega,t)$
in the right-hand sides of \eqref{eq:preWs}, with $\phi$
satisfying \eqref{eq:WsQ:0.2}. It is obvious from \eqref{tldha} that
$\bar W$ is uniquely determined as a function of 
$(W, Z, Y, V, \omega,t)$ from the third equation of \eqref{eq:preWs}. In particular, since $Z$ and $W$ are bounded, 
$\bar W=O(\lambda^k)$. This implies that 
$\bar Z= O(||\partial_W \phi_z\|) \|\bar W\|= O(\hat\lambda^k)$.
Then, using \eqref{tldha}, one can see that given any $W$, the right-hand sides of the first, second and fourth equations define a contracting map of the region
$$\|Z\| \leq C \lambda^{-k} \hat\lambda^k,\quad \|Y,V\|\leq C\lambda^k,$$
for any sufficiently large $C$. Therefore, these equations define $(Z,Y,V)$
as a function of $W$, and this function takes values in this region.

Thus, the preimage of  $(\bar Z,\bar Y,\bar V)=\phi(\bar W)$ is 
a graph of a function of $W$, as in the required form \eqref{eq:WsQ:0.1}. In order to obtain estimates \eqref{eq:WsQ:0.2} for the derivatives of this function, one
differentiates both sides of \eqref{eq:preWs}, and uses the corresponding estimates in \eqref{tldha}.
\end{proof}

\subsection{Local unstable invariant manifold $W^u_{loc}(Q)$}
We next investigate the unstable manifold of the fixed point $Q$ of $T_k$. Namely, in order to locate the intersection  $W^u_{loc}(Q)\cap W^s(O)$, we need a bound from below on the size of a connected piece of $W^u_{loc}(Q)$. More specifically, we estimate the extent in the $Y$ direction of the intersection of $W^u_{loc}(Q)$ with  a locally invariant center-stable manifold $W^{cs}$ of the map $T_k$ (see Lemma \ref{lem:WuQ}).  

To prove the existence of $W^{cs}$ we use the following 

\begin{thm*}[{\citep[Theorem 4.3]{SSTC1}}]
Let $D_1$ and $D_2$ be some convex closed subsets of some Banach spaces. Suppose a map $T$ (not necessarily single-valued) is defined by the condition that for every $(x,\bar y)\in D_1\times D_2$ one has
\begin{equation*}\label{eq:sstc}
(\bar x,\bar y)=T(x,y)\quad\Leftrightarrow\quad \bar x =F(x,\bar y),\quad y = G(x,\bar y),
\end{equation*}
where $(F,G):D_1\times D_2 \to D_1\times D_2$ is a smooth map.
Let $F$ and $G$ be  such that 
\begin{equation*}\label{eq:sstc0}
\sqrt{\sup_{(x,\bar y)\in D_1\times D_2}\left\{
\left\|\dfrac{\partial F}{\partial x}\right\|
\left\|\dfrac{\partial G}{\partial \bar y}\right\|
\right\}}
+
\sqrt{
\left\|\dfrac{\partial F}{\partial \bar y}\right\|_\circ
\left\|\dfrac{\partial G}{\partial x}\right\|_\circ
}
<1,
\quad
\left\|\dfrac{\partial F}{\partial x}\right\|_\circ
+
\sqrt{
\left\|\dfrac{\partial F}{\partial \bar y}\right\|_\circ
\left\|\dfrac{\partial G}{\partial x}\right\|_\circ
}<1,
\end{equation*}
where $\|\cdot\|$ denotes the matrix norm induced by the norms of the corresponding Banach spaces and 
$$
\left\|\dfrac{\partial F}{\partial \bar y}\right\|_\circ=\sup_{x,\bar y\in D_1\times D_2}\left\|\dfrac{\partial F}{\partial \bar y}\right\|
\quad\mbox{and}\quad
\left\|\dfrac{\partial G}{\partial  x}\right\|_\circ=\sup_{x,\bar y\in D_1\times D_2}\left\|\dfrac{\partial G}{\partial  x}\right\|.
$$
Then, the map $T$ has an attracting invariant  manifold $\mathcal{M}^*$ which is given by the graph $x=w^*(y)$ of some $C^1$ function $w^*$, and is a $C^1$-limit
of iterations by $T$ of the manifold $x=C$ for any constant $C\in D_1$.
\end{thm*}

Here attracting means that points outside $\mathcal{M}^*$ converge to it exponentially by forward iterations. Considering $T^{-1}$ and interchanging $F,x$ and $G,\bar y$, we immediately obtain
\begin{cor} \label{corrr}
Let $(F,G):D_1\times D_2 \to D_1\times D_2$ be a smooth map. Suppose 
\begin{equation*}\label{eq:sstc}
(\bar x,\bar y)=T(x,y)\quad\Leftrightarrow\quad \bar x =F(x,\bar y),\quad y = G(x,\bar y),
\end{equation*}
and
\begin{equation}\label{eq:ssct2}
\sqrt{\sup_{(x,\bar y)\in D_1\times D_2}\left\{
\left\|\dfrac{\partial G}{\partial \bar y}\right\|
\left\|\dfrac{\partial F}{\partial  x}\right\|
\right\}}
+
\sqrt{
\left\|\dfrac{\partial G}{\partial x}\right\|_\circ
\left\|\dfrac{\partial F}{\partial \bar y}\right\|_\circ
}
<1,
\quad
\left\|\dfrac{\partial G}{\partial \bar y}\right\|_\circ
+
\sqrt{
\left\|\dfrac{\partial G}{\partial x}\right\|_\circ
\left\|\dfrac{\partial F}{\partial \bar y}\right\|_\circ
}<1.
\end{equation}
Then, the map $T$ has a repelling invariant manifold of the form $y=w^*(x)$ for some $C^1$ function $w^*$. This manifold is a $C^1$ limit of preimages of $y=C$ for any constant $C\in D_2$. 
\end{cor}

In our setting, $T_k(Z,Y,W,V)=(\bar Z,\bar Y,\bar W,\bar V)$ if the relation \eqref{eq:maps:Tk} is fulfilled. Thus, we can apply the above corollary by putting $x=(Z,Y,W),\; y= V.$ 
For that, we extend \eqref{eq:maps:Tk} to a map $D_1\times D_2 \to D_1\times D_2 $ with $D_1=\mathbb{R}^{d^s+1}$ and $D_2=[-\delta,\delta]^{d^u-1}$. This can be achieved easily by considering some smooth bump function $\xi$ which equals to identity on $[-c\delta,c\delta]^{d^s+1}$ and vanishes outside $[-\delta,\delta]^{d^s+1}$ for some constant $c\in (0,1)$, and define the extension of \eqref{eq:maps:Tk} as 
\begin{equation*}
\begin{aligned}
\bar Z &= \lambda^k \alpha_1 Z+Y+\lambda^k\beta_1W+  h_1(\xi(Z, Y, W), \bar V)=:F_1(Z, Y, W, \bar V),\\
\bar Y &= \hat\mu + \alpha_2\lambda^k\gamma^k Z  + d\gamma^kY^2\varphi_2(Z, Y, W)+  \gamma^k h_2( \xi(Z, Y, W), \bar V)=:F_2(Z, Y, W, \bar V),\\
\bar W &=  \lambda^k \alpha_3 Z+\lambda^k\beta_3W+  h_3( \xi(Z, Y, W), \bar V)=:F_3(Z, Y, W, \bar V),\\
V&= \lambda^k \alpha_4 Z+\lambda^k\beta_4W+  h_4( \xi(Z, Y, W), \bar V)=:G(Z, Y, W, \bar V).
\end{aligned}
\end{equation*}   
Since the estimates in \eqref{eq:maps:Tk_nonlinear} and \eqref{eq:h_df}  imply
\begin{align*}
&\left\|\dfrac{\partial F}{\partial  x}\right\|_\circ=
\left\|\dfrac{\partial (\bar Z,\bar Y,\bar W)}{\partial (Z,Y,W)}\right\|_\circ=O(\lambda^k|\gamma|^k+\delta|\gamma|^k),
\qquad
\left\|\dfrac{\partial G}{\partial  x}\right\|_\circ=
\left\|\dfrac{\partial V}{\partial (Z,Y,W)}\right\|_\circ=O(\lambda^k+\delta),\\
&\left\|\dfrac{\partial F}{\partial  \bar y}\right\|_\circ=
\left\|\dfrac{\partial (\bar Z,\bar Y,\bar W)}{\partial \bar V}\right\|_\circ=O(\delta)+O(\gamma^k\hat\gamma^{-k}),
\qquad
\left\|\dfrac{\partial G}{\partial  \bar y}\right\|_\circ=
\left\|\dfrac{\partial V}{\partial \bar V}\right\|_\circ=O(\gamma^{-k}),
\end{align*}
the fulfillment of condition  \eqref{eq:ssct2} follows immediately for all sufficiently small $\delta$ and sufficiently large $k$. Therefore, Corollary
\ref{corrr} gives us the existence of a manifold  $V=w^{cs}(Z,Y,W)$
such that the part of this manifold inside the region
\begin{equation*}\label{eq:Pi'}
\Pi':=[-c\delta,c\delta]\times [-c\delta,c\delta]\times [-c\delta,c\delta]^{d^s-1}\times [-c\delta,c\delta]^{d^u-1}\subset \Pi
\end{equation*}
is a locally invariant manifold $W^{cs}$ of $T_k$. Since it is the limit
of the preimages of the plane $\{V=V_Q\}$ which is tangent to 
the cone field  $\mathcal{C}^{cs}$ given by \eqref{eq:cone_cs},
the backward invariance of $\mathcal{C}^{cs}$ implies that the tangents to $W^{cs}$ lie in  $\mathcal{C}^{cs}$ (up to increasing the coefficient $K$ in \eqref{eq:cone_cs} slightly if necessary). 

By decreasing $\delta$, one can take  $\Pi'$ as the new $\Pi$ and obtain the following
\begin{lem}\label{lem:csmanifold}
For each sufficiently large $k$, there exists
a $C^1$-smooth repelling manifold $W^{cs}$ in $\Pi$,
which is locally invariant in the sense that for any point $M\in W^{cs}$, its iterates by $T_k$ lie in $W^{cs}$ as long as they belong to $\Pi$. Moreover, the invariant manifold is given by 
$\{V=w^{cs}(Z,Y,W)\}$ for some function $w^{cs}$ defined on $[-\delta,\delta]\times [-\delta,\delta]\times [-\delta,\delta]^{d^s-1}$ and such that
$w^{cs}(Z_Q,Y_Q,W_Q)=V_Q$. The tangents
to the manifold lie in the cone field $\mathcal{C}^{cs}$ given by \eqref{eq:cone_cs}, so $\|w^{cs}\|_{C^1}=O(\delta)$.
\end{lem} 

We are now in the position to describe $W^u_{loc}(Q)$.

\begin{lem}\label{lem:WuQ}
For all sufficiently large $k$, there exists a connected piece $W^u_{loc}(Q)$ of $W^u(Q)\cap \Pi$, which has the form
\begin{equation}\label{eq:WuQ:0}
W=W_Q+\phi^u(Z,Y,V),
\end{equation}
where $\phi^u$ is a $C^{r-2}$ function (with respect to variables and parameters) satisfying
\begin{equation}\label{eq:WuQ:00}
\phi^u(Z_Q,Y_Q,V_Q)=0 \quad\mbox{and}\quad
\|\phi^u\|_{C^1}=O(\delta).
\end{equation} 
The function $\phi^u$ is defined on the region $\{(Z,Y)\in D(V), \|V\|\leq \delta\}$ where $D(V)$ is, at every $V$, a connected open subset of $[-\delta,\delta]\times [-\delta,\delta]$ such that either all regions $D(V)$ intersect both $\{Y=K_1\delta^3\}$ and $\{Y=K_2\delta^2\}$, or all intersect both $\{Y=-K_1\delta^3\}$ and $\{Y=-K_2\delta^2\}$, where $K_{1,2}$ are some positive constants independent of $\delta$ and $k$.
\end{lem} 

\begin{proof}
The invariance of the cone field  $\mathcal{C}^{cu}$ given by \eqref{eq:cone_cu} implies that there is a piece $w^u$ of $W^u(Q)$ which contains $Q$ and has the form \eqref{eq:WuQ:0} with $\phi^u(Z,Y,V)$ defined on a small neighborhood of $(Z_Q,Y_Q,V_Q)$ and satisfying \eqref{eq:WuQ:00}. By the invariance of $\mathcal{C}^{cu}$, further iterations of $w^u$ have the same form, so
 it suffices to show that some forward iterate of $w^u$ contains points whose $(Z,Y,V)$-coordinates cover the desired region 
$\{(Z,Y)\in D(V)\}$.

Let us first investigate the iterations of $w^u\cap W^{cs}=:w^c$, where $W^{cs}$ is the center-stable manifold given by Lemma \ref{lem:csmanifold}. Since $w^u$ and $W^{cs}$ are tangent to $\mathcal{C}^{cu}$ and, respectively, $\mathcal{C}^{cs}$, the piece $w^c$ has the form
\begin{equation}\label{eq:WuQ:1}
(W,V)=w^c(Z,Y)=:(w^c_1(Z,Y),w^c_2(Z,Y)),
\end{equation}
where $(W_Q,V_Q)=w^c(Z_Q,Y_Q)$. The partial derivatives of the functions $w^c_{1,2}$ have order $O(\delta)$; this estimate persist under the  forward iteration by $T_k$ (by the forward invariance of the cone field
$\mathcal{C}^{cu}$ and the manifold $W^{cu}$).

We claim that there exists $n$ such that  $T^i_k(w^c)\subset \{|Z|< \delta^{3/2}, |Y|< \delta^2\}$ for $i=0,\dots,n-1$,
and
$T^n_k(w^c)$  intersects the boundary 
$\{|Y|=\,\delta^2\}$. In what follows we prove the claim.

By \eqref{eq:maps:Tk}, for any surface of the form $(W,V)=w(Z,Y)=(w_1(Z,Y),w_2(Z,Y))$, the restriction
of $T_k$ to this surface is given by $(Z,Y)\mapsto (\bar Z, \bar Y)$ where
\begin{align*}
\bar Z &=  \lambda^k \alpha_1 Z+Y+\lambda^k\beta_1 w_1(Z,Y)+ h_1(Z,Y,w_1(Z,Y),w_2(\bar Z,\bar Y)),\\
\bar Y &= \hat\mu+ \alpha_2 \lambda^k\gamma^k Z  +d\gamma^kY^2+ \gamma^k h_2(Z,Y,w_1(Z,Y), w_2(\bar Z,\bar Y)).
\end{align*}
Since the derivatives of $w^c$ and its iterates are $O(\delta)$,
we have from \eqref{eq:maps:Tk_nonlinear} and \eqref{eq:h_df} that
$$\det \frac{\partial (\bar Z,\bar Y)}{\partial (Z,Y)}= - \alpha_2 (\lambda\gamma)^k (1+ O(\delta))$$
for the restriction of $T_k$ on $T^i_k(w^c)$ for any $i$ (as long as $T^i_k(w^c)$ lies in $\Pi$).
Since $|\lambda\gamma|>1$ and $\alpha_2$ is bounded away from zero, it follows that
$|\det \frac{\partial(\bar Z,\bar Y)}{\partial (Z,Y)}| >1$ for large enough $k$ and small $\delta$. Therefore,
the area of $T^i_k(w^c)$ grows exponentially with $i$ until $T^i_k(w^c)$ intersects the boundary of
$\Pi$.

It follows that there exists $n$ such that $T^n_k(w^c)$ intersects for the first time the union of $\{Z=\pm\delta^{3/2}\}$ and 
$\{Y=\pm\delta^2\}$. In fact, $T^n_k(w^c)$ must intersect $\{Y=\pm\delta^2\}$ and lie inside $\{|Z| < \delta^{3/2}\}$. To see this, note from formula \eqref{eq:maps:Tk} that, as long as the $Y$-coordinates of points of $T_k^i(w^c)$ are bounded by 
$\delta^2$, the $Z$-coordinates of points of $T_k^{i+1}(w^c)$ are of order $O(\delta^2)\ll \delta^{3/2}$. The claim is proven.

We can now finish the proof of the lemma. Take the surface $T^n_k(w^c)$ we just constructed. It is a graph of a function 
$(Z,Y)\mapsto (W,V)$ whose domain $\tilde D$ is connected and intersects $\{|Y| =\delta^2\}$. Note that all backward iterations of
this surface lie in the region $\Pi_Y$ defined by \eqref{eq:PiY}. Note also that the existence of the forward-invariant cone field $\mathcal{C}^{uu}$ in Lemma \ref{lem:cones}  guarantees the existence of an invariant strong-unstable  $C^0$ lamination $\mathcal{L}^{uu}$ in the region $\Pi_Y$, which consists of $C^{r}$ leaves through points whose entire backward orbits stay in $\Pi_Y$. In particular, there is a leaf through each point of $T^n_k(w^c)$, and the union $w^n$ of these leaves is a connected piece of 
$W^u(Q)$.

The leaves form a $C^0$ foliation of $w^n$. Each leaf $\ell^{uu}$ is the graph of some function $(Z,Y,W)=\ell^{uu}(V)$, and it is defined at $V\in [-\delta,\delta]^{d^u-1}$ as long as $(\ell^{uu}(V),V)\in \Pi_Y$. Since the tangents to $\ell^{uu}$ lie in $\mathcal{C}^{uu}$, it follows from 
\eqref {eq:coneuu_df} that the change in $(Z,Y,W)$ as $V$ varies in a range of order $\delta$ along any such leaf is $O(\delta^3)$,
provided that $k$ is taken large enough. Therefore, the manifold $w^n$ is a graph of a function $(Z,Y,V)\mapsto W$ whose domain 
contains the set $\{(Z,Y)\in D(V), \|V\| \leqslant \delta\}$ where $D$ is a homeomorphic image of $\tilde D$
(the domain of the function defining the surface $T^n_k(w^c)$) by the holonomy map of the foliation, which is $O(\delta^3)$
close to the identity. Since $\tilde D$ is connected and contains a point with the coordinate $Y$ such that $|Y|=\delta^2$ and
the point $(Z_Q,Y_Q)$ with the $Y$-coordinate $O(\lambda^k)=O(\delta^3)$ (if $k$ is large enough), the statement of the lemma follows.\end{proof} 

\subsection{Local extended stable manifold $W^{sE}_{loc}(Q)$}\label{sec:Wse}
A necessary step in the proof of Theorem \ref{thm:sf} is stabilizing a heterodimensional cycle involving $O$ and $Q$.
For that, we need to verify non-degeneracy conditions of \citep{LT21} for this heterodimensional cycle.  Particularly, the extended stable manifold of $Q$ is considered in condition NC1, see Section \ref{sec:HDCcon}.

According to Remark \ref{saddle1}, when $t$ is bounded away from zero, the nearest to the unit circle unstable multiplier
$\nu$ of the point $Q$ is real. The corresponding extended stable subspace of $DT_k|_{Q}$ is the direct sum of the eigenspace
corresponding to the multipliers that are less than $1$ in the absolute value (this subspace lies in the cone $\mathcal{C}^{ss}$ of \eqref{eq:cone_ss}) and the eigenvector corresponding to the multiplier $\nu$. The relation  \eqref{eigen:4} along with the fact that any eigenvector of $\nu$ also lies in $\mathcal{C}^{cs}$ of \eqref{eq:cone_cs}   implies that, for $t\neq 0$,  $|z - t y|\leq K\delta$ and $\|v\| \leq 2K\delta ((1+t)|y|+\|w\|)$ for all sufficiently large $k$. Obviously, these two inequalities are also satisfied by vectors in $\mathcal{C}^{ss}$. Thus,   for $t$ bounded away from zero, 
the vectors $(z,y,w,v)$ in the extended stable subspace satisfy
\begin{equation}\label{conewse}
|z - t y|+\|v\| \leq K'\delta (\|w\|+|y|).
\end{equation}
for some constant $K'$, independent of $k$. The local extended stable manifold $W^{sE}_{loc}(Q)$ is an invariant manifold which is tangent at $Q$ to the extended stable subspace, see e.g.  \cite{SSTC1}. The global extended stable manifold $W^{sE}(Q)$ is the union of all backward orbits that start in 
$W^{sE}_{loc}(Q)$. The manifold $W^{sE}(Q)$ contains the stable manifold  $W^{s}(Q)$. 
\begin{lem}\label{lem:WseQ}
Let $t_-,t_+$ be as in \eqref{tpm}. For $t\in(t_-,t_+)$    bounded away from zero, the tangents to $W^{sE}_{loc}(Q)$
at the points  $(Z,Y,W,V)\in W^{s}_{loc}(Q)$ such that
\begin{equation}\label{yboq}
Y - Y_Q= \delta O(\lambda^k),
\end{equation}
lie in the cone defined by \eqref{conewse}.
\end{lem}
\begin{proof}
Note that at $Q$ (hence for the points close $Q$) the tangent to
$W^{sE}_{loc}(Q)$ does belong to the cone \eqref{conewse}. Therefore, since $W^{sE}(Q)$ is the union of backward orbits starting at
$W^{sE}_{loc}(Q)$, it is enough to check that the cone \eqref{conewse}
is invariant under the action of the derivative $DT_k^{-1}$.

Take a point $P\in W^{s}_{loc}(Q)$ whose $Y$-coordinate 
satisfy \eqref{yboq}. Suppose that $(\bar z,\bar y,\bar w,\bar v)=DT_k|_{P}  (z, y, w, v)$ and 
that $(\bar z,\bar y,\bar w,\bar v)$ lies in the  cone \eqref{conewse}, that is, 
\begin{equation}\label{eq:saddle110}
|\bar z - t \bar y|+\|\bar v\| \leq K'\delta (\|\bar w\|+|\bar y|).
\end{equation}
We need to show that $(z, y, w, v)$ also lies in the cone.

Note that $Y^{new}=\delta O(\lambda^k)$ by \eqref{newcoor} and \eqref{yboq}. Hence, by resolving $(Z,Y,\bar W, V)^{new}$ as functions of $(\bar Z,\bar Y, W,\bar V)^{new}$ through  \eqref{eq:preWs} and differentiating these functions,  one can find  that
\begin{equation}\label{referred:2}
\|(y -  \bar z, z - ty, \bar w, v)\| =  O(\delta) \|(\bar z, \bar y, w, \bar v)\|.
\end{equation}
Indeed, one first substitute the second equation in \eqref{eq:preWs} into the first one for $Y^{new}$; then, since the right-hand side of the new  $Y^{new}$-equation has small first derivative with resepect to $Y^{new}$ (due to our condition $Y^{new}=\delta O(\lambda^k)$ that makes the derivative coming from $\tilde h_1$ by the chain rule small), one finds $Y^{new}$ as a function of $(\bar Z,Z, W,\bar V)^{new}$; next, combining this expression for $Y^{new}$  with the remaining equations in \eqref{eq:preWs}, we get the new formula with  variables $(\bar Z,\bar Y, W,\bar V)^{new}$ as arguments; finally, differentiating this formula after sorting terms yields \eqref{referred:2}

The relation \eqref{referred:2} implies $\bar z\sim y$ and $|z - ty|+\|\bar w\|+\|v\| = O(\delta)(|\bar z|+|\bar y|+\|w\|+\|\bar v\|)$.
Taking into account that  $\bar y\sim \bar z/t$ and $\bar v =O(\delta)(\|\bar w\|+|\bar y|)$, according to \eqref{eq:saddle110}, one readily sees that  $(z, y, w, v)$ indeed lies in the cone \eqref{conewse} for an appropriate $K'$.
\end{proof}

\subsection{Local strong-unstable manifold $W^{uu}_{loc}(Q)$}\label{sec:WuuQ}
We also collect some facts on $W^{uu}_{loc}(Q)$, which are used for verifying the non-degeneracy condition NC2 in Section \ref{sec:HDCcon}.
The strong-unstable subspace of $DT_k|_{Q}$ is the eigenspace corresponding to all  eigenvalues of $DT_k|_{Q}$ outside the unit circle, except for the one with the smallest absolute value. As  mentioned in Remark~\ref{saddle1}, for $t$ bounded away from zero, the two nearest to the unit circle unstable multipliers $\nu$ and $\hat\nu$ of the point $Q$ are real, and $|\nu| < |\hat\nu|$ (see \eqref{eigen:3}). Hence, the strong-unstable subspace of $DT_k|_{Q}$ is
the direct sum of the eigenspace lying in $\mathcal{C}^{uu}$ (see \eqref{eq:coneuu_df}) and the eigenvector corresponding to the multiplier $\hat\nu$. Since this eigenvector lies in $\mathcal{C}^{cu}$ given by \eqref{eq:cone_cu} and
its $(z,y)$-component is given by \eqref{eigen:4},
the strong-unstable subspace of $DT_k|_{Q}$ lies in the cone
\begin{equation}\label{zywveq}
\|(z,w)\|\leq K (\lambda\gamma)^{-k}\|(y,v)\|
\end{equation}
for some constant $K$. The strong-unstable manifold $W^{uu}(Q)$ is a submanifold of $W^u(Q)$, tangent to the strong-unstable subspace at $Q$. We call the local strong-unstable manifold the connected piece of $W^{uu}(Q)$ that contains $Q$
and lies in the domain of the map $T_k$. It is obtained by iterating a sufficiently small connected piece of 
$W^{uu}(Q)$ around $Q$. In what follows we find a formula for it.

First note that the cone field defined by \eqref{zywveq} in a small neighborhood of $Q$ is invariant under $DT_k$.  By the formula \eqref{locmapq} for the map $T_k$, the expression of $Y_Q$ in \eqref{eq:Q_coor}, and the estimates \eqref{hagain}, the derivative
$DT_k$ takes a vector $(z,y,w,v)$ to a vector $(\bar z,\bar y, \bar w,\bar v)$ such that
\begin{equation}\label{conepointq}
\begin{aligned}
|\bar z| &=  |y| + O(\lambda^k)\|(z,y,w,\bar v)\|, \\
|\bar y| &=  |b_1\alpha^*\lambda^k\gamma^k (z - t y (1+ Y^{new}/Y_Q+O(W))|+ O(\hat\lambda^k\gamma^k)\|(w,\bar v)\|,\\
\|(\bar w, v)\| &= O(\lambda^k)\|(z,y,w,\bar v)\|.
\end{aligned}
\end{equation}

This implies that there exists a constant $K>0$ such that the image of the cone \eqref{zywveq}
by $DT_k$ is again in this cone, provided that the factor $1+Y^{new}/Y_Q+O(W)= Y/Y_Q + O(W)$ is bounded away from zero.
Moreover, $DT_k$ is expanding in the cone  \eqref{zywveq} because
\begin{equation}\label{expan+u}
\|v\| = O(\lambda^k) \|\bar v,y\|, \qquad |y|  = O((\lambda\gamma)^{-k})|\bar y|+ O((\hat\lambda/\lambda)^k) \|\bar v\|.
\end{equation}

Therefore, for any $C>0$ and all large $k$, the image of any manifold, which lies in the region $Y/Y_Q>C$, contains the point $Q$
and has its tangents lying in the cone \eqref{zywveq}, is again a manifold with tangents in this cone. Since a small piece of 
$W^{uu}_{loc}(Q)$ near $Q$ is a manifold of this type, the iterations of this piece satisfy \eqref{expan+u} as long as their preimages stay at $Y/Y_Q>C$. Thus, $W^{uu}_{loc}(Q)$ is given by a function of the form
\begin{equation} \label{eqwuuq}
(Z,W)=\phi^{uu}(Y,V)
\quad\mbox{with}\quad
 \frac{\partial \phi^{uu}}{\partial (Y,V)} = O(|\lambda\gamma|^{-k}).
\end{equation}

Let us now estimate the size of $W^{uu}_{loc}(Q)$. Note that it contains an invariant submanifold $W^{uuu}_{loc}(Q)$ which is
tangent at $Q$ to the eigenspace of $DT_k|_{Q}$ lying  in the cone $\mathcal{C}^{uu}$ of \eqref{eq:coneuu_df}. This manifold 
can be defined by the equation
\begin{equation}\label{uuuu}
Y=Y_Q + \phi^{uuu}(V-V_Q),
\end{equation}
where
\begin{equation}\label{uuuest}
\frac{\partial \phi^{uuu}}{\partial V} = O(\delta^2).
\end{equation}
It divides $W^{uu}_{loc}(Q)$ into two parts. 
When the multiplier $\hat\nu$ is positive -- this is the case when 
\begin{equation}\label{tsign}
t b_1 \alpha^* (\lambda\gamma)^k < 0,
\end{equation}
see \eqref{eigen:3}, these parts are invariant, each, by $T_k$. 
Denote as $W^{uu+}_{loc}(Q)$ the invariant part corresponding to
$$\frac{Y}{Y_Q} > 1 + \frac{1}{Y_Q} \phi^{uuu}(V-V_Q).$$

\begin{lem}\label{lem:WuuQ} For $t\in (t_-,t_+)$ bounded away from zero and such that \eqref{tsign} holds,
for all sufficiently large $k$, the piece $W^{uu+}_{loc}(Q)$ defined by \eqref{eqwuuq}  extends up to  $Y = \delta \cdot {\rm sign}(Y_Q)$.
\end{lem}
\begin{proof} 
Take any $C>0$ such that \eqref{eqwuuq} is valid for $Y/Y_Q>C$.
By \eqref{uuuest}, the part of $W^{uuu}_{loc}(Q)$ corresponding to $\|V-V_Q\|\leq {C_1}{\delta^{-2}}\lambda^{k}$ for some constant $C_1>0$ lies in the region $Y/Y_Q>C$. Therefore, when iterating the small piece of $W^{uu+}_{loc}(Q)$ in order to obtain the entire $W^{uu+}_{loc}(Q)$, the points in $W^{uu+}_{loc}(Q)$ will remain at $Y/Y_Q>C$ as long as $\|V-V_Q\|\leq {C_1}{\delta^{-2}}\lambda^{k}$. So, the tangents to $W^{uu+}_{loc}(Q)$ lie in the cone \eqref{zywveq}, as required, and the expansion
given by \eqref{expan+u} implies that the domain of the function $\phi^{uu}$ that defines $W^{uu+}_{loc}(Q)$ in \eqref{eqwuuq}  
indeed extends to the entire region
$$\frac{\delta}{|Y_Q|} \geq  \frac{Y}{Y_Q} > 1 + \frac{1}{Y_Q} \phi^{uuu}(V-V_Q), \qquad \|V\|\leq \delta.$$
\end{proof}

\section{Heterodimensional cycle involving $O$ and $Q$}\label{sec:createhdc}

In this section, we   prove 
\begin{thms}\label{thm:sfsimple}
Suppose $|\lambda\gamma|>1$ and conditions C1--C5 are satisfied. Assume that $L$ is a saddle-focus of type (2,1) and has index $d^u$. Then, for any proper unfolding $\{f_{\mu,\omega}\}$ with   $f_{0,\omega^*}=f$, there exists a sequence $(\mu_k,\omega_k)\to (0,\omega^*)$ such that $f_{\mu_k,\omega_k}$ has a heterodimensional cycle involving the continuation of $L$ and a periodic orbit $\hat L$ of index $d^u+1$. Moreover, this cycle unfolds as $\mu$ varies from $\mu_k$, or $\omega$ varies from  $\omega_k$.
\end{thms}
This is done in Sections \ref{sec:nontr} and \ref{sec:trans} by creating a heterodimensional cycle involving $O$ and the 
point $Q$ described in the previous section, which amounts to showing that the intersections $W^u(O)\cap W^s(Q)$ and $W^u(Q)\cap W^s(O)$ are non-empty.   

After that, we show in Section \ref{sec:robusthdc} that this heterodimensional cycle satisfies the non-degeneracy conditions in \citep{LT21}. This enables us to  apply \citep[Theorem 7]{LT21} to obtain robust heterodimensional dynamics involving blenders, thus proving  Theorem \ref{thm:sf} for the saddle-focus case.

Recall that given an orbit of a homoclinic tangency $\Gamma$, we choose two points, 
$M^-\in W^u_{loc}(O)$ and $M^+\in W^s_{loc}(O)$ such that the global map $T_1$ acts from a small neighborhood of $M^-$ to a small neighborhood of $M^+$, and
$M^{\pm}$ lie in a sufficiently small neighborhood of $O$ where the local map is brought to the form \eqref{eq:maps:T0} and  \eqref{eq:maps:T0_nonlinear}. 

\begin{defi}[Accompanied homoclinic tangencies]
We call the orbit of a homoclinic tangency {\em accompanied} if in any neighborhood of
$M^-$ there exist points of a transverse intersection of $W^s(O)$ with $W^u_{loc}(O)$, and among these points there are
those with the coordinate $y>y^-$ and those with $y<y^-$, where $y^-$ is the $y$-coordinate of $M^-$.
\end{defi}

The two lemmas below (proven in Section \ref{sec:twolem}) imply that arbitrarily close to $(\mu,\omega)=(0,\omega^*)$
there exist parameter values for which $f_{\mu,\omega}$ has an accompanied homoclinic tangency satisfying conditions C1--C5
and the family $f_{\mu,\omega}$ is a proper unfolding of this new tangency. This means that it is enough to prove the theorem
for accompanied tangencies only.  Recall that $d$ is the coefficient in \eqref{eq:maps:T1_cross}.
\begin{lem}\label{lem:tangency1}
Let $\omega/(2\pi)$ be either irrational or equal to $p/q$ with $q\geqslant 6$. If conditions C1--C4 are satisfied and 
$$d y^- >0,$$
then the tangency is accompanied. 
\end{lem}

\begin{lem}\label{lem:tangency2}
Let $\omega/(2\pi)$ be either irrational or equal to $p/q$ with $q\geqslant 6$. For any smooth one-parameter family $f_{\mu}$ with $f_0$ satisfying conditions C1--C5, there exists parameter values accumulating on $\mu=0$,
for which $f_{\mu}$ has a homoclinic tangency with $d y^->0$. Moreover, the newly created tangency satisfies conditions  C1--C5 and splits as $\mu$ varies.
\end{lem}


\subsection{Non-transverse intersection of $W^u(O)$ with $ W^s(Q)$}\label{sec:nontr} 
 In this subsection, we prove
\begin{prop}\label{prop:nontr}
Let $f_{\mu,\omega}$ be a proper family, and let the homoclinic tangency at $(\mu=0,\omega=\omega^*)$
be accompanied. Let the periodic point $Q$ and the function $\mu_k(\omega,t)$ be given by Lemma \ref{lem:index2}. There exist an infinite subset $\mathcal{K}\subset \mathbb N$ associated with functions $\omega_k(t)$ for $k\in\mathcal{K}$,
which are  defined for $t\in (t_-,t_+)$ given by \eqref{tpm} and converge to $\omega^*$ as $k\to\infty$, such that the intersection $W^u(O)\cap W^s(Q)$ is non-empty when $\mu=\mu_k(\omega_k(t),t)$ and $\omega=\omega_k(t)$.
Moreover, for any fixed $t$, this intersection unfolds as either $\mu$ or $\omega$ changes.
The functions $\omega_k$ are monotone with respect to $t$:
\begin{equation}\label{omtd}
\frac{d\omega_k(t)}{dt}= \Omega_k \lambda^k + o(\lambda^k),
\end{equation}
for a constant $\Omega_k\neq 0$.
\end{prop}

Recall that the intersection $W^u(O)\cap W^s(Q)$ found at $(\mu_k,\omega_k)$  is non-transverse, as $\dim W^u(O)+\dim W^s(Q) < \dim \mathcal{M}$. So, by saying that the intersection unfolds we simply mean that it disappears with non-zero velocities as the parameter changes. 

\begin{proof}
We divide the proof into several steps.

\noindent\textbf{(1) Plan of proof.}
Since the tangency is accompanied by transverse homoclinics, there exist pieces
of $W^u(O)$ which intersect $W^s_{loc}(O)$ transversely and sufficiently close to $M^+$. 
Choose such piece $W^{u*}$; the transversality means that it is a smooth manifold of the form  
$$
(x_1,x_2,u)=(x_1^*,x_2^*,u^*) + O(y,v),
$$
for some constant $(x_1^*,x_2^*,u^*)$. The closeness to $M^+$ means $(x_1^*,x_2^*,u^*)$ is close to $(x_1^+,x_2^+,u^+)$, which by \eqref{eq:x+2} implies 
\begin{equation}\label{eq:xstar}
(x_1^*,x_2^*)\neq 0.
\end{equation}
By \eqref{eq:maps:T0k}, the image $S_j:=T_0^j(W^{u*})$ is given by 
\begin{equation}
\begin{aligned}
\tilde x_1 =\lambda^j (x_1^* \cos j\omega- x_2^* \sin j\omega)+O(\hat\lambda^j),\quad
\tilde x_2&=\lambda^j (x_1^* \sin j\omega+ x_2^* \cos j\omega)+O(\hat\lambda^j),\quad
\tilde u =O(\hat\lambda^j),
\end{aligned}
\end{equation}
where the terms $O(\hat\lambda^j)$ are functions of $(\tilde y,\tilde v)$ with $\tilde y-y^- \in[-\delta,\delta]$ and $\tilde v-v^- \in[-\delta,\delta]^{d^u-1}$, and of the parameters $(\mu,\omega)$, and the first derivatives of these functions with respect to the variables and parameters are also $O(\hat\lambda^j)$. We  prove the proposition by finding an intersection point $P$ of $T_1(S_j)$ with $ W^s_{loc}(Q)$.\\

\noindent\textbf{(2) Formula for $T_1(S_j)$.} 
 By \eqref{eq:maps:T1_cross}, the image $T_1(S_j)$ satisfies
\begin{equation}\label{eqsj8}
\begin{aligned}
 x_1-x_1^+&=\lambda^j A_1
 + b_1(\tilde y-y^-)+O(v)+O((\tilde y-y^-)^2)+ O(\hat\lambda^j),\\
 x_2-x_2^+&=\lambda^j A_2
 +O(v)+O((y-y^-)^2)+O(\hat\lambda^j),\\
 y&=\mu+\lambda^j A_3
+  O(v)+d (\tilde y-y^-)^2 +O((\tilde y-y^-)^3)+O(\hat\lambda^j),\\
 u-u^+&=b_4(\tilde y-y^-)+O(v)+O((\tilde y-y^-)^2)+O(\lambda^j),\\
 \tilde v-v^-&=b_5(\tilde y-y^-)+O(v)+O((\tilde y-y^-)^2)+O(\lambda^j),
\end{aligned}
\end{equation}
where we denote
\begin{equation}\label{eq:AB}
\begin{aligned}
&A_i=a_{i1}(x_1^* \cos j\omega- x_2^* \sin j\omega) + a_{i2}(x_1^* \sin j\omega+ x_2^* \cos j\omega), \quad i=1,2,\\
&A_3=c_1(x_1^* \cos j\omega- x_2^* \sin j\omega) + c_2(x_1^* \sin j\omega+ x_2^* \cos j\omega);
\end{aligned}
\end{equation}
the $O(\hat\lambda^j)$ and $O(\lambda^j)$ terms are smooth functions of $(\tilde y,\tilde v)$ and $v$, and of $(\mu,\omega)$;
the first partial derivatives of the $O(\hat\lambda^j)$ are also $O(\hat\lambda^j)$ and the first derivatives of 
$O(\lambda^j)$ with respect to $(y,\tilde v,v)$ are also $O(\lambda^j)$, whereas the first derivatives with respect to $(\mu,\omega)$
are $O(j\lambda^j)$ (since $\lambda$ is a smooth function of the parameters, $\partial_{(\mu,\omega)} \lambda^j= O(j\lambda^j)$).

Lemma \ref{lem:WsQ} provides a description for $W^s_{loc}(Q)$ in coordinates $(Z,Y,W,V)$ 
on $\Pi^+\cap T_0^{-k}\Pi^-$, introduced in Section
\ref{norforfir}. To find an intersection with $T_1(S_j)$, we rewrite the above equation \eqref{eqsj8} in the same
coordinates $(Z,Y,W,V)$. 

We take $k$ large enough such that 
\begin{equation}\label{eq:j}
\hat\lambda^j \gg \lambda^{-k}\hat\lambda^k.
\end{equation}
Apply the coordinate transformation \eqref{eq:maps:coor1}
to \eqref{eqsj8}. Note that in doing so, we express $v$ as a function of $(x,\tilde y,u,\tilde v,\eps)$, as it is given 
by the last equation in \eqref{eq:maps:T0k}. For the conformity with the bi-focus case, we only use these estimates
for $v$:
\begin{equation}\label{eq:weakerv}
v=O(\gamma^{-k}),\quad 
\dfrac{\partial v}{\partial (x,\tilde y,u,\tilde v,\mu)}=O(\gamma^{-k}),\quad
\dfrac{\partial v}{\partial \omega}=O(k\gamma^{-k}),
\end{equation}
which are weaker than those given by \eqref{eq:maps:T0k}. The resulting equations
for $T_1(S_j)$ in the Shilnikov coordinates on $\Pi^+\cap T_0^{-k}\Pi^-$ are as follows
(we use \eqref{eq:j} and \eqref{eq:lambdagamma}): 
\begin{equation*}\label{eq:W^uO1}
\begin{aligned}
 X_1&=\lambda^j A_1+b_1 s+O(s^2)+O(\hat\lambda^j),\\
 X_2&=\lambda^j A_2+O(\gamma^{-k})+O(s^2)+O(\hat\lambda^j),\\
\gamma^{-k} Y&=\mu +\lambda^j A_3+ d s^2+ O(s^3) + O(\hat\lambda^j),\\
 U&=b_4 s+O(s^2)+O(\lambda^j),
\end{aligned}
\end{equation*}
where we denote $s=\tilde y-y^-$; the $O(\cdot)$ terms are functions of $V$ and $s$, and of the parameters. The convention for the
derivatives of $O(\hat\lambda^j)$ and $O(\lambda^j)$ are the same as before (in equation\eqref{eqsj8}). The derivatives
of $O(s^2)$ with respect to $(V,\mu,\omega)$ are $O(s^2)$ as well, and the derivative with respect to $s$ is $O(s)$.
The derivatives of the $O(s^3)$ term with respect to all variables and parameters are estimated as $O(s^2)$ (because the minimal smoothness we operate is $C^3$).

Applying the remaining coordinate transformations \eqref{eq:maps:coor2}, \eqref{eq:maps:coor2.1} and \eqref{eq:maps:coor3} of Section \ref{norforfir}, we obtain the following formula for $T_1(S_j)$:
\begin{align}\label{nontr1}
\begin{split}
Z&=\lambda^j\frac{\alpha^*A_1+\beta^* A_2}{b_1\alpha^*} + s +O(s^2)+O(\hat\lambda^j),\\
\gamma^{-k} Y&=\mu+\lambda^j A_3+d s^2 + O(s^3)+O(\hat \lambda^j),\\
W&=O(s^2)+O(\lambda^j),
\end{split}
\end{align}
where the same convention for the derivatives of $O(\cdot)$ terms as in \eqref{eq:W^uO1} holds, and the 
$O(\cdot)$ terms depend on $V$ and $s$.  Recall that
 $Z$ and $Y$ are center coordinates, $W$ and $V$ stand for the strong-stable, and, respectively, strong-unstable coordinates.
\\

\noindent\textbf{(3) Finding the non-transverse intersection of $T_1(S_j)$ with $W^s_{loc}(Q)$.} 
By Lemma \ref{lem:WsQ}, the local stable manifold $W^s_{loc}(Q)$ has the form
\begin{equation}\label{nontr2}
Z=Z_Q+\phi^s_z(W;\omega,t),\qquad
Y=Y_Q+\phi^s_y(W;\omega,t),\qquad
V=V_Q+\phi^s_v(W;\omega,t),
\end{equation}
where $\phi^s$ satisfies \eqref{eq:WsQ:0.2} and the coordinates of $Q$ are given by Lemma \ref{lem:index2}.

The sought non-transverse intersection of $W^u(O)$ and $W^s(Q)$ corresponds to a solution $(s,W,\omega,t)$ to the system consisting of equations \eqref{nontr2} and \eqref{nontr1} at $\mu=\mu_k(\omega,t)$ given by \eqref{eq:mu_Q}. 
In what follows we solve this system. Note that all the coefficients in \eqref{nontr1} are functions of the parameters
$\mu$ and $\omega$, but we can replace them by their values at $\mu=0$ (the difference is absorbed by
the $O(\hat\lambda^j)$ terms, as $\mu=\mu_k=O(\lambda^{2k})$). So we assume these coefficients to be functions of $\omega$ only. Since the $O(s^2)$, $O(\hat\lambda^j)$ and $O(\lambda^j)$ terms depend on $\mu$, and $\mu=\mu_k$ depends on $t$, these terms should be considered as $t$-dependent, and the derivative with respect to $t$ is estimated (by the chain rule) as $\lambda^{2k}O(s^2)$, $O(\lambda^{2k}\hat\lambda^j)$ and
$O(\lambda^{2k}\lambda^j)$, respectively, since $\partial_t\mu_k = O(\lambda^{2k})$ by \eqref{eq:t_deri}.

Substituting the $V$-equation in \eqref{nontr2} into the $W$-equation in \eqref{nontr1}, one expresses $W$ as a function of $s,\omega,t$:
\begin{equation}\label{nontr:7a}
W=O(s^2) + O(\lambda^j).
\end{equation}
(This expression is different from the last one in \eqref{nontr1} as the latter also depends on $V$.)
Then one expresses $V$ from the $V$-equation in \eqref{nontr2} as a function of $(s,\omega,t)$ and substitutes the result into the right-hand sides of the $Z$ and $Y$ equations in \eqref{nontr1}, so they become functions of $(s,\omega,t)$ only. Replacing $\mu$ by $\mu_k(t,\omega)$ from \eqref{eq:mu_Q} gives then
\begin{equation}\label{zy7a}
\begin{array}{l}\displaystyle
Z=\lambda^j\frac{\alpha^*A_1+\beta^* A_2}{b_1\alpha^*} + s +O(s^2)+O(\hat\lambda^j),\\ \displaystyle
\gamma^{-k} Y = 
\dfrac{(b_1\alpha^*)^2}{2d}\lambda^{2k} (t-\dfrac{t^2}{2})+\lambda^j A_3+d s^2 +O(s^3) +O(\hat \lambda^j).
\end{array}
\end{equation}
We remind that the $O(\hat\lambda^j)$ terms here have the first derivative with respect to $(s,\omega)$ bounded
by $O(\hat\lambda^j)$ and the derivative with respect to $t$ bounded by $\lambda^{2k}O(\hat\lambda^j)$.
Thus, all terms in the formula \eqref{eq:mu_Q} for $\mu_k$, except for the one we keep as the first term
in the right-hand side of the $Y$ equation, are absorbed by $O(\hat\lambda^j)$ due to \eqref{eq:j}.
The derivatives of the $O(s^2)$ and $O(s^3)$ terms with respect to $t$ are $\lambda^{2k}O(s^2)$.

Since the derivatives \eqref{eq:WsQ:0.2} of  \eqref{nontr2} are small, one can substitute 
\eqref{nontr:7a} and \eqref{eq:Q_coor}  into \eqref{nontr2} and obtain by the Implicit Function Theorem that
\begin{equation}\label{zy8}
Z=-\dfrac{b_1\alpha^*}{2d}\lambda^k t+\tilde\phi_z(s,t,\omega), \qquad Y=\tilde\phi_y(s,t,\omega),
\end{equation}
where we further find from   \eqref{eq:t_deri} that
$$|\tilde\phi_{z,y}|+\left\|\frac{\partial\tilde\phi_{z,y}}{\partial(s,\omega)}\right\|=O(\hat\lambda^j), \qquad
\frac{\partial\tilde\phi_{z,y}}{\partial t}= \lambda^k O(s^2) + \lambda^k O(\lambda^j).$$
Equating $Z$ in this formula with $Z$ in \eqref{zy7a} gives
\begin{equation}\label{seqt}
s=-\lambda^j\frac{\alpha^*A_1+\beta^* A_2}{b_1\alpha^*} -\dfrac{b_1\alpha^*}{2d}\lambda^k t +\sigma(\omega,t),
\end{equation}
where, using \eqref{eq:j} and \eqref{eq:lambdagamma}
$$|\sigma|+\left|\frac{\partial\sigma}{\partial\omega}\right|=O(\hat\lambda^j), \qquad 
\frac{\partial\sigma}{\partial t} = 
\lambda^k O(\hat \lambda^j).$$

Substitution of \eqref{seqt} and the $Y$ equation of \eqref{zy8} into the $Y$ equation of \eqref{zy7a} gives the following relation between $\omega$ and $t$ which corresponds to the sought intersection of $T_1(S_j)$ with $W^u_{loc}(Q)$:
$$A_3 =- (\alpha^*A_1+\beta^* A_2)\lambda^k t + g(t,\omega)$$
where, using \eqref{eq:j} and \eqref{eq:lambdagamma}, we obtain
\begin{equation}\label{gesteem}
|g| + \left|\frac{\partial g}{\partial\omega}\right| = O(\lambda^{-j}\hat\lambda^j), \qquad  \frac{\partial g}{\partial t}= 
\lambda^k O(\lambda^{-j}\hat\lambda^j).
\end{equation}
Finally, by \eqref{eq:AB}, we  get
\begin{equation}\label{nontr:7aa}
(C_1(t,\omega) x_1^*+C_2(t,\omega)x_2^*)\cos j\omega+(C_2(t,\omega) x_1^*- C_1(t,\omega)x_2^* )\sin j\omega - 
g(t,\omega) = 0,
\end{equation}
where
\begin{equation}\label{ccap}
C_i=c_i + t \lambda^k(a_{i1} \alpha^* + a_{i2}\beta^*), \qquad i=1,2,
\end{equation}
with $\alpha^*$, $\beta^*$ given by \eqref{eq:maps:alphabeta}.

We may write equation \eqref{nontr:7aa} as
$$(c_1x_1^*+c_2x_2^*)\cos j\omega+(c_2 x_1^*- c_1 x_2^* )\sin j\omega = O(\lambda^{-j}\hat\lambda^j),$$
since $k$ is chosen such that \eqref{eq:j} is satisfied, implying that $O(\lambda^k)$ terms are absorbed by 
$O(\lambda^{-j}\hat\lambda^j)$; the first derivatives of the right-hand side are also estimates by 
$O(\lambda^{-j}\hat\lambda^j)$. By \eqref{eq:a_31} and \eqref{eq:xstar}, we have
\begin{equation*}\label{nontr:a_31}
(c_1^2+c_2^2)((x_1^*)^{2}+(x_2^*)^{2})\neq 0,
\end{equation*}
which, after sorting terms, is just
$$(c_1x_1^*+c_2x_2^*)^2+(c_2x_1^*-c_1x_2^* )^2\neq 0.$$
So, equation \eqref{nontr:7aa} can be further reduced to
\begin{equation}\label{nontr:15}
\sqrt{(c_1^2+c_2^2)((x_1^*)^{2}+(x_2^*)^{2})}\sin(j\omega+\varphi)=O(\lambda^{-j}\hat\lambda^j),
\end{equation}
where
\begin{equation}\label{varpeq}
\varphi=\arctan \dfrac{c_1x_1^*+c_2x_2^*}{c_2x_1^* - c_1x_2^* }. 
\end{equation}

Since the right-hand side of \eqref{nontr:15} is small with derivatives at large $j$, this equation can be resolved for all sufficiently large $k$ as  
\begin{equation}\label{nontr:17}
\omega_k(t;j)= j^{-1}(n_{k,j}\pi-\varphi)+O(\lambda^{-j}\hat\lambda^j),
\end{equation}
where $n_{k,j}$ can be any integer. In particular, we take 
$$n_{k,j}=\left[\frac{j\omega^*}{\pi}\right]+i_{k,j},$$
where $i_{k,j}$ is taken from a bounded set of integers described below. Since one can write $[{j\omega^*}/{\pi}]={j\omega^*}/{\pi}+O(1)$, substituting the above $n_{k,j}$ into \eqref{nontr:17} shows that   $\omega_k(t;j)\to\omega^*$ as $j\to\infty$ and $k\to\infty$, as desired.

The choice of $i_{k,j}$ is determined  by the requirement that the value of
$\alpha^*$ stays away from zero for $\omega=\omega_k$ (this is our standing assumption underlying all the computations,
see \eqref{eq:alpha*}). By \eqref{eq:maps:alphabeta}, this is equivalent to $(k\omega_k \bmod \pi)$ staying bounded away
from $-\arctan(c_2/c_1)$. In fact (we will use it below), we can choose $i_{k,j}$ such that $(k\omega_k \bmod \pi)$ stays away
from any finite set of given values between zero and $\pi$. Indeed, let the number of values to avoid be $(M-1)$. Take $k=(MN+1) \left[j/M\right]$ with the integer $N$ large enough (so that \eqref{eq:j} holds). Then the consecutive values of $k\omega_k$ corresponding to $i_{k,j} =1, \ldots, M$ differ by $(N\pi + \pi/M +O(k^{-1}))$, i.e., they take $M$ different
values $\!\!\!\mod\pi$. Since the cardinality of the forbidden set is less than $M$, it follows that for at least one choice of $i_{k,j}$ the corresponding value of 
$k\omega_k \!\!\!\mod \pi$ stays bounded away from this set.

We have shown that for the values of $\omega(t)$ given by \eqref{nontr:17} there exists an intersection of 
$W^u(O)$ with $W^s(Q)$.
Since the derivative of the left-hand side of \eqref{nontr:15} with respect to $\omega$ is non-zero at $\omega=\omega_k$, it follows that the intersection of $W^u(O)$ with $W^s(Q)$ splits as $\omega$ changes. This intersection also unfolds when $\mu$ changes from $\mu_k$, since $d\mu_k/d\omega =k\lambda^k(\alpha^* x^+_1+\beta^*x^+_2)+o(\lambda^k)\neq 0$ by \eqref{eq:mu_Q}.
Thus, to finish the proof of the proposition, it remains to prove \eqref{omtd}.\\

\noindent\textbf{(4) Computing the derivative $d\omega_k/dt$.}
We have 
$$\frac{d}{dt}\omega_k(t;j)= - \left.\frac{\partial G}{\partial t} \middle/ \frac{\partial G}{\partial \omega}\right.,$$
where we denote as $G$ the left-hand side of \eqref{nontr:7aa}, and the partial derivatives are computed at 
$\omega=\omega_k(t;j)$.
By \eqref{nontr:15} and \eqref{nontr:17},
$$\frac{\partial G}{\partial \omega}= j\sqrt{(c_1^2+c_2^2)((x_1^*)^{2}+(x_2^*)^{2})}\cos(j\omega+\varphi)
+O(\lambda^{-j}\hat\lambda^j)$$
$$\qquad\qquad\qquad\qquad= (-1)^{n_{k,j}} j \sqrt{(c_1^2+c_2^2)((x_1^*)^{2}+(x_2^*)^{2})}+O(\lambda^{-j}\hat\lambda^j)\neq 0$$
at $\omega=\omega_k(t;j)$, and, by \eqref{gesteem}, \eqref{ccap}, \eqref{eq:j},
$$\frac{\partial G}{\partial t}= \lambda^k (\alpha^*A_1+\beta^* A_2) + \lambda^k O(\lambda^{-j}\hat\lambda^j).$$
where $A_{1,2}$ are given by \eqref{eq:AB} and $\alpha^*,\beta^*$ are given by \eqref{eq:maps:alphabeta}. Thus to prove \eqref{omtd} we need to show that the quantity
$$\alpha^*A_1+\beta^* A_2=
(c_1A_1 +c_2A_2)\cos k\omega_k + (c_2A_1-c_1A_2)\sin k\omega_k$$
stays bounded away from zero as $j,k\to \infty$.

As we mentioned, we can choose $\omega_k$ such that $k\omega_k\!\!\!\mod\pi$ is bounded away from any given value, so
to prove that $(c_1A_1 +c_2A_2)\cos k\omega_k + (c_2A_1-c_1A_2)\sin k\omega_k$ is bounded away from zero, it is enough to show that 
\begin{enumerate}[nosep]
\item $A_1$ and $A_2$ tend to a limit as $j,k\to \infty$, and 
\item $(c_1A_1 +c_2A_2)^2 + (c_2A_1-c_1A_2)^2$ stays bounded away from zero.
\end{enumerate}
 By \eqref{eq:a_31},
it is the same as to show that $A_1^2+A_2^2$ tends to a non-zero value.  By \eqref{varpeq} and \eqref{nontr:17},
$$\lim_{j,k\to\infty}\tan (j\omega_k) =- \frac{c_1x_1^*+c_2x_2^*}{c_2x_1^* - c_1x_2^*},$$
so, 
$$
\lim_{j,k\to\infty} \sin^2(j\omega_k) = \dfrac{(c_1 x_1^*+c_2x_2^*)^2}{(c_1^2+c_2^2)((x^*_1)^2+(x^*_2)^2)}
\quad\mbox{and}\quad
\lim_{j,k\to\infty} \cos^2(j\omega_k) = \dfrac{(c_2 x_1^*-c_1x_2^*)^2}{(c_1^2+c_2^2)((x^*_1)^2+(x^*_2)^2)}.
$$
By \eqref{eq:AB}, one has
\begin{align*}
\lim_{j,k\to\infty} \sqrt{c_1^2+c_2^2}\; |A_1|= |c_2 a_{11} - c_1a_{12}| \sqrt{(x_1^*)^2+(x_2^*)^2},\\
\lim_{j,k\to\infty} \sqrt{c_1^2+c_2^2}\; |A_2|= |c_2 a_{21} - c_1a_{22}| \sqrt{(x_1^*)^2+(x_2^*)^2}.
\end{align*}
One sees from \eqref{eq:x+2} and \eqref{eq:det} that the second limit is non-zero by taking the transverse homoclinic point sufficiently close to $M^+$ so that $(x^*_1,x^*_2)$ is sufficiently close to $(x^+_1,x^+_2)\neq 0$. 
\end{proof}

We now verify a transversality condition which will be discussed in Section \ref{sec:HDCcon}.
Let $W_{loc}^{sE}(Q)$ be the extended stable manifold defined in Section \ref{sec:Wse}. For brevity, let us write $\mu_k(t)$ for $\mu_k(\omega_k(t),t)$ given by Proposition \ref{prop:nontr}

\begin{lem}\label{lem:gc1}
At $(\mu,\omega)=(\mu_k(t),\omega_k(t))$, for $t$ lying in the interval $(t_-,t_+)$ given by \eqref{tpm} and bounded away from zero, the invariant manifold $W^u(O)$ is transverse to $W_{loc}^{sE}(Q)$ at the point of their intersection.
\end{lem}

\begin{proof}
Recall that the intersection point $P$ of $W^s(Q)\cap W^u(O)$ is found as the intersection $W^s_{loc}(Q)$ with $T_1(S_j)$, where $S_j$ is a piece of
 $W^u(O)$. Below we find the tangent spaces of $W^{sE}_{loc}(Q)$ and $T_1(S_j)$ at $P$.

Since $P$ lies on $W^s_{loc}(Q)$, it follows from \eqref{eq:WsQ:0.2} that its $Y$-coordinate is $Y_Q+O(\lambda^k)$. Thus, Lemma \ref{lem:WseQ} guarantees that for all sufficiently large $k$ the tangent vectors to $W^{sE}_{loc}(Q)$ at $P$ lie in the cone \eqref{conewse}.

On the other hand, since $\alpha^*A_1+\beta^*A_2\neq 0$ at the moment of intersection (see the end of the proof of Proposition \ref{prop:nontr}), we find that $s\lambda^{-j}$ in \eqref{seqt} is bounded away from zero and infinity. It follows from formula \eqref{nontr1} that $T_1(S_j)$ for such $s$ is the graph of a function $(Y,V)\mapsto (Z,W)$ such that
$$\frac{\partial (Z,W)}{\partial (Y,V)}=o(1)_{k,j\to+\infty},$$
so it is transverse to the cone \eqref{conewse} when $t$ is bounded away from zero.
\end{proof}

\subsection{Transverse intersection of $W^u(Q)$ with $ W^s(O)$}\label{sec:trans}
Theorem \ref{thm:sfsimple} immediately follows from Proposition \ref{prop:nontr} and the statement below.

\begin{prop}\label{prop:tran} Let $\{f_{\mu,\omega}\}$ be a proper family.
If the homoclinic tangency at $\mu=0$ is accompanied, then the unstable manifold $W^u(Q)$ of the periodic point $Q$ obtained at $\mu=\mu_k$ for all sufficiently large $k$ in Lemma \ref{lem:index2} has a transverse intersection with $W^s(O)$.
\end{prop}

\begin{proof} 
Since the homoclinic tangency is accompanied, one can choose $\delta>0$ such that $W^s(O)$ intersects $W^u_{loc}(O)$ 
at a point $M^*_+=(x=0,u=0,y=y^*_+,v=v^*_+)$ such that  
\begin{equation}\label{eq:tran:4+}
y^*_+ -y^-\in (2K_1\delta^3,\dfrac{K_2}{2}\delta^2)
\quad\mbox{and}\quad
\|v^*_+ - v^-\|<\dfrac{1}{2}\delta,
\end{equation}
where $y^-$ is the $y$-coordinate of the homoclinic tangency point $M^-$, and $K_{1,2}$ are constants from 
Lemma \ref{lem:WuQ}.
We denote the small connected piece of $W^s(O)$ through $M^*_+$ as $W^{s+}$.
 Also, possibly choosing a different small
value for $\delta$, one finds an intersection point $M^*_-=(x=0,u=0,y=y^*_-,v=v^*_-)$ of the intersection of $W^s(O)$ with $W^u_{loc}(O)$  such that  
\begin{equation}\label{eq:tran:4-}
y^*_- -y^-\in (-2K_1\delta^3,-\dfrac{K_2}{2}\delta^2)
\quad\mbox{and}\quad
\|v^*_- - v^-\|<\dfrac{1}{2}\delta,
\end{equation}
The 
corresponding small connected piece of $W^s(O)$ will be denoted as $W^{s-}$.

We establish the existence, for all sufficiently large $k$, of a transverse intersection of $W^u_{loc}(Q)$ with $T_0^{-k}(W^{s+})$ or
$T_0^{-k}(W^{s-})$. We will look for the intersections with $T_0^{-k}(W^{s+})$ for those $k$ for which the regions $D(V)$ in Lemma \ref{lem:WuQ} lie in $Y>0$, and the intersections with $T_0^{-k}(W^{s-})$
for those $k$ for which the regions $D(V)$ lie in $Y<0$. Since the computations are the same in both cases, 
we do them for the manifolds $T_0^{-k}(W^{s\pm})$ simultaneously.

By the transversality, if we denote as $(\tilde x,\tilde u, \tilde y, \tilde v)$ the coordinates near $M^-$, the equation of $W^{s\pm}$ is
$$(\tilde y,\tilde v)=\psi_\pm(\tilde x,\tilde u),$$
where $\psi_\pm$ is a smooth function, defined for small $(\tilde x,\tilde u)$, such that $\psi_\pm(0,0)=(y^*_\pm,v^*_\pm)$. The transverse intersection of $W^{s\pm}$ with $W^u_{loc}$ persists at small changes of the parameters, implying that the equations for $W^{s\pm}$ keep
their form, with $\psi_\pm$ depending smoothly on the parameters as well, and the inequalities \eqref{eq:tran:4+} or \eqref{eq:tran:4-} also persist.

It follows then from \eqref{eq:maps:T0k} that if $(x,u,y,v)=T_0^{-k}(\tilde x,\tilde u, \tilde y, \tilde v)$ for a point $(\tilde x,\tilde u, \tilde y, \tilde v)\in W^{s\pm}$, then
\begin{equation}\label{referred:3}
(\tilde y,\tilde v) = (y^*_\pm,v^*_\pm) + \tilde\psi_\pm(x,u),
\end{equation}
where $\tilde\psi_\pm= O(\lambda^k)$ along with its derivative with respect to $(x,u)$.
Going to the $(Z,Y,W,V)$-coordinates used in Lemma \ref{lem:WuQ} (defined by the composition of the coordinate transformations \eqref{eq:maps:coor1}, \eqref{eq:maps:coor2}, \eqref{eq:maps:coor2.1}, \eqref{eq:maps:coor3}), we obtain from \eqref{referred:3} the equation for $T_0^{-k}W^{s\pm}$ in the form
\begin{equation}\label{eq:tran:3}
(Y,V)=(y^* , v^*) + \hat\psi_\pm(Z,W),
\end{equation}
where 
$$
y^*=y^*_\pm-y^-\in \pm(2K_1\delta^3,\dfrac{K_2}{2}\delta^2)
\quad\mbox{and}\quad
\|v^*\|=\|v^*_\pm-v^- \|+ O(y^*_+-y^-)<\dfrac{2}{3}\delta$$
 by \eqref{eq:tran:4+}, and  the smooth function $\hat\psi_\pm$ is defined for $(Z,W)\in [-\delta,\delta]\times[-\delta,\delta]^{d^s-1}$ and is estimated as
$O(\lambda^k)$ along with the derivative with respect to $(Z,W)$. 

Thus, there exists an interval of $Z$ values for which the points with $(Y=y^*,V=v^*)$ lie in the interior of the definition domain $D(v^*)$ of the function $\phi^u(Z,Y,V)$ in Lemma \ref{lem:WuQ}. Since $\hat\psi$ is $C^1$-small at all large $k$,
it follows that the system of equations comprised by 
the equation \eqref{eq:tran:3} for $T_0^{-k}(W^{s\pm})$ and the equation
$$W=W_Q+ \phi^u(Z,Y,V)$$
for $W^u_{loc}(Q)$ has a unique solution $(Y,V,W)$ for all $Z$ from this interval. This gives the sought curve of 
the transverse intersection of $T_0^{-k}(W^{s+})$ or $T_0^{-k}(W^{s-})$  with $W^u_{loc}(Q)$.
\end{proof}

\begin{rem}\label{lem:gc2}
The proof shows that the $Y$-coordinate of the points in the
intersection curve $T^{-k}(W^{s+})\cap W^u_{loc}(Q)$ or
$T^{-k}(W^{s-})\cap W^u_{loc}(Q)$ satisfy $|Y|\in(2K_1\delta^3, K_2\delta^2/2)$.
It then follows immediately from Lemma \ref{lem:WuuQ} that when $t\in(t_-,t_+)$ is bounded away from zero and
condition \eqref{tsign} holds, the strong unstable manifold $W^{uu}_{loc}(Q)$ intersects this curve,
and is transverse to $T^{-k}(W^{s+})$ or $T^{-k}(W^{s-})$ given by \eqref{eq:tran:3}.
\end{rem}

\subsection{Accumulating transverse homoclinic points: proofs of Lemmas \ref{lem:tangency1} and \ref{lem:tangency2}}\label{sec:twolem}

\subsubsection{Proof of Lemma \ref{lem:tangency1}.}
We will use the formulas for $T_0$ and $T_1$ in the original coordinate $(x_1,x_2,y,u,v)$.  Let $l_1=\tilde y-y^-,\ l_2=\tilde v-v^-$. According to  \eqref{eq:maps:T1_cross}, the image $T_1(W^u_{loc}(O)\cap \Pi^-)$ is given by
\begin{equation}\label{eq:maps:T1_cross2}
\begin{aligned}
x_1 - x_1^+ &= b_1 l_1 + a_{15} v +O(l_1^2+v^2), \\
x_2 - x_2^+ &=  a_{25} v +O(l_1^2+v^2),  \\
y &= \mu  +d l_1^2+ c_5v+O(l_1^3)+O(|l_1v|+v^2),  \\
u - u^+ &=b_4l_1+  a_{45} v +O(l_1^2+v^2),  \\
l_2 &=b_5l_1+  a_{55}v +O(l_1^2+v^2).
\end{aligned}
\end{equation}
Since $W^s_{loc}(O)$ is give by $\{y=0,v=0\}$, a point $(x_1,x_2,y,u,v)\in T_1(W^u_{loc}(O)\cap \Pi^-)$ is a homoclinic point of $O$ if 
\begin{equation}\label{eq:tan1:1}
T_0^k(x_1,x_2,y,u,v)=(x'_1,x'_2,y',u',v')
\quad\mbox{and}\quad
T_1(x'_1,x'_2,y',u',v')=(\bar x_1,\bar x_2, 0,\bar u, 0)
\end{equation}
for some $(x'_1,x'_2,y',u',v')\in\Pi^-$,
namely, the point corresponds to a solution $(l_1, v)$ to the system consisting of  \eqref{eq:maps:T1_cross2} and \eqref{eq:tan1:1}. It follows from \eqref{eq:maps:T0k} and \eqref{eq:maps:T1_cross} that these coordinates satisfy
\begin{equation}\label{eq:maps:T0knew}
\begin{aligned}
 &x'_{1}= \lambda^k x_1\cos k\omega  - \lambda^k x_2\sin k\omega + O(\hat\lambda^k),\\ 
& x'_{2} = \lambda^k x_1  \sin k\omega  + \lambda^k x_2 \cos k\omega  + O(\hat\lambda^k),\\ 
& y   = \gamma^{-k} y' +  O(\hat\gamma^{-k}),\quad
 u'   = O(\hat{\lambda}^k) ,\quad
 v  = O(\hat \gamma^{-k}),
\end{aligned}
\end{equation}
and
\begin{equation}\label{eq:maps:T1_crossnew}
\begin{aligned}
\bar{x}_1 - x_1^+ &=a_{11} x'_1 +a_{12} x'_2 +b_1(y' - y^-) + a_{14} u' + a_{15}  \bar{v} + \dots, \\
\bar{x}_2 - x_2^+ &=a_{21} x'_1 +a_{22} x'_2  + a_{24} u' + a_{25}  \bar{v} + \dots, \\
\bar y &= \mu + c_1 x'_1 +c_2 x'_2 +d(y' - y^-)^2+ c_4 u' + c_5  \bar{v} + \dots, \\
\bar{u} - u^+ &=a_{41} x'_1 +a_{42} x'_2 +b_4(y' - y^-) + a_{44} u' + a_{45}  \bar{v} + \dots ,\\
v' - v^- &=a_{51} x'_1 +a_{52} x'_2 +b_5(y' - y^-) + a_{54}u' + a_{55}  \bar{v} + \dots ,
\end{aligned}
\end{equation}
where the dots denote quadratic terms in the Taylor expansions, except the $(y' - y^-)^2$ term in the third equations. 

Let us solve the above-mentioned system. First, set $\bar{v}=0$ in the fifth equation of \eqref{eq:maps:T1_crossnew}, and substitute it into the $y$- and $v$-equations in \eqref{eq:maps:T0knew}. This gives
\begin{equation}\label{eq:tan1:3}
v=O(\hat\gamma^{-k})
\end{equation}
as a function of $x_1,x_2,y,u$. Next, we solve for $l_1$ from the first equation of \eqref{eq:maps:T1_cross2} as
\begin{equation}\label{514}
l_1=b^{-1}_1(x_1-x^+_1)-b^{-1}_1a_{15}v+O((x_1-x^+_1)^2+v^2).
\end{equation}
Combining this with \eqref{eq:tan1:3} and the remaining equations of \eqref{eq:maps:T1_cross2}, we obtain
\begin{equation}\label{eq:tan1:2}
\begin{aligned}
l_1&=b^{-1}_1(x_1-x^+_1)+O((x_1-x^+_1)^2)+O(\hat\gamma^{-k}),\\
x_2&=x_2^++O((x_1-x_1^+)^2)+O(\hat\gamma^{-k}),\\
y&=\mu+\dfrac{d}{b_1^2}(x_1-x_1^+)^2+O((x_1-x_1^+)^3)+O(\hat\gamma^{-k}),\\
u&=u^++\dfrac{b_4}{b_1}(x_1-x_1^+)+O((x_1-x_1^+)^2)+O(\hat\gamma^{-k}),\\
l_2 &= \dfrac{b_5}{b_1}(x_1-x_1^+)+O((x_1-x_1^+)^2)+O(\hat\gamma^{-k}),
\end{aligned}
\end{equation}
where the $O(\hat\gamma^{-k})$ terms are functions of $x_1$.

After setting 
\begin{equation}\label{515}
X=x_1-x_1^+ \quad \mbox{and}\quad  Y=y'-y^-,
\end{equation}
the third equation of \eqref{eq:tan1:2} and the $y$-equation of \eqref{eq:maps:T0knew} gives
\begin{equation*}
Y=y'-y^-=\gamma^k y+o(1)_{k\to\infty}-y^-=\gamma^k\mu -y^-+\dfrac{\gamma^kd}{b_1^2}X^2+O(\gamma^kX^3)+o(1)_{k\to \infty},
\end{equation*}
or,
\begin{equation}\label{eq:52}
0=\mu-\gamma^{-k}y^- -\gamma^{-k}Y+\frac{d}{b_1^2}X^2+O(X^3)+o(\gamma^{-k}).
\end{equation}

To find $X$ we need one more equation. Note that the $y$-coordinates $\bar y$ of points on $T_1\circ T_0^k\circ T_1(W^u_{loc}(O))$ can be found by directly combining \eqref{eq:maps:T0knew}, \eqref{eq:maps:T1_crossnew} with using \eqref{eq:maps:T1_cross2} and \eqref{eq:maps:alphabeta}. Alternatively, it can be found by using the $Y$-equation in \eqref{eq:maps:Tk_1} and then returning to the original coordinates by the \eqref{eq:maps:coor1}. In either way, we obtain

\begin{equation*}
\bar y=\mu+\alpha^*\lambda^kx_1^++\beta^*\lambda^kx_2^+ +\alpha^*\lambda^k X +d Y^2+O(\lambda^k X^2)+o(\lambda^k)+O(Y^3),
\end{equation*} 
where the left hand side is a function of $X,Y,\bar v$. For points in $W^s_{loc}(O)$, namely, those with $\bar y=\bar v=0$, one further has
\begin{equation}\label{eq:51}
0=\mu+\alpha^*\lambda^kx_1^++\beta^*\lambda^kx_2^+ +\alpha^*\lambda^k X +d Y^2+O(\lambda^k X^2)+o(\lambda^k)+O(Y^3).
\end{equation} 
The sought transverse homoclinic points now correspond to the non-degenerate solutions to the system consisting of \eqref{eq:52} and  \eqref{eq:51}. 

Let us find the solutions. Note that there are infinitely many $k\in 2\mathbb{N}$ such that $d(\alpha^* x_1^++\beta^*x_2^+)<0$. To see this, use \eqref{eq:maps:alphabeta} to write
\begin{equation}\label{eq:q}
\begin{aligned}
\alpha^* x_1^++\beta^*x_2^+
=\sqrt{\sigma_1^2+\sigma_2^2}\sin (k\omega +\eta)
\end{aligned}
\end{equation}
where $\sigma_1=c_1x_1^+ +c_2x_2^+$ and $\sigma_2=-c_1x_2^+ +  c_2x_1^+$. We need to show that $\sin (k\omega +\eta)$ with $k\in 2\mathbb N$ can have different signs, and it suffices to show that $k'\omega \bmod \pi$ with $k'\in \mathbb N$ can take at least three different values in $[0,1)$. But this is guaranteed when $\omega/\pi$ is irrational or equals $p/q$ with $q\geqslant 3$.

We take such $k\in 2\mathbb N$ and apply the scaling
\[(X,Y)\mapsto \left(b_1\sqrt{\frac{y^-}{d}}\gamma^{-\frac{k}{2}}U,\sqrt{\frac{\alpha^* x_1^+ +\beta^*x_2^+}{-d}}\lambda^{\frac{k}{2}}V\right).\] 
Taking $\mu=0$, we rewrite \eqref{eq:52} and \eqref{eq:51} as
\begin{align*}
1=U^2+o(1)_{k\to \infty},\qquad
1=V^2+o(1)_{k\to\infty}.
\end{align*}
Obviously, the solutions are $(U,V)=(\pm 1+o(1)_{k\to \infty}, \pm 1+o(1)_{k\to \infty})$, which lead to four solutions in the $(X,Y)$-variables:
\begin{equation}\label{eq:transcoor}
\begin{aligned}
(X_k^1,Y_k^1)&=\left(b_1\gamma^{-\frac{k}{2}}\sqrt{\frac{y^-}{d}}+ o(\gamma^{-\frac{k}{2}}),\lambda^{\frac{k}{2}}\sqrt{\frac{\alpha^* x_1^+ +\beta^*x_2^+}{-d}} +o(\lambda^{\frac{k}{2}})\right),\\
(X_k^2,Y_k^2)&=\left(b_1\gamma^{-\frac{k}{2}}\sqrt{\frac{y^-}{d}}+ o(\gamma^{-\frac{k}{2}}),-\lambda^{\frac{k}{2}}\sqrt{\frac{\alpha^* x_1^+ +\beta^*x_2^+}{-d}} +o(\lambda^{\frac{k}{2}})\right),\\
(X_k^3,Y_k^3)&=\left(-b_1\gamma^{-\frac{k}{2}}\sqrt{\frac{y^-}{d}}+ o(\gamma^{-\frac{k}{2}}),\lambda^{\frac{k}{2}}\sqrt{\frac{\alpha^* x_1^+ +\beta^*x_2^+}{-d}} +o(\lambda^{\frac{k}{2}})\right),\\
(X_k^4,Y_k^4)&=\left(-b_1\gamma^{-\frac{k}{2}}\sqrt{\frac{y^-}{d}}+ o(\gamma^{-\frac{k}{2}}),-\lambda^{\frac{k}{2}}\sqrt{\frac{\alpha^* x_1^+ +\beta^*x_2^+}{-d}} +o(\lambda^{\frac{k}{2}})\right).
\end{aligned}
\end{equation}
All these solutions are non-degenerate, since the Jacobian of the system given by \eqref{eq:52} and \eqref{eq:51} at these solutions is 
$$\dfrac{4d^2}{b_1^2}XY+o(\lambda^{\frac{k}{2}}\gamma^{-\frac{k}{2}})=\pm \dfrac{4d}{b_1}\lambda^{\frac{k}{2}}\gamma^{-\frac{k}{2}}\sqrt{-y^-(\alpha^*x_1^++\beta^*x_2^+)}+o(\lambda^{\frac{k}{2}}\gamma^{-\frac{k}{2}}),$$
which is non-zero.

The corresponding transverse homoclinic points in $T_1(W^u_{loc}(O))$ are given by
\begin{align*}
N_k^1&=(x_1^++\tilde{X}+o(\gamma^{-\frac{k}{2}}),x_2^1,(\tilde{Y}+y^-)\gamma^{-k}+o(\gamma^{-k}),u^1,v^1),\\
N_k^2&=(x_1^++\tilde{X}+o(\gamma^{-\frac{k}{2}}),x_2^2,(-\tilde{Y}+y^-)\gamma^{-k}+o(\gamma^{-k}),u^2,v^2),\\
N_k^3&=(x_1^+-\tilde{X}+o(\gamma^{-\frac{k}{2}}),x_2^3,(\tilde{Y}+y^-)\gamma^{-k}+o(\gamma^{-k}),u^3,v^3),\\
N_k^4&=(x_1^+-\tilde{X}+o(\gamma^{-\frac{k}{2}}),x_2^4,(-\tilde{Y}+y^-)\gamma^{-k}+o(\gamma^{-k}),u^4,v^4),
\end{align*}
where $\tilde{X}=b_1\gamma^{-\frac{k}{2}}\sqrt{\dfrac{y^-}{d}}$ and $\tilde{Y}=\lambda^{\frac{k}{2}}\sqrt{\dfrac{\alpha^* x_1^+ +\beta^*x_2^+}{-d}}$, and the coordinates $x_2^j,u^j,v^j$ $(j=1,2,3,4)$  not written explicitly can be found from \eqref{eq:tan1:2}. Denote $\hat{N}_k^j=T_1^{-1}(N_k^j)$. Since all these points lie on $W^u_{loc}(O)$, their coordinates are $(0,0,\hat{y}_k^j,0,\hat{v}_k^j)$, which, by the first and last equations of \eqref{eq:tan1:2}, are given by
\begin{align*}
\hat{y}_k^j-y^-=\dfrac{X_k^j}{b_1}+O((X_k^j)^2)+O(\hat{\gamma}^{-k})\quad\mbox{and}\quad
\hat{v}_k^j - v^- = \dfrac{b_5}{b_1}X_k^j+O((X_k^j)^2)+O(\hat\gamma^{-k}).
\end{align*}
Thus, we have $\hat{N}_k^{j}\to M^-$ as $k\to\infty$ and $\hat{y}^{1,2}_k$ and $\hat{y}^{3,4}_k$ lie on different sides of $y^-$.

\subsubsection{Proof of Lemma \ref{lem:tangency2}.}
We will first find new homoclinic tangencies which has $dy^->0$, then show that they unfold as $\mu$ varies, and finally prove that these tangencies are generic, namely, conditions C1--C5 are satisfied. The proof is accordingly divided into three parts.

\noindent {\bf (1) Creation of secondary homoclinic tangencies.} As discussed in the proof of Lemma \ref{lem:tangency1}, homoclinic points in $W_{loc}^u(O)$ with the global map $T_1\circ T_0^k\circ T_1$ correspond to solutions to the system consisting of \eqref{eq:52} and \eqref{eq:51}. Particularly, quadratic homoclinic tangencies corresponds to solutions with multiplicity two. In what follows we find such solutions.

On the one hand, after eliminating $\mu$ from \eqref{eq:52} and \eqref{eq:51}, we have
\begin{equation}\label{eq:505}
\lambda^k(\alpha^*x_1^++\beta^* x_2^+)+\gamma^{-k}y^-+\alpha^*\lambda^kX+\gamma^{-k}Y-\dfrac{d}{b_1^2}X^2+dY^2+O(\lambda^kX^2)+O(X^3)+O(Y^3)+o(\lambda^k)=0.
\end{equation}
On the other hand, solutions of multiplicity two implies that the Jacobian of the system vanishes. So, we have
\begin{align}\label{eq:56}
\mathrm{Jacobian} =&-\gamma^{-k}(\alpha^* \lambda^k+O(\lambda^kX)+o(\lambda^k))-\dfrac{4d^2}{b_1^2}\big(Y+O(Y^2)+o(\lambda^k)\big)\big(X+O(X^2)+o(\gamma^{-k})\big)\nonumber \\
=&-\alpha^* \lambda^k\gamma^{-k}-\dfrac{4d^2}{b_1^2}\hat{X}\hat{Y}+O(\lambda^k\gamma^{-k}\hat{X})+o(\lambda^k\gamma^{-k})=0,
\end{align}
where we denoted the last two brackets in the first line by
\begin{equation}\label{eq:571}
\hat{Y}:=Y+O(Y^2)+o(\lambda^k) \quad  \mbox{and} \quad  \hat{X}:=X+O(X^2)+o(\gamma^{-k}).
\end{equation}
With these notations, \eqref{eq:505} can be recast as 
\begin{equation}\label{eq:55}
\lambda^k(\alpha^*x_1^++\beta^* x_2^+)+\gamma^{-k}y^-+\alpha^*\lambda^k\hat{X}+\gamma^{-k}\hat{Y}-\dfrac{d}{b^2_1}\hat{X}^2+d \hat{Y}^2+O(\hat{X}^3)+O(\hat{Y}^3)+o(\lambda^k)=0.
\end{equation}

We now solve the system consisting of \eqref{eq:56} and \eqref{eq:55}. By taking $k$ such that $d(\alpha^* x_1^++\beta^* x_2^+)<0$ (as we did in the proof of Lemma \ref{lem:tangency1}), we apply the scaling
\begin{equation}\label{eq:57}
(\hat{X},\hat{Y})\mapsto \sqrt{\dfrac{-(\alpha^* x_1^++\beta^* x_2^+)}{d}}\left(\dfrac{-\alpha^* b_1^2}{4d(\alpha^* x_1^++\beta^* x_2^+)} \lambda^{\frac{k}{2}}\gamma^{-k}U,\lambda^{\frac{k}{2}}V \right),
\end{equation}
which takes \eqref{eq:56} and \eqref{eq:55} to (recall that $|\lambda\gamma|>1$)
\begin{align*}
1+UV =o(1)_{k\to \infty},\qquad
1-V^2=o(1)_{k\to \infty}
\end{align*}
with solutions $(U^1,V^1)=(1+o(1)_{k\to \infty},- 1+o(1)_{k\to \infty})$ and $(U^2,V^2)=(-1+o(1)_{k\to \infty},1+o(1)_{k\to \infty})$.

Returning to the $(X,Y)$-variables, we have  
\begin{equation}\label{eq:coornewHT}
\begin{aligned}
(\hat{X}_k^1,\hat{Y}_k^1)&=\sqrt{\dfrac{-(\alpha^* x_1^++\beta^* x_2^+)}{d}}\left(\dfrac{-\alpha^* b_1^2\lambda^{\frac{k}{2}}\gamma^{-k}}{4d(\alpha^* x_1^++\beta^* x_2^+)}+o(\lambda^{\frac{k}{2}}\gamma^{-k}),-\lambda^{\frac{k}{2}}+o(\lambda^{\frac{k}{2}})\right),\\
(\hat{X}_k^2,\hat{Y}_k^2)&=\sqrt{\dfrac{-(\alpha^* x_1^++\beta^* x_2^+)}{d}}\left(\dfrac{\alpha^* b_1^2\lambda^{\frac{k}{2}}\gamma^{-k}}{4d(\alpha^* x_1^++\beta^* x_2^+)}+o(\lambda^{\frac{k}{2}}\gamma^{-k}),\lambda^{\frac{k}{2}}+o(\lambda^{\frac{k}{2}})\right).
\end{aligned}
\end{equation}
It is easy to check that these solutions are non-degenerate so they indeed have multiplicity two. By \eqref{eq:52}, the corresponding $\mu$ values to the solutions have the same form as 
\begin{equation}\label{eq:tan2:mu}
\mu^j_k=\gamma^{-k}y^-(1+o(1)_{k\to \infty}), \quad j=1,2.
\end{equation}
By \eqref{515} and \eqref{eq:571}, the tangency points on $T_1(W^u_{loc}(O))$ are given by
\begin{equation}\label{coorM}
M_k^j=\big(\hat{X}_k^j+x_1^++o(\gamma^{-k}),x_{2,k}^j, \gamma^{-k}(\hat{Y}_k^i+y^-+o(\lambda^k)),u_k^j,v_k^j\big)
\end{equation}
with $x_{2,k}^j,\ u_k^j,\ v_k^j$ obtained from \eqref{eq:tan1:3} and \eqref{eq:tan1:2}.

Denote $M_k^{-,j}=T_1^{-1}(M_k^j)$. Since these points lie on $W^u_{loc}(O)$, they have coordinates $(0,0,\hat{y}_k^{-,j},0,\hat{v}_k^{-,j})$. It follows from the first and fifth equations of \eqref{eq:tan1:2} that 
\begin{equation}\label{eq:65}
\hat{y}_k^{-,j}-y^-=\dfrac{\hat{X}_k^j}{b_1}+O((\hat{X}_k^j)^2)+o(\hat{\gamma}^{-k}) \quad\mbox{and}\quad
\hat{v}_k^{-,j}-v^-=O(\lambda^{\frac{k}{2}}\gamma^{-k}),
\end{equation} 
which shows that $M_k^-\to M^-$ as $k\to\infty$. 

Now, let us check that at least one of these homoclinic tangency points satisfies $dy^->0$. Recall that the global map corresponding to the secondary tangency point $M^{-,j}_k$ is $T_1\circ T_0^k\circ T_1$. Let us denote the corresponding ``$d$'' by $d_k^j$. By \eqref{eq:maps:T1_crossnew}, we have
\begin{align}\label{eq:58}
d_{k}^j&=\dfrac{1}{2}\dfrac{\partial^2 \bar{y}}{\partial \tilde{y}^2}(0,0,\hat{y}^{-,j}_k,0,\hat{v}^{-,j}_k) \nonumber\\
&=\dfrac{1}{2}\dfrac{\partial }{\partial \tilde{y}}\left(\dfrac{\partial \bar{y}}{\partial x'_{1}}\dfrac{ \partial x'_{1}}{\partial \tilde{y}}+ \dfrac{\partial \bar{y}}{\partial x'_{2}}\dfrac{ \partial x'_{2}}{\partial \tilde{y}} +\dfrac{\partial \bar{y}}{\partial y'}\dfrac{ \partial y'}{\partial \tilde{y}} +\dfrac{\partial \bar{y}}{\partial u'}\dfrac{ \partial u'}{\partial \tilde{y}}  +\dfrac{\partial \bar{y}}{\partial \bar{ v}}\dfrac{ \partial \bar{v}}{\partial \tilde{y}}\right)\Big|_{(0,0,\hat{y}^{-,j}_k,0,\hat{v}^{-,j}_k)}
\end{align}
with
\begin{equation}\label{eq:508}
\begin{aligned}
\dfrac{\partial \bar{y}}{\partial y'}\Big|_{T_0^k\circ T_1(0,0,\hat{y}^{-,j}_k,0,\hat{v}^{-,j}_k)}
&=2d(y'-y^-)+O((y'-y^-)^2)+O(\lambda^k)\\
&=2d\hat Y^j_k+O((\hat Y^j_k)^2)+O(\lambda^k),\\[10pt]
\dfrac{\partial^2 \bar{y}}{\partial y'^2}\Big|_{T_0^k\circ T_1(0,0,\hat{y}^{-,j}_k,0,\hat{v}^{-,j}_k)}
&=2d +O(\hat Y^j_k) +O(\lambda^k),
\end{aligned}
\end{equation}
where we use $\hat{Y}^j_k=y'-y^-+O((y'-y^-)^2)+o(\lambda^k) $ by \eqref{515} and \eqref{eq:571}.

To calculate ${\partial y'}/{\partial \tilde{y}}$, 
one needs to consider the composition $ T_0^k\circ T_1$. Thus, we rewrite  \eqref{eq:maps:T1_crossnew} with the appropriate notations as
\begin{equation*}
\begin{aligned}
{x}_1 - x_1^+ &=a_{11} \tilde x_1 +a_{12} \tilde x_2 +b_1(\tilde y- y^-) + a_{14}\tilde  u + a_{15}  {v} + \dots, \\
{x}_2 - x_2^+ &=a_{21} \tilde x_1 +a_{22} \tilde x_2  + a_{24} \tilde u + a_{25}  {v} + \dots, \\
 y &= \mu + c_1 \tilde x_1 +c_2 \tilde x_2 +d(\tilde y - y^-)^2+ c_4 \tilde u + c_5  {v} + \dots, \\
{u} - u^+ &=a_{41} \tilde x_1 +a_{42} \tilde x_2 +b_4(\tilde y- y^-) + a_{44} \tilde u+ a_{45}  {v} + \dots ,\\
v' - v^- &=a_{51} \tilde x_1 +a_{52} \tilde x_2 +b_5(\tilde y- y^-) + a_{54}\tilde u + a_{55}  {v} + \dots ,
\end{aligned}
\end{equation*}
and combine it  with \eqref{eq:maps:T0knew} to get
\begin{equation}\label{eq:60}
\begin{aligned}
\dfrac{\partial y' }{\partial \tilde{y}}\Big|_{(0,0,\hat{y}^{-,j}_k,0,\hat{v}^{-,j}_k)}
&=\dfrac{\partial y' }{\partial x_{1}}\dfrac{ \partial x_{1}}{\partial \tilde{y}}+\dfrac{\partial y'}{\partial x_{2}}\dfrac{ \partial x_{2}}{\partial \tilde{y}}+ \dfrac{\partial y' }{\partial y}\dfrac{ \partial y}{\partial \tilde{y}} +\dfrac{\partial y'}{\partial u}\dfrac{ \partial u}{\partial \tilde{y}}+\dfrac{\partial y' }{\partial v'}\dfrac{ \partial v'}{\partial \tilde{y}}\\
&=2d\gamma^k(\tilde{y}-y^-)+o((\tilde{y}-y^-)^2)+o(1)_{k\to\infty}\\
&=\dfrac{2d\gamma^k}{b_1} \hat{X}^j_{k}+o((\hat{X}^j_{k})^2)+o(1)_{k\to\infty},\\
\dfrac{\partial^2 y'}{\partial \tilde{y}^2}\Big|_{(0,0,\hat{y}^{-,j}_k,0,\hat{v}^{-,j}_k)}
&=2d\gamma^k+o(\hat{X}^j_{k})+o(1)_{k\to\infty},
\end{aligned} 
\end{equation}
where we use $\hat{X}^j_k=\tilde y-y^- + O((\tilde y-y^-)^2)+o(\gamma^{-k})$ by the first equation in \eqref{eq:tan1:2} and \eqref{eq:571}. 

Similarly, we obtain
\begin{equation}\label{eq:59}
 \dfrac{\partial x'_1}{\partial \tilde{y}}=O(\lambda^k), \dfrac{\partial x'_{2}}{\partial \tilde{y}}=O(\lambda^k),  \dfrac{\partial u'}{\partial \tilde{y}}=o(\lambda^k),
  \end{equation}
  where the derivatives are all evaluated at the point $(0,0,\hat{y}^{-,j}_k,0,\hat{v}^{-,j}_k)$.
Note that by definition $ \bar v$ is independent of $\tilde y$ and hence  ${\partial \bar v}/{\partial \tilde y}=0$ in \eqref{eq:58}.
Combining \eqref{eq:58} with \eqref{eq:508} -- \eqref{eq:59}  and using  \eqref{eq:coornewHT}, we obtain
\begin{align*}\label{eq:61}
d_{k}^j
=&\dfrac{1}{2}\left[\dfrac{\partial^2 \bar{y}}{\partial y'^2}\left(\dfrac{\partial y'}{\partial \tilde{y}}\right)^2+\dfrac{\partial \bar{y}}{\partial y' }\dfrac{\partial ^2 y' }{\partial \tilde{y}^2}\right]+O(\lambda^k)\nonumber\\
=&\big(d+O(Y^j_k)+O(\lambda^k)\big)\left(\dfrac{2d}{b_1}\gamma^k \hat{X}^j_{k}+o((\hat{X}^j_{k})^2)+o(1)_{k\to\infty}\right)^2\nonumber\\
&+d\hat{Y}^j_k\big(2d\gamma^k+o(\hat{X}^j_k)+o(1)_{k\to\infty}\big)+o(\lambda^k)\nonumber\\
=&
(-1)^{j}2d^2\lambda^{\frac{k}{2}}\gamma^k\sqrt{\dfrac{-(\alpha^* x_1^++\beta^* x_2^+)}{d}}\left(1+o(1)_{k\to\infty}\right),
\end{align*}
which implies that for one of $M_k^{-,1}$ and $M_k^{-,2}$ we have $d_{k}^j y^->0$. 

\noindent {\bf (2) Unfolding of secondary homoclinic tangencies.}
Denote $T_1^{new}=T_1\circ T_0^k\circ T_1$. For any system $g$ close to $f$, we define $\mu^{new}(g)$ as the functional measuring the (signed)
distance between $T_1^{new}(0,0,y^{-}(g),0,v^{-})$ and $W^s_{loc}(O)$. To show that the newly found tangencies unfolds as $\mu$ varies, it suffices to prove that for each sufficiently large $k$, the parameter $\mu(g)$ satisfies $$\dfrac{\partial \mu^{new}(f_{\mu})}{\partial \mu}\Big|_{\mu=\mu_k}\not=0,$$  
where $\mu_k$ are given by \eqref{eq:tan2:mu}.

Let us find the $y$-equation of $T_1^{new}$. Combining \eqref{eq:maps:T1_cross} and \eqref{eq:maps:T0knew}, yields
\begin{equation}\label{64}
\begin{aligned}
x'_1 &= \lambda^k(x_1^+\cos k\omega  - x_2^+\sin k\omega)+\lambda^k \cos k\omega [a_{11}\tilde{x}_1 +a_{12}\tilde{x}_2+b_1(\tilde{y}-y^-)]\\ &\quad -\lambda^k \sin k\omega [a_{21}\tilde{x}_1 +a_{22}\tilde{x}_2+a_{24}\tilde{u} +a_{25} v]+O(\hat{\lambda}^k),  \\
x'_2 &= \lambda^k(x_1^+\sin k\omega  + x_2^+\cos k\omega)+\lambda^k \sin k\omega [a_{11}\tilde{x}_1 +a_{12}\tilde{x}_2+b_1(\tilde{y}-y^-)]\\ &\quad +\lambda^k \cos k\omega [a_{21}\tilde{x}_1 +a_{22}\tilde{x}_2+a_{24}\tilde{u} +a_{25} v]+O(\hat{\lambda}^k),  \\
 \gamma^{-k }y' &= \mu + c_1 \tilde{x}_1 +c_2 \tilde{x}_2 +d(\tilde{y} - y^-)^2+ c_4 \tilde{u} + c_5 v + O(\hat{\gamma}^{-k}), \\
u' &= O(\hat{\lambda}^k) ,\qquad v = O(\hat{\gamma}^{-k}) .
\end{aligned}
\end{equation}
Substituting \eqref{64} into the third equation of \eqref{eq:maps:T1_crossnew}, we obtain
\begin{align}\label{640}
\bar{y}
&=\mu+c_1\big[\lambda^k( x_1^+\cos k\omega -x_2^+\sin k\omega )+\lambda^k\cos k\omega(a_{11}\tilde{x}_1 +a_{12}\tilde{x}_2+b_1(\tilde{y}-y^-)+a_{14}\tilde{u}) \nonumber\\
&\quad -\lambda^k\sin k\omega(a_{21}\tilde{x}_1 +a_{22}\tilde{x}_2+a_{24}\tilde{u})\big] + c_2\big[\lambda^k(x_1^+\sin k\omega  +x_2^+\cos k\omega )\nonumber\\
&\quad +\lambda^k\sin k\omega(a_{11}\tilde{x}_1 +a_{12}\tilde{x}_2+b_1(\tilde{y}-y^-)+a_{14}\tilde{u})  +\lambda^k\cos k\omega(a_{21}\tilde{x}_1 +a_{22}\tilde{x}_2+a_{24}\tilde{u})\big]\nonumber\\
&\quad +d\big[\mu\gamma^k+c_1\tilde{x}_1\gamma^k +c_2\tilde{x}_2\gamma^k+d\gamma^k(\tilde{y}-y^-)^2+c_4\gamma^k\tilde{u}+O(\gamma^k\hat{\gamma}^{-k})-y^-\big]^2+O(\hat{\lambda}^k)+c_5\bar{v}.
\end{align}

By definition, $\mu^{new}$ equals the right hand side of the above equation after substituting $(\tilde{x}_1,\tilde{x}_2, \tilde{y},\tilde{u},\bar{v})=(0,0,y^{-,j}_{k},0,0)$ (see \eqref{eq:65}) into it, which is given by
\begin{align*}
\mu^{new}
=\mu&+c_1\lambda^k\big[(x_1^+\cos k\omega-x_2^+\sin k\omega)+b_1\cos k\omega(y^{-,j}_{k}-y^-)\big]\\
&+c_2\lambda^k\big[(x_1^+\sin k\omega+x_2^+\cos k\omega)+b_1\sin k\omega(y^{-,j}_{k}-y^-)\big]\\
&+d\big[\mu\gamma^k+d\gamma^k(y^{-,j}_{k}-y^-)^2+O(\gamma^k\hat{\gamma}^{-k})-y^-\big]^2+O(\hat{\lambda}^k).
\end{align*}
The sought derivative is
\begin{equation}\label{eq:650}
\dfrac{\partial \mu^{new}}{\partial \mu}=1+2 d(\mu\gamma^k+d\gamma^k(y^{-,j}_{k}-y^-)^2+O(\gamma^k\hat{\gamma}^{-k})-y^-)(\gamma^k+O(\gamma^k\hat{\gamma}^{-k})).
\end{equation}
Before checking this derivative we notice that, since ${\partial\bar{y}}/{\partial\tilde{y}}=0$ at the tangency point, one can find from \eqref{640} that
\begin{equation*}\label{eq:66}
 \begin{aligned}
 0=&\left.\dfrac{\partial\bar{y}}{\partial\tilde{y}}\right|_{(0,0,y^{-,j}_{k},0,0)}\\
=&b_1\alpha^*\lambda^k+4b^2_3\gamma^{k}(y^{-,j}_{k}-y^-)(\mu\gamma^k+d\gamma^k(y^{-,j}_{k}-y^-)^2+O(\gamma^k\hat{\gamma}^{-k})-y^-)+O(\hat \lambda^k),
\end{aligned}
\end{equation*}
where the last equality used \eqref{eq:maps:alphabeta}. Since $\alpha^*\neq 0,b_1\neq 0$ by assumption, combining the above equation with \eqref{eq:650} and using \eqref{eq:65}, yields
\[\dfrac{\partial \mu^{new}}{\partial \mu}\Big|_{\mu=\mu_k}=1-\dfrac{b_1\alpha^*\lambda^k(\gamma^k+O(\gamma^k\hat{\gamma}^{-k}))}{2d\gamma^k \hat X^j_k+O(\hat\lambda^k)}\not=0,\]
where $\hat X^j_k=O(\lambda^\frac{k}{2}\gamma^{-k})\neq 0$ by \eqref{eq:coornewHT}.

\noindent {\bf (3) Genericity conditions.}  We now prove that the secondary homoclinic tangencies satisfy conditions C1--C5 introduced in Section \ref{sec:gc}. Condition C2 is shown in the first part of the proof. 
Condition C3 is also fulfilled since $T_1^{new}=T_1\circ T_0^k\circ T_1$, where $T_1$ satisfies C3 and the strong-stable and strong-unstable bundles near $O$ are  invariant under the action of 
$\D T_0^k$.

Recall that condition C5 states that the position vector $(x^+_1,x^+_2)$ of the tangency and the derivative vector $\beta$ are not parallel, see the end of Section \ref{sec:gc}. Here the tangency points are given by $T_1^{new}(M^{-,j}_k)$, where $M^{-,j}_k\in W^u_{loc}(O)$ are given by \eqref{eq:65}. So, the position vector is the $x$-coordinates $(\bar x_1,\bar x_2)$ of $T_1^{new}(M^{-,j}_k)$, and the derivative vector is $\beta=(\partial \bar x_1/\partial \tilde y, \partial \bar x_2/\partial \tilde y)$. Let us first find $\beta$. Substituting \eqref{64} into the first two equations of \eqref{eq:maps:T1_crossnew}, we obtain
\begin{equation}\label{eq:x1x2}
\begin{aligned}
\bar{x}_1
&=x_1^++a_{11}\big[\lambda^k(x_1^+\cos k\omega-x_2^+\sin k\omega)+\lambda^k\cos k\omega(a_{11}\tilde{x}_1+a_{12}\tilde{x}_2+b_1(\tilde{y}-y^-))\\
&\quad -\lambda^k\sin k\omega(a_{21}\tilde{x}_1+a_{22}\tilde{x}_2+a_{24}\tilde{u}+a_{25} v)+O(\hat{\lambda}^k)\big]+a_{12}\big[\lambda^k(x_1^+\sin k\omega+x_2^+\cos k\omega)\\
&\quad +\lambda^k\sin k\omega(a_{11}\tilde{x}_1+a_{12}\tilde{x}_2+b_1(\tilde{y}-y^-))
+\lambda^k\cos k\omega(a_{21}\tilde{x}_1+a_{22}\tilde{x}_2+a_{24}\tilde{u}+a_{25} v)+O(\hat{\lambda}^k)\big]\\
&\quad +b_1\gamma^k\big[\mu+c_1\tilde{x}_1+c_2\tilde{x}_2+d(\tilde{y}-y^-)^2+c_4\tilde{u}+c_5 v +O(\hat{\gamma}^{-k})-\gamma^{-k}y^-\big]+O(\hat{\lambda}^k)+a_{15}\bar{v},\\[5pt]
\bar{x}_2
&=x_2^++a_{21}\big[\lambda^k(x_1^+\cos k\omega-x_2^+\sin k\omega)+\lambda^k\cos k\omega(a_{11}\tilde{x}_1+a_{12}\tilde{x}_2+b_1(\tilde{y}-y^-))\\
&\quad -\lambda^k\sin k\omega(a_{21}\tilde{x}_1+a_{22}\tilde{x}_2+a_{24}\tilde{u}+a_{25} v)+O(\hat{\lambda}^k)\big]+a_{22}\big[\lambda^k(x_1^+\sin k\omega+x_2^+\cos k\omega)\\
&\quad +\lambda^k\sin k\omega(a_{11}\tilde{x}_1+a_{12}\tilde{x}_2+b_1(\tilde{y}-y^-))
+\lambda^k\cos k\omega(a_{21}\tilde{x}_1+a_{22}\tilde{x}_2+a_{24}\tilde{u}+a_{25} v)+O(\hat{\lambda}^k)\big]\\
&\quad +O(\hat{\lambda}^k)+a_{15}\bar{v}.
\end{aligned}
\end{equation}
Thus,
\begin{align}\label{eq:661}
\dfrac{\partial \bar{x}_1}{\partial \tilde{y}}
&= a_{11}b_1\lambda^k\cos k\omega +a_{12}b_1\lambda^k\sin k\omega +2b_1d\gamma^k(\tilde{y}-y^-),\nonumber\\
\dfrac{\partial \bar{x}_2}{\partial \tilde{y}}
&= a_{21}b_1\lambda^k\cos k\omega +a_{22}b_1\lambda^k\sin k\omega.
\end{align}
Substituting $\tilde y= y^{-,j}_k$ into \eqref{eq:661}, with using \eqref{eq:coornewHT} and \eqref{eq:65}, gives
\begin{equation}\label{eq:662}
\left(\dfrac{\partial \bar{x}_1}{\partial \tilde{y}},\dfrac{\partial \bar{x}_2}{\partial \tilde{y}}\right)\Big|_{M_k^{-,j}}=\left(\pm\sqrt{\dfrac{-(\alpha^*x_1^++\beta^*x_2^+)}{d}} \dfrac{\alpha^*b_1^2\lambda^{\frac{k}{2}}}{2(\alpha^*x_1^++\beta^*x_2^+)}+o(\lambda^{\frac{k}{2}})+o(\gamma^k\hat{\gamma}^{-k}),O(\lambda^k)\right).
\end{equation}

To obtain the $x$-coordinates $(\bar{x}_1,\bar{x}_2)$ of $T_1^{new}(M^{-,j}_k)$, we substitute in the right hand sides of \eqref{eq:x1x2} the $\tilde x,\tilde y,\tilde u$ coordinates of $M^{-,j}_k$ (see \eqref{eq:65}) with setting $\bar v=0$, and get
\begin{equation}\label{eq:663}
(\bar{x}_1,\bar{x}_2)=(x_1^++b_1y^-+o(1)_{k\to\infty},\;x_2^++o(1)_{k\to \infty}).
\end{equation}  
Since $x^+_2\neq 0$ by \eqref{eq:x+2},  the two vectors \eqref{eq:662} and \eqref{eq:663} are not parallel for all sufficiently large $k$, and hence condition C5 is satisfied. Equation \eqref{eq:663} immediately implies condition C4.

Finally, let us check condition C1. Suppose for the contradiction that it
 does not hold. Then, there are two linearly independent vectors in the intersection of the tangent spaces at $M^{-,j}_k$ of $W^u_{loc}(O)$ and $(T_1^{new})^{-1}(W^s_{loc}(O))$. Since the weak unstable dimension is 1, some linear combination of these vectors, which is also a common tangent vector, must lie in the strong-unstable subspace. This contradicts  the transversality between the strong-unstable foliation and $(T_1^{new})^{-1}(W^{sE}_{loc}(O))$ given by condition C3.
The proof of Lemma \ref{lem:tangency2} is now complete.

\section{Creating robust heterodimensional dynamics: proof of Theorem \ref{thm:sf} for the saddle-focus case}\label{sec:robusthdc}

Theorem \ref{thm:sfsimple} implies Theorem \ref{thm:sf} (for the saddle-focus case) by virtue of the following stabilization result:
\begin{thms}[{\citep[Theorem 7]{LT21}}]\label{thm:lt21}
Let a $C^r$ system have a non-degenerate heterodimensional cycle $H$ involving two hyperbolic periodic points $O_1$ and $O_2$ with $ind(O_1)+1=ind(O_2)$ and let the center-stable multiplier of $O_1$ be $\lambda_{1,1}=\lambda_{1,2}^*=\lambda e^{i\omega}$ with $\omega \in (0,\pi)$ and the center-unstable multiplier of $O_2$ be real $\nu$. If the quantities $\ln\lambda/\ln|\nu|$, $\omega/{\pi}$ and $1$ are rationally independent, then, in any neighborhood of $H$, the system  has  $C^1$-robust heterodimensional dynamics involving a cs-blender and a cu-blender. Moreover,  for any one-parameter family unfolding this cycle, there exist intervals of parameter values converging to $0$ for which
the point $O_1$ is homoclinically related to the cs-blender.
\end{thms}

By Theorem \ref{thm:sfsimple}, there is a heterodimensional cycle involving $O$ and $Q$ for $f_{\mu_k,\omega_k}$ for every sufficiently large $k$, and, moreover, this cycle unfolds as $\omega$ changes from $\omega_k$. Thus, Theorem \ref{thm:sf} is obtained by applying the above result to the family $f_\eps:=f_{\mu,\omega}$ near $(\mu,\omega)=(\mu_k,\omega_k)$, with $O_1=O$ and $O_2=Q$, provided that the multiplier $\nu$ of $Q$ is real, the required rational independence condition is satisfied, and the heterodimensional cycle found in Theorem \ref{thm:sfsimple} is non-degenerate in the sense of \citep{LT21}.  Note that, by Lemma \ref{lem:index2} and Proposition \ref{prop:nontr}, the parameter values $(\mu_k(t),\omega_k(t))$ are functions of $t$.  In what follows, we show that the above conditions are indeed satisfied at least for some values of $t$.

\subsection{Adjustment of the central multipliers}
Recall that we denote the center-unstable multipliers of $Q$ by $\nu$ and $\hat \nu$. When $t$ is bounded away from zero, these multipliers are real and $|\nu|<|\hat\nu|$, see \eqref{eigen:3}.
The multipliers are functions of $t$, and  also depend on $k$ since $Q$ does.

Note also that the absolute value $\lambda$ of the complex leading multiplier of $O$ depends on
the parameters $\mu$ and $\omega$, so it is also a function of $t$ and $k$. By \eqref{eq:t_deri} and \eqref{omtd},
\begin{equation}\label{lamtd}
\frac{d}{dt} \lambda(t) = \dfrac{\partial \lambda}{\partial \mu}\dfrac{\partial \mu_k}{\partial t}
+ \dfrac{\partial \lambda}{\partial \omega}\dfrac{\partial \omega_k}{\partial t}
= O(\lambda^k).
\end{equation}

\begin{lem}\label{lem:irra}
For all sufficiently large $k$, arbitrarily close to $t=t_+$ and $t=t_-$ given by \eqref{tpm}, there exist values of $t$ such that the quantities $\ln \lambda(t) /\ln |\nu(t)|$, $\omega_k(t)/\pi$ and 1 are rationally independent.
\end{lem}

\begin{proof} Denote $A(t)=\ln \lambda(t) /\ln |\nu(t)|$, $B(t)=\omega_k(t)/\pi$. By \eqref{lamtd} and \eqref{eigen:3}, the derivative $A'(t)$ is bounded away from zero for all
large $k$. Therefore, the line $(A(t),B(t))$ is a smooth curve in the $(A,B)$ plane. Since different straight lines
$$i_1 A + i_2 B = i_3$$
with integers $i_{1,2,3}$ intersect at non-zero angles, a smooth curve can lie in the union of these lines over all $i \in \mathbb{Z}^3\setminus\{0\}$ only if it entirely lies in one of them. Thus, if $A(t)$, $B(t)$ and $1$
are rationally dependent for all $t$ from some interval, then there exist integers $(i_1,i_2,i_3)$, not all zero, such that 
\begin{equation}\label{irra:1}
i_1A(t)+i_2 B(t)+i_3=0
\end{equation}
for all $t$ from the interval under consideration -- with the same $i_{1,2,3}$. 

So, let \eqref{irra:1} hold for an interval approaching $t=t_-$ or $t=t_+$. Since $\nu(t)\to \pm 1$ as $t\to t_\pm$ by Lemma \ref{lem:index2}, we have that $A(t_\pm)=\infty$. Thus, $i_1=0$, but
this would mean that $B(t)$ is constant, i.e., $\frac{d}{dt} \omega_k=0$, which contradicts to \eqref{omtd}. 
\end{proof}

\subsection{Non-degeneracy conditions for heterodimensional cycles}\label{sec:HDCcon}
The heterodimensional cycle involving $O$ and $Q$ is non-degenerate if the non-degeneracy conditions
C1--C4 in  \citep[Sections 2.2 and 2.3]{LT21} are satisfied, as described below (where we relabel them as NC1--NC4). We assume that $t$ is bounded away from zero, so the multipliers $\nu$ and $\hat\nu$ of $Q$ are real.

Recall, as we discussed in Section \ref{sec:gc}, that there exist a strong-stable foliation $\mathcal{F}^{ss}_O:=\mathcal{F}^{ss}$ on $W^s_{loc}(O)$ and a local extended unstable manifold $W^{uE}_{loc}(O)$, whose forward iterations give the global extended unstable manifold $W^{uE}(O)$. Similarly, there exists a strong-unstable foliation 
$\mathcal{F}^{uu}_Q$ on $W^u_{loc}(Q)$, including, as a leaf, the manifold $W^{uu}_{loc}(Q)$ defined in Section \ref{sec:WuuQ}. There is also a local extended stable manifold $W_{loc}^{sE}(Q)$ whose backward iterations give the global extended stable manifold $W^{sE}(Q)$, as discussed in Section \ref{sec:Wse}. The first three non-degeneracy conditions of \citep[Section 2.2]{LT21}, adapted to our setting, are as follows.

\noindent\textbf{NC1.} 
At a point of the (non-transverse) intersection between $W^u(O)$ and $W^s(Q)$, the manifold $W^{uE}(O)$ is transverse to  $W^{s}(Q)$ and $W^{u}(O)$ is transverse to $W^{sE}(Q)$.\\[5pt]
\noindent\textbf{NC2.} 
At a point of the (transverse) intersection between  $W^u(Q)$ and $W^s(O)$, 
the leaf of the foliation 
$\mathcal{F}^{uu}_Q$ through this point is not tangent to
$W^s(O)$ and the leaf of $\mathcal{F}^{ss}_O$ through this point is not tangent to
$W^u(Q)$. \\[5pt]
\noindent\textbf{NC3.} The intersection point in \textbf{NC2}  belongs to neither $W^{ss}(O)$ nor $W^{uu}(Q)$.

Note that $\mathcal{F}^{ss}_O$ is extended to a small neighborhood of $O$ as a continuous invariant foliation with smooth leaves whose tangents lie in the cone $\mathcal{C}^{ss}$. So, $W^s_{loc}(Q)$ is a leaf of this extended foliation and the first part of condition C3 implies the first part of NC1, since $Q$ is close to $M^+$. The second part of NC1 follows from Lemma \ref{lem:gc1}. 

Remark~\ref{lem:gc2} implies that for $t$ satisfying \eqref{tsign} the manifold $W^s(O)$ is transverse to all leaves of $\mathcal{F}^{uu}_Q$ which are close to $W^{uu}_{loc}(Q)$. Hence the first part of NC2 is satisfied by any point on the curve of the intersection of $W^s(O)$ and $W^{u}_{loc}(Q)$. Since the leaves of
$\mathcal{F}^{ss}_O$ lie in the cone $\mathcal{C}^{ss}$ of \eqref{eq:cone_ss}, the second part of NC2 is immediate because the tangent to $W^u_{loc}(Q)$ lies in the cone $\mathcal{C}^{cu}$ of \eqref{eq:cone_cu}, see Lemma \ref{lem:WuQ}.

To verify NC3 for $t$ satisfying \eqref{tsign}, remind that the transversality of $W^{uu}_{loc}(Q)$ to
$W^s(O)$ given by Remark \ref{lem:gc2} implies that $W^{uu}_{loc}(Q)$ intersects the curve of
the intersection of $W^s(O)$ and $W^{u}_{loc}(Q)$ only at one point. Therefore, there are points on this curve which are close to $W^{uu}_{loc}(Q)$ (so NC2 still holds) but do not belong to $W^{uu}_{loc}(Q)$.
Finally, recall that, as shown in the proof of Proposition \ref{prop:tran}, the pieces $W^{s\pm}$ of $W^s(O)$, whose backward iterations create the intersection with
$W^{u}_{loc}(Q)$,  are some backward iterations of a piece of $W^s_{loc}(O)$ near $M^-\not \in W^{ss}(O)$ (see the proof of Lemma \ref{lem:tangency1}). Thus, we also have that
the intersection points of $W^s(O)$ and $W^{u}_{loc}(Q)$ do not belong to $W^{ss}(O)$.\\

The case where the multiplier $\nu$ of $Q$ is real corresponds to the ``saddle-focus cycle'' of \citep[Section 2.3]{LT21}, and the 4th non-degeneracy condition  reads in this case as follows.

\noindent\textbf{NC4.} In the coordinates where the map near $O$ assumes the form \eqref{eq:maps:T0}, the $x$-vector component of the tangent vector $p$ to the curve $W^s_{loc}(O)\cap W^u(Q)$ at some point $P\notin W^{uu}(Q)$ on this curve is not parallel to the vector $(x^+_1,x^+_2)$.

Let us show that condition NC4 is fulfilled. Recall that, in the proof of Proposition \ref{prop:tran}, we create an intersection of $W^u_{loc}(Q)$ with one of  the preimages $T_0^{-k}(W^{s+})$ and $T_0^{-k}(W^{s-})$, where $W^{s\pm}$ are connected pieces of $W^s(O)$ containing the accompanying transverse homoclinic points near $M^-$ given by Lemma \ref{lem:tangency1}. Let us assume that the intersection is with $W^{s+}$, and the reasoning for $W^{s-}$ is the same. As shown in the proof of Lemma \ref{lem:tangency1}, there is a small piece $\hat W^{s+}$ of $W^s_{loc}(O)$ near $M^+$ such that $W^{s+}=T_1^{-1}\circ T_0^{-i}\circ T_1^{-1}(\hat W^{s+})$ for some integer $i>0$ which can be taken as large as we want. Since $T_0^{-k}(W^{s+})\cap W_{loc}^u(Q)$ does not lie in $W^{uu}(Q)$ by Remark~\ref{lem:gc2}, we can pick 
a point $P'\in T_0^{-k}(W^{s+})\cap (W_{loc}^u(Q)\setminus W^{uu}(Q))$. We claim that the image $P:= T_1\circ T^i_0 \circ T_1\circ T^k_0(P')$ is the sought point  in condition NC4, which we prove as follows.

We fix a large $i$ and denote $F_k:=T_1\circ T^i_0 \circ T_1\circ T^k_0$. The differential $\D F_k$ is the product of two linear maps $\D T_0^k$ and $\D G_i$ with $ G_i:=T_1\circ T_0^i\circ T_1$. The map $G_i$ acts from a small neighborhood of a point in $W^u_{loc}$ to a neighborhood of a  point in $W^s_{loc}$, where these two points belong to the same transverse homoclinic orbit. By the transversality and by  condition C3, the image of the $x$-subspace by $\D G_i$ has non-zero projection to the $x$-subspace.
The differential $\D T^k_0$ gives a rotation in the $x$-component to the angle $k\omega$ by \eqref{eq:maps:T0}, and  $(k\omega \bmod 2\pi)$ takes at least three different values since $\omega\in (0,\pi)$. Sine $\D G_i$ is $O(\delta)$-close to a constant linear map, it follows that, for any vector $p'$ with a non-zero $x$-component, the angle of the $x$-component of its image $\D F_k(p')$ gets in the disjoint $O(\delta)$-neighborhoods of at least three different non-zero values. This implies that  there exist infinitely many $k$ such that the $x$-component of $\D F_k(p')$  is not parallel to $(x^+_1,x^+_2)$.   Now, to conclude the proof, it suffices to show that the tangent vector to the curve $F_k^{-1}(W^s_{loc}(O)\cap W^u(Q))\subset T_0^{-k}(W^{s+})$ at $P'$ has a non-zero $x$-component. But this is immediate from \eqref{eq:tran:3} and \eqref{eq:maps:coor3}.

\section{Bi-focus case}\label{sec:df}
In this section, we prove Theorem \ref{thm:sf} for the bi-focus case. Let $O$ be a bi-focus with central multipliers $\lambda_1=\lambda^*_2=\lambda e^{i\omega_1}$ $(0<\omega_1<\pi)$ and $\gamma_1=\gamma^*_2=\gamma e^{i\omega_2}$ $(0<\omega_2<\pi)$. We assume that $\lambda\gamma>1$ and conditions C1--C5 are fulfilled. The theorem will be proved through a reduction to the saddle-focus case, which is achieved by introducing suitable coordinates such that the first-return map in this case has the same form as \eqref{eq:maps:Tk}.

\subsection{First-return map and its normal form}\label{sec:dfreturn}
As in the saddle-focus case, there exist $C^r$ coordinates \citep[Lemma 6]{GST08} such that the local map $T_0$ assumes the same form as \eqref{eq:maps:T0} with taking $y=(y_1,y_2)\in \mathbb{R}^2$. More specifically, the $\bar y$-equation in \eqref{eq:maps:T0} is replaced by
\begin{align*}
\bar y_1& = \gamma  y_1\cos\omega - \gamma y_2\sin\omega  + p_{3,1}(x,y,u,v),\\ 
\bar y_2& = \gamma y_1  \sin\omega  + \gamma y_2 \cos\omega +  p_{3,2}(x,y,u,v),
 \end{align*}
where we did not write the dependence on parameters explicitly and the function $p_3:=(p_{3,1},p_{3,2})$ satisfies \eqref{eq:maps:T0_nonlinear}. Then, by \citep[Lemma 7]{GST08}, for any point $(x_{1},x_{2},y_1,y_2,u,v)\in \Pi^+$, it holds that $(\tilde x_{1},\tilde x_{2},\tilde y_{1},\tilde y_{2},\tilde u,\tilde v)= T_0^k(x_{1},x_{2},y_1,y_2,u,v)\in \Pi^-$ if and only if the points satisfy the relation
\begin{equation}\label{eq:df:T0k}
\begin{aligned}
\tilde x_{1} &= \lambda^k x_1\cos k\omega_1  - \lambda^k x_2\sin k\omega_1 + O(\hat{\lambda}^k),\\ 
 \tilde x_{2}& = \lambda^k x_1  \sin k\omega_1  + \lambda^k x_2 \cos k\omega_1  + O(\hat{\lambda}^k),\\ 
y_1 &= \gamma^{-k} \tilde y_{1}\cos k\omega_2 + \gamma^{-k}\tilde y_{2}\sin k\omega_2+ O(\hat \gamma^{-k}),\\
y_2 &= -\gamma^{-k} \tilde y_{1}\sin k\omega_2 + \gamma^{-k}\tilde y_{2}\cos k\omega_2+ O(\hat \gamma^{-k}),\\
 \tilde u  & = O(\hat{\lambda}^k) ,\qquad
 v  = O(\hat \gamma^{-k}),
\end{aligned}
\end{equation}
where $\hat\lambda \in (|\lambda_3|,\lambda)$ and $\hat\gamma\in (\gamma,|\gamma_3|)$ are chosen to satisfy \eqref{eq:lambdagamma} and
\begin{equation}\label{eq:lambdagamma2}
\hat \gamma^{-1}>\lambda\gamma^{-1}.
\end{equation}
Here the $O(\cdot)$ terms are functions of $x_1,x_2,\tilde y_1,\tilde y_2,u,\tilde v$, and their first derivatives with respect to coordinates have the same forms.

Using conditions C1--C3, we apply \citep[Corollary 1]{GST08}. It says that the global map $T_1$ from $\Pi^-$ to $\Pi^+$ satisfy that for any point $(\tilde x_{1 },\tilde x_{2 },\tilde y_{1 },\tilde y_{2 },\tilde u ,\tilde v )\in \Pi^-$, we have $( x_{1}, x_{2},y_1, y_2, u, v) = T_1(\tilde x_{1 },\tilde x_{2 },\tilde y_{1 },\tilde y_{2 },\tilde u ,\tilde v )$ if and only if
\begin{equation}\label{eq:df:T1_cross}
\begin{aligned}
{x}_1 - x_1^+ &=a_{11} \tilde x_{1 } +a_{12} \tilde x_{2 } +b_{11}(\tilde y_{1 } - y_1^-)+b_{12}y_2 + a_{14} \tilde u  + a_{15}  v + \dots, \\
{x}_2 - x_2^+ &=a_{21} \tilde x_{1 } +a_{22} \tilde x_{2 }  +b_{22} y_2+ a_{24} \tilde u  + a_{25}  v + \dots, \\
 y_1 &=\mu + c_1 \tilde x_{1 } +c_2 \tilde x_{2 } + c_4\tilde  u  + c_5  v +d(\tilde y_{1 } - y_1^-)^2
 + \dots, \\
{u} - u^+ &=a_{41} \tilde x_{1 } +a_{42} \tilde x_{2 } +b_{41}(\tilde y_{1} - y_1^-)+b_{42} y_2 + a_{44} \tilde u  + a_{45}  v + \dots ,\\
\tilde y_{2} - y^-_2
&=a_{51} \tilde x_{1 } +a_{52} \tilde x_{2 } +b_{51}(\tilde y_1 - y_1^-)+ 
b_{52}y_2+a_{54} \tilde u+a_{55} v
+ \dots ,\\
\tilde v - v^-&=a_{61} \tilde x_{1 } +a_{62} \tilde x_{2 } +b_{61}(\tilde y_1  - y_1^-)+ 
b_{62} y_2
+a_{64}\tilde u+a_{65} v+ \dots ,
\end{aligned}
\end{equation}
where $b_{11}\neq 0,d\neq 0, c^2_{1}+c^2_{2}\neq 0, b_{52}\neq 0, \det a_{65}\neq 0$, and the dots denote the second and higher order terms, except that 
the dots in the $y_1$-equation stand for $O(|\tilde  y - y^-|^3 + (\tilde x, \tilde u,y_2,v)^2 + 
|\tilde  y - y^-| \cdot \|\tilde x, \tilde u,y_2,v \|)$. Condition C3 implies 
\begin{equation}\label{eq:gc3_df}
(x^+_1)^2 + (x^+_2)^2\neq 0
\quad\mbox{and}\quad
(y^-_1)^2 + (y^-_2)^2\neq 0.
\end{equation}
Here the vector $\beta$ in condition C5  is $(b_{11},0)$, so condition C5 means
\begin{equation}\label{eq:gc6_df}
x^+_2\neq 0.
\end{equation}

The first return map $T_k:= T_1\circ T^k_0$ is such that $(\bar x,\bar y,\bar u,\bar v)=T_k(x,y,u,v)$ if and only if there is $(\tilde x,\tilde y,\tilde u,\tilde v)\in \Pi^-$ such that 
\begin{equation*}
(x,y,u,v)
\xmapsto{T^k_0}(\tilde x,\tilde y,\tilde u,\tilde v)
\xmapsto{T_1}(\bar x,\bar y,\bar u,\bar v).
\end{equation*}
As in the saddle-focus case, after the Shilnikov coordinate transformation 
\begin{equation}\label{eq:df:coor1}
\begin{aligned}
&X_1 = {x_1} - x_1^+,\quad
X_2 = {x}_2 - x_2^+,\quad
Y_1 = \tilde y_1- y^-_1,\quad
Y_2 = \tilde y_{2} - y_2^-,\\
&U = u - u^+,\quad
V=\tilde  v-v^-.
\end{aligned},
\end{equation}
the map $T_k:(X,Y,U,V)\mapsto (\bar X,\tilde Y,\bar U,\tilde V)$ is defined on
\begin{equation*}
(X,Y,U,V)\in[-\delta,\delta]^2\times[-\delta,\delta]^2\times[-\delta,\delta]^{d^s-2}\times[-\delta,\delta]^{d^u-2}=:\hat\Pi.
\end{equation*}

Let us set $(x,y,u,v)=(\bar x,\bar y,\bar u,\bar v)$ in \eqref{eq:df:T1_cross} and combine it with \eqref{eq:df:T0k}. 
As in Section \ref{sec: Shilnikov}, we find that $(\bar X_1,\bar X_2,\bar Y_1,\bar Y_2,\bar U,\bar V)=T_k(X_1,X_2,Y_1,Y_2,U,V)$ if and only if
\begin{equation}\label{eq:df:Tk0}
\begin{aligned}
&\bar X_1 = b_{12}\mu  + \lambda^k \hat\alpha_1(X_1+ x_1^+) +\lambda^k \hat \beta_1(X_2+x_2^+) + b_{11} Y_1 
 + \hat h_1(X_1,X_2,Y_1,U,\bar Y_2,\bar V),
\\
&\bar X_2 = b_{22}\mu +\lambda^k \hat\alpha_2(X_1+ x_1^+) +\lambda^k \hat \beta_2(X_2+x_2^+) 
+ \hat h_2(X_1,X_2,Y_1,U,\bar Y_2,\bar V),
\\
&\bar Y_1\cos k\omega_2 +\bar Y_2\sin k\omega_2 = \gamma^k \mu-y^* + \lambda^k\gamma^k \alpha^*(X_1+ x_1^+) +\lambda^k\gamma^k  \beta^*(X_2+x_2^+) \\
&\qquad\qquad\qquad\qquad\qquad\qquad\qquad\qquad\qquad\qquad\quad
+ d\gamma^k Y_1^2   + \gamma^k\hat h_3(X_1,X_2,Y_1,U,\bar Y_2,\bar V),
\\
&\bar U =b_{42}\mu  + \lambda^k \hat\alpha_4(X_1+ x_1^+) +\lambda^k \hat \beta_4(X_2+x_2^+) + b_{41} Y_1
+ \hat h_4(X_1,X_2,Y_1,U,\bar Y_2,\bar V),
\\
&Y_2=b_{52}\mu +
\lambda^k \hat\alpha_5(X_1+ x_1^+) +\lambda^k \hat \beta_5(X_2+x_2^+) + b_{51} Y_1
+ \hat h_5(X_1,X_2,Y_1,U,\bar Y_2,\bar V),
\\
 &V=b_{62}\mu  +
\lambda^k \hat\alpha_6(X_1+ x_1^+) +\lambda^k \hat \beta_6(X_2+x_2^+) + b_{61} Y_1
+ \hat h_6(X_1,X_2,Y_1,U,\bar Y_2,\bar V),
\end{aligned}
\end{equation}
where the  terms $b_{i2}\mu$ appear in the following way: after replacing $y_2$ by $\bar y_2$  in \eqref{eq:df:T1_cross}, the terms $b_{i2}  y_2$ become
\begin{align*}
b_{i2} \bar y_2 &= b_{i2}\gamma^{-k} ((\bar Y_1+y^-_1)\cos k\omega_2 -(\bar Y_2+y^-_2) \sin k\omega_2) + O(\hat\gamma^{-k}),
\end{align*}
and the coefficients are given by
\begin{equation}\label{eq:df:alphabeta}
\begin{aligned}
\alpha^* &= c_{1}\cos k\omega_1+c_{2}\sin k\omega_1,
\qquad
 &\beta^* &= -c_{1}\sin k\omega_1+c_{2}\cos k\omega_1,
\\
\hat \alpha_i &= a_{i1}\cos k\omega_1+a_{i2}\sin k\omega_1,
\qquad
&\hat \beta_i &= -a_{i1}\sin k\omega_1+a_{i2}\cos k\omega_1,
\quad (i=1,2,4,5,6)
\\
y^* &= y^-_1\cos k\omega_2 +y^-_2\sin k\omega_2.&&
\end{aligned}.
\end{equation}

Comparing with   $\hat{h}_3$  in \eqref{ybarexp}, the function  $\hat{h}_3$ here has  extra terms involving $\bar Y_2$. They  come from the terms involving $y_2$ in the dots of the $y_1$-equation in \eqref{eq:df:T1_cross}. We see that the contribution made by products of $y_2$ with other variables is of order $O(\hat\lambda^k)$, and that made by powers of $y_2$ is of order $\gamma^{-2k}$ -- smaller than the term $d_3 v=O(\hat\gamma^{-k})$ (by the assumption $\lambda\gamma>1$ and the last condition in \eqref{eq:df:T1_cross}). It follows that the estimates for the derivatives of $\hat h_3$  are the same as in the saddle-focus case, except for $\partial \hat h_3/\partial Y_2$ which is only present in the bi-focus case. Thus, with the assumption $|\lambda\gamma|>1$, \eqref{eq:lambdagamma} and \eqref{eq:lambdagamma2}, we obtain
\begin{equation*}\label{eq:df:h}
\begin{array}{l}
\hat h_3=O(\hat\lambda^k)+O(\lambda^k)Y_1+O(Y^3_1),\qquad
\dfrac{\partial \hat h_3}{\partial (X,U)} 
=O(\hat\lambda^k)+O(\lambda^k)Y^2_1,\\
\dfrac{\partial \hat h_3}{\partial Y_1} 
= O(\lambda^k)+O(\gamma^k\hat\gamma^{-k})Y_1+O(Y^2_1),
\qquad
\dfrac{\partial \hat h_3}{\partial (\bar Y_2,\bar V)}=O(\gamma^{-k})Y_1+O(\hat\gamma^{-k}),\\
\dfrac{\partial \hat h_3}{\partial \mu}
=O(k\gamma^k\hat\gamma^{-k})+o(Y^2_1), \qquad 
\dfrac{\partial \hat h_3}{\partial \omega}=k\gamma^k\hat\gamma^{-k} 
O(|\mu| + \lambda^k+ Y^2_1 ) + o(Y^2_1).
\end{array}
\end{equation*}
and, for $i=1,2,4,5,6$, 
\begin{equation*}\label{eq:df:h2}
\begin{array}{l}
\hat h_i 
=O(\hat \lambda^k)+O(\lambda^k)Y_1+O(Y^2_1),\qquad
\dfrac{\partial \hat h_i}{\partial(X,U)} 
=O(\hat\lambda^k)+O(\lambda^k)Y_1,\\
\dfrac{\partial \hat h_i}{\partial Y_1}
=O(\lambda^k)+O(Y_1), \qquad
\dfrac{\partial \hat h_i}{\partial (\bar Y_2,\bar V)}
=O(\gamma^{-k}),\\
\dfrac{\partial \hat h_i}{\partial \mu}
=O(k\gamma^k\hat\gamma^{-k})+O(Y^2_1), \qquad
\dfrac{\partial \hat h_i}{\partial \omega}
=O(k\lambda^k\gamma^k\hat\gamma^{-k})+O(Y^2_1).
\end{array}
\end{equation*}

We will further consider $k$ such that 
\begin{equation}\label{eq:k}
\cos k\omega_2
\quad\mbox{and}\quad
\cos k\omega_2+ b_{51}\sin k\omega_2\quad\mbox{are bounded away from zero.}
\end{equation}
This can be always achieved since $\omega_2\in(0,\pi)$ implies that $(k\omega_2 \bmod{2\pi})$ takes at least three different values.

We start the reduction to the (2,1) saddle-focus case by introducing the new variable
\begin{equation}\label{eq:df:coor2}
Y=Y_1\cos k\omega_2 + Y_2\sin k\omega_2.
\end{equation}
We find $ Y_2$ from the fifth equation in \eqref{eq:df:Tk0} as a function of $X_1,X_2,Y,U,\bar Y_2,\bar V$:
\begin{equation*}\label{eq:df:Y_2}
 Y_2=\dfrac{1}{1+b_{51}\tan k\omega_2}
\left(
\lambda^k \hat\alpha_5(X_1+ x_1^+) +\lambda^k \hat \beta_5(X_2+x_2^+) + b_{51} Y\cos^{-1} k\omega_2
\right)
+ \tilde h_5.
\end{equation*}
Substituting this with \eqref{eq:df:coor2} into the remaining equations, yields
\begin{equation}\label{eq:df:Tk2}
\begin{aligned}
\bar X_1 &=b_{12}\mu  + \lambda^k \tilde\alpha_1(X_1+ x_1^+) +\lambda^k \tilde \beta_1(X_2+x_2^+) + \tilde b_{11} Y
 + \tilde h_1(X_1,X_2,Y,U,\bar Y_2,\bar V),
\\
\bar X_2 &= b_{22}\mu  +\lambda^k \tilde\alpha_2(X_1+ x_1^+) +\lambda^k \tilde \beta_2(X_2+x_2^+) 
+ \tilde h_2(X_1,X_2,Y,U,\bar Y_2,\bar V),
\\
\bar Y &= \gamma^k \mu-y^*+ \lambda^k\gamma^k \alpha^*(X_1+ x_1^+) +\lambda^k\gamma^k\beta^*(X_2+ x_2^+)  
+ \tilde d\gamma^k Y^2  
 + \gamma^k\tilde h_3(X_1,X_2,Y,U,\bar Y_2,\bar V),
\\
\bar U &=b_{42}\mu + \lambda^k \tilde\alpha_4(X_1+ x_1^+) +\lambda^k \tilde \beta_4(X_2+x_2^+) +  \tilde b_{41} Y
+ \tilde h_4(X_1,X_2,Y,U,\bar Y_2,\bar V),
\\
 Y_2&=b_{52}\mu  +
\lambda^k \tilde\alpha_5(X_1+ x_1^+) +\lambda^k \tilde \beta_5(X_2+x_2^+) +  \tilde b_{51} Y
+ \tilde h_5(X_1,X_2,Y,U,\bar Y_2,\bar V),
\\
V&=b_{62}\mu  +
\lambda^k \tilde\alpha_6(X_1+ x_1^+) +\lambda^k \tilde \beta_6(X_2+x_2^+) +  \tilde b_{61} Y
+ \tilde h_6(X_1,X_2,Y,U,\bar Y_2,\bar V),
\end{aligned}
\end{equation}
where 
\begin{equation}\label{eq:df:alphabeta2}
\begin{array}{l}
\tilde b_{51}=\dfrac{b_{51}}{\cos k\omega_2+b_{51}\sin k\omega_2},\quad

\tilde \alpha_5 = \dfrac{\hat\alpha_5\cos k\omega_2}{\cos k\omega_2+b_{51}\sin k\omega_2},\quad
\tilde \beta_5 = \dfrac{\hat\beta_5\cos k\omega_2}{\cos k\omega_2+b_{51}\sin k\omega_2},\\[10pt]

\tilde b_{i1}=\dfrac{b_{i1}(1-\tilde b_{51}\sin k\omega_2)}{\cos k\omega_2}=\dfrac{b_{i1}}{\cos k\omega_2+b_{51}\sin k\omega_2},\quad

\tilde d=\dfrac{d}{(\cos k\omega_2+b_{51}\sin k\omega_2)^2},\\[10pt]
\tilde \alpha_j 
= \hat \alpha_j-\dfrac{b_{j1}\hat\alpha_5\sin k\omega_2}{\cos k\omega_2+b_{51}\sin k\omega_2},\quad
\tilde \beta_j 
=\hat \beta_j-\dfrac{b_{j1}\hat\beta_5\sin k\omega_2}{\cos k\omega_2+b_{51}\sin k\omega_2},
\end{array}
\end{equation}
for $i=1,4,6$ and $j=1,2,4,6$, and functions $\tilde h$ satisfy the same estimates as $\hat h$, with replacing $Y_1$ by $Y$.

One readily finds that, after denoting $V^{new}=(Y_2,V)$  in formula \eqref{eq:df:Tk2}, the above formula for $T_k$ assumes the form
\begin{equation*}
\begin{aligned}
\bar X_1 &=b_{12}\mu+ \lambda^k \tilde\alpha_1(X_1+ x_1^+) +\lambda^k \tilde \beta_1(X_2+x_2^+) + b_1 Y + \tilde h_1(X,Y,U,\bar V),\\
\bar X_2 &= b_{22}\mu+\lambda^k \tilde\alpha_2(X_1+ x_1^+) +\lambda^k \tilde \beta_2(X_2+x_2^+) + \tilde h_2(X,Y,U,\bar V),\\
\bar Y &= \gamma^k \mu-y^- + \lambda^k\gamma^k \alpha^*(X_1+ x_1^+) +\lambda^k\gamma^k  \beta^*(X_2+x_2^+)
 + \tilde d\gamma^k Y^2 + \gamma^k\tilde h_3(X,Y,U,\bar V),\\
\bar U &=b_{42}\mu+ \lambda^k \tilde\alpha_4(X_1+ x_1^+) +\lambda^k \tilde \beta_4(X_2+x_2^+) + b_4 Y + \tilde h_4(X,Y,U,\bar V),\\
V&=b'_{52}\mu+\lambda^k \tilde\alpha'_5(X_1+ x_1^+) +\lambda^k \tilde \beta'_5(X_2+x_2^+) + b_5 Y + \tilde h'_5(X,Y,U,\bar V),
\end{aligned}
\end{equation*}
where $b'_{52}=(b_{52},b_{62})^T,\tilde \alpha'_5=(\tilde{\alpha}_5,\tilde\alpha_6)^T,\tilde \beta'_5=(\tilde{\beta}_5,\tilde\beta_6)^T,\tilde h'_5=(\tilde{h}_5,\tilde h_6)^T$, and the functions $\tilde h_i$ $(i=1,2,3,4)$ and $\tilde h'_5$ satisfy the same estimates as the functions $\hat h_i$ and, respectively, $\hat h_5$ in \eqref{ybarexp} and \eqref{eq:maps:Tk_1}. Using the rescaling
\begin{equation}\label{eq:rho}
\mu=\lambda^k\rho,
\end{equation} 
we can  apply a series of coordinate transformations as in Subsections \ref{sec:coor1}--\ref{sec:coor3}, where, in particular, the terms $b_{i2}\mu=O(\lambda^k)$ are removed as constants by the transformation \eqref{eq:maps:coor2.1}.
 Finally, formula \eqref{eq:df:Tk2} assumes the form
\begin{equation}\label{eq:df:Tk}
\begin{aligned}
\bar Z &= \lambda^k \alpha_1 Z+b_1\alpha^*Y+\lambda^k\beta_1W+ h_1(Z,Y,W,\bar V),\\
\bar Y &= \hat\rho + \lambda^k\gamma^k Z  +\tilde d\gamma^kY^2+ \gamma^k h_2(Z,Y,W,\bar V),\\
\bar W &=  \lambda^k \alpha_3 Z+\lambda^k\beta_3W+ h_3(Z,Y,W, \bar V),\\
V&= \lambda^k \alpha_4 Z+\lambda^k\beta_4W+ h_4(Z,Y,W,\bar V),
\end{aligned}
\end{equation}  
where $(Z,Y,W,V)\in \mathbb{R}\times \mathbb{R}\times \mathbb{R}^{d^s-1}\times \mathbb{R}^{d^u-1}$, the coefficients $\alpha$ and $\beta$, as in the saddle-focus case \eqref{eq:maps:Tk_coef}, are linear combinations of the ones in \eqref{eq:df:Tk2}, 
\begin{equation*}\label{eq:df:Tk_const}
\begin{aligned}
\hat\rho= \lambda^k\gamma^k \rho -y^* +\lambda^k\gamma^k(\alpha^* x^+_1 + \beta^* x^+_2)+O(\hat\lambda^k\gamma^k ),
\end{aligned}
\end{equation*}
and the functions $h$ satisfy the same estimates as \eqref{eq:maps:Tk_nonlinear} except that $\partial h/\partial \bar V$ are replaced by \eqref{eq:h_df}, and that the derivatives with respect to the new variables $\rho$ are given by
\begin{equation*}
\begin{aligned}
\dfrac{\partial  h_{i}}{\partial \rho}
&=\dfrac{\partial  h_{i}}{\partial \mu}
\dfrac{\partial  \mu}{\partial \rho}
=O(\lambda^k\gamma^k\hat\gamma^{-k})+O(\lambda^kY^2),\quad i=1,3,4,\\
\dfrac{\partial  h_2}{\partial \rho}
&=\dfrac{\partial  h_2}{\partial \mu}
\dfrac{\partial  \mu}{\partial \rho}
=O(\lambda^k\gamma^k\hat\gamma^{-k})+o(\lambda^kY^2).
\end{aligned}
\end{equation*}



\subsection{Reduction to the saddle-focus case: proof of Theorem \ref{thm:sf} for the bi-focus case}
On the one hand, we note that proofs of all the results for the saddle-focus case, except for Lemmas \ref{lem:tangency1} and \ref{lem:tangency2} (which are used to prove Theorem \ref{thm:sfsimple}), are based on the analysis of the first-return map \eqref{eq:maps:Tk}, under the generic conditions C1--C5 and with  the weaker estimates \eqref{eq:h_df}. 
On the other hand,  the first-return map \eqref{eq:df:Tk} for the bi-focus case has the same form  as  \eqref{eq:maps:Tk}, and its estimates are the same as  \eqref{eq:maps:Tk_nonlinear} and \eqref{eq:h_df} (with rescaling $\mu$ as in \eqref{eq:rho}).  It follows that the above-mentioned results for the saddle-focus case also hold for the bi-focus case for all sufficiently large $k$ satisfying \eqref{eq:k}, if the corresponding $\rho$ values in those results are well-defined, that is, they are uniformly bounded for all $k$. This requirement on $\rho$ is easily justified by noting from \eqref{eq:mu_Q} that we will consider $\rho_k=\mu_k/\lambda^k=O(1)$.
 
 Thus, Theorem \ref{thm:sf} for the bi-focus case will be proven once we show that the transverse homoclinic points  given by Lemmas \ref{lem:tangency1}  and \ref{lem:tangency2} also exist in the bi-focus case. In fact, we now have a stronger result replacing these two lemmas.

\begin{lem}\label{lem:transhomo_df}
Let $\omega_1/\pi,\omega_2/\pi,1$ be rationally independent. If conditions C1--C4 are satisfied, then the tangency is accompanied.
\end{lem}

\begin{proof}
We follow the procedure in the proof of Lemma \ref{lem:tangency1}, and the Implicit Function Theorem will be used in the same way as before without further explanations. The focus is on finding and solving a system of equations whose solutions correspond to the desired homoclinic points, which is a counterpart to the one consisting of \eqref{eq:52} and \eqref{eq:51}.

Setting $l_1= \tilde y_1 - y^-_1,l_2=\tilde y_2 - y^-_2,l_3=\tilde v-v^-$, one finds from \eqref{eq:df:T1_cross} that the image $T_1(W^u_{loc}(O)\cap \Pi^-)$ is given by
\begin{equation}\label{eq:df:tran1}
\begin{aligned}
x_1 - x_1^+ &= b_{11} l_1 + O(|y_2|+|v|+l_1^2), \\
x_2 - x_2^+ &=  O(|v|+y^2_2+l_1^2),  \\
y_1 &= \mu +d l_1^2+ O(l_1^3) + O(|v|+|vl_1|+|y_2l_1|+y^2_2),  \\
u - u^+ &=b_{41}l_1+  O(|y_2|+|v|+l_1^2),  \\
l_2 &=b_{51} l_1+ O(|y_2|+|v|+l_1^2),\\
l_3 &=b_{61} l_1+ O(|y_2|+|v|+l_1^2).
\end{aligned}
\end{equation}
Expressing $l_1$ as a function of $x_1$ from the first equation, we obtain
\begin{equation}\label{eq:df:tran2}
 l_1 = b^{-1}_{11} (x_1 - x_1^+) + O(|y_2|+|v|+(x_1 - x_1^+)^2).
\end{equation}
Substituting the last two equations of \eqref{eq:df:T1_cross}, with setting $\tilde y=\tilde v=0$, into the $y$- and $v$-equations in \eqref{eq:df:T0k}, yields $y_2=O(\gamma^{-k})$ and $v=O(\hat\gamma^{-k})$ as  functions of $x_1,x_2,y_1,u$. Combining this with \eqref{eq:df:tran2} and the $y_1$- and $l_2$-equations in \eqref{eq:df:tran1}, yields 
\begin{equation}\label{eq:df:tran2.1}
\begin{aligned}
y_1&=\mu+\frac{d}{b_{11}^2}(x_1-x^+_1)^2+O((x_1-x^+_1)^3)+o(\gamma^{-k}),\\
 l_1 &= b^{-1}_{11} (x_1 - x_1^+) + O((x_1 - x_1^+)^2)+O(\gamma^{-k}),\\
 l_2 &=b^{-1}_{11}b_{51}(x_1 - x_1^+) + O((x_1 - x_1^+)^2)+O(\gamma^{-k}).
\end{aligned}
\end{equation}
Combining the first equation in \eqref{eq:df:tran2.1} with the $y_1$-equation in  \eqref{eq:df:T0k}, with applying  coordinate transformations \eqref{eq:df:coor1} and \eqref{eq:df:coor2}, yields
\begin{equation}\label{eq:df:tran3}
0=\mu-\gamma^{-k}y^*-\gamma^{-k}Y+\frac{d}{b_{11}^2}X^2+O(X^3)+o(\gamma^{-k}),
\end{equation}
where $y^*$ is given by \eqref{eq:df:alphabeta}. 

We now find the $y_1$-coordinates of points $(\bar x_1,\bar x_2,\bar y_1,\bar y_2,\bar u,\bar v)\in T_1\circ T^k_0\circ T_1(W^u_{loc}(O))$ by combining \eqref{eq:df:T0k} and \eqref{eq:df:T1_cross}, with using \eqref{eq:df:tran1} and \eqref{eq:df:tran2}. (Equivalently, one can just use the formula \eqref{eq:df:Tk2} and then undo the coordinate transformation $(\bar y,\bar v)\mapsto(\bar Y,\bar V)$). We obtain
\begin{equation*}
\bar y_1 =  \mu+\lambda^k \alpha^* x_1^++\lambda^k  \beta^* x_2^+ + \lambda^k \alpha^*X_1 
+ d Y^2 +o(\lambda^k)+ O(Y^3),
\end{equation*}
where the $o(\lambda^k)$ term is a function of $X_1,Y,\bar y_2,\bar v$. Such point belongs to $W^s_{loc}(O)$ if and only if $\bar y_1=\bar y_2=\bar v=0$, which implies
\begin{equation}\label{eq:df:tran4}
0 =  \mu+\lambda^k \alpha^* x_1^++\lambda^k  \beta^* x_2^+ + \lambda^k \alpha^*X_1 
+ d Y^2 +o(\lambda^k)+ O(Y^3),
\end{equation}
where  $o(\lambda^k)$ now only depends on $X_1,Y$.

Transverse homoclinic points correspond to non-degenerate solutions to the system consisting of \eqref{eq:df:tran3} and \eqref{eq:df:tran4}. Comparing these two equations with \eqref{eq:52} and \eqref{eq:51}, and recalling the discussion below \eqref{eq:q}, one sees that the remaining computation can be carried out exactly as in the saddle-focus case if $k$ can be taken such that \eqref{eq:k} is satisfied,  $dy^*>0$, and $d(\alpha^*x_1^++\beta^*x_2^+)<0$. The existence of such $k$ is obvious since the set $\{(k\omega_1/2\pi,k\omega_2/2\pi)\}$ is dense in $\mathbb{R}^2$, due to the rational independence among $\omega_1/\pi,\omega_2/\pi,1$.
\end{proof}



\section{Acknowledgments}
We thank Pablo Barrientos, Sergey Gonchenko and Dmitrii Mints for useful discussions. We also thank Shuhei Hayashi and Shuntaro Tomizawa for their comments after reading the first version of the manuscript.
Dongchen Li was supported by the Science Fund Program for Excellent Young Scientists (Overseas), the New Cornerstone Science Foundation, and the ERC project 677793 StableChaoticPlanetM. Xiaolong Li was supported by the NSFC Grants 11701199 and 12331005. Katsutoshi Shinohara was supported by the JSPS KAKENHI Grants 21K03320 and 24K06793. Dmitry Turaev was supported by the Leverhulme Trust.

%

\end{document}